\documentclass[11pt,
]{article}

\usepackage{amsmath,amsthm,amssymb,%amstext,%amsfonts,
latexsym,mathrsfs} 
\usepackage[a]{esvect}
\usepackage[bbgreekl]{mathbbol}

\usepackage{mparhack}
\usepackage{xcolor}

\usepackage[english,greek]{babel}
\languageattribute{greek}{polutoniko}
\usepackage{marginnote}

\usepackage{makeidx}
\makeindex

\usepackage
[%hypertex%,
colorlinks,linkcolor=magenta,citecolor=blue,
]
{hyperref}
\newcommand{\nek}{\newcommand}
\nek{\renek}{\renewcommand}
\nek{\vyk} [1] {}

\makeatletter
\if@twoside%
\else\fi%
\makeatother  
\DeclareFontFamily{U}{sdd}{\skewchar\font127 }
\DeclareFontShape{U}{sdd}{m}{n}{ <->[0.9]bbold10}{}%
\DeclareMathAlphabet{\nmi}{U}{sdd}{m}{n}%

\nek{\ubf}{\fontseries{b}\selectfont}
\nek{\bfit}{\bfseries\itshape}
\nek{\bftt}{\ttfamily\bfseries\upshape\selectfont}

\makeatletter
\@definecounter{section}
\makeatother
\renek{\thesection}{\Roman{section}}
\nek{\gla}{\clearpage\section}

\nek{\parf}{\subsection}
\renek{\thesubsection}{\arabic{subsection}}
\nek{\punk}{\subsubsection}
\renek{\thesubsubsection}
{\thesubsection\alph{subsubsection}}

\theoremstyle{plain}

\newtheorem{theore}             {Theorem} [subsection]
\newtheorem{corollary}  [theore]{Corollary}
\newtheorem{prop}  [theore]{Proposition}
\newtheorem{lemma}      [theore]{Lemma}
\newtheorem{claim}      [theore]{Claim}

\theoremstyle{definition}

\newtheorem{definition} [theore]{Definition}
\newtheorem{vopr}       {Question}
\newtheorem{prim}       [theore]{Example}
\newtheorem{zam}        [theore]{Remark}
\newtheorem*{prF}   {Proof}    
\newtheorem*{ack}   {Acknowledgement}    
\newtheorem{ggi}   [theore] {Blanket Assumption}  

\nek{\thsp}{\hspace{0.1ex}}
\nek{\bvo}{\begin{vopr}}
\nek{\evo}{\end{vopr}}
\nek{\back}{\begin{ack}}
\nek{\eack}{\end{ack}}
\nek{\bpro}{\begin{prop}}
\nek{\epro}{\end{prop}}
\nek{\bcor}{\begin{corollary}}
\nek{\ecor}{\end{corollary}}
\nek{\bdf} {\begin{definition}}
\nek{\edf} {\qed\end{definition}}
\nek{\eDf} {\end{definition}}
\nek{\bgg} {\begin{ggi}}
\nek{\egg} {\qed\end{ggi}}
\nek{\ble} {\begin{lemma}}
\nek{\ele} {\end{lemma}}
\nek{\bcl} {\begin{claim}}
\nek{\ecl} {\end{claim}}
\nek{\bpri}{\begin{prim}}
\nek{\epri}{\qed\end{prim}}
\nek{\eprI}{\end{prim}}
\nek{\bte} {\begin{theore}}
\nek{\ete} {\end{theore}}
\nek{\baq} {\begin{aq}}
\nek{\eaq} {\end{aq}}
\nek{\bre}{\begin{zam}}
\nek{\ere}{\qed\end{zam}}
\nek{\eRe}{\end{zam}}
\nek{\bpf} {\begin{prF}} 
\nek{\epf} {\qed\end{prF}} 
\nek{\epg} {\end{prF}} 
\nek{\epF}[1]{\hfill\hbox{$\square$\ ({\small#1})}\end{prF}}

\nek{\qeD}[1]{\hfill\hbox{$\square$\ ({\small#1})}}

\nek{\ben}{\begin{enumerate}\itemsep=0.2em}
\nek{\een}{\end{enumerate}}
\nek{\bde}{\begin{description}\itemsep=0.2em}
\nek{\ede}{\end{description}}
\nek{\bit}{\begin{itemize}\itemsep=0.2em}
\nek{\eit}{\end{itemize}}
\nek{\bay}{\begin{array}}
\nek{\eay}{\end{array}}
\nek{\bce}{\begin{center}}
\nek{\ece}{\end{center}}
\nek{\bur}{\begin{equation}}
\nek{\eur}{\end{equation}}
\nek{\bqu}{\begin{quotation}\noindent}
\nek{\equ}{\end{quotation}}
\nek{\tre} [1] {\mathop{\tt tree}(#1)}
\nek{\baz} [1] {\mathop{\tt base}(#1)}
\nek{\lh} [1] {\mathop{\tt lh}(#1)}
\nek{\siz} [1] {\mathop{\tt size}(#1)}
\nek{\card} {\mathop{\tt card}}
\nek{\dom} {\mathop{\tt dom}}
\nek{\clos} {\mathop{\tt clos}}
\nek{\num} {\mathop{\tt num}}
\nek{\ran} {\mathop{\tt ran}}
\nek{\rank} {\mathop{\tt rank}}
\nek{\tlim} {\mathop{\tt lim}}
\nek{\tsup} {\mathop{\tt sup}}
\nek{\tmax} {\mathop{\tt max}}
\nek{\sep}{{\tt Sep}}
\nek{\red}{{\tt Red}}
\nek{\TC}  {{\rm TC}\hspace{0.4ex}}
\nek{\hc}  {\mathrm{HC}}
\nek{\rL} {{\mathbf L}}
\nek{\rV} {{\mathbf V}}

\nek{\ZFC} {\text{\ubf ZFC}}
\nek{\ZFL} {\text{\ubf ZFL}}
\nek{\zfc} {\ZFC}
\nek{\zfcm} {\zfc^-}
\nek{\al}  {\alpha}
\nek{\ga}  {\gamma}
\nek{\Ga}  {\Gamma}
\nek{\da}  {\delta}
\nek{\Da}  {\Delta}
\nek{\kpa} {\kappa}
\nek{\la}  {\lambda}
\nek{\ve}  {\varepsilon}
\nek{\vpi} {\varphi}
\nek{\vpo} {\overline\vpi}
\nek{\sg}  {\sigma}
\nek{\Sg}  {\Sigma}
\nek{\om}  {\omega}
\nek{\Om}  {\Omega}
\nek{\lom} {^{<\om}} 
\nek{\omi} {\om_1}
\nek{\oli} {\omi^\rL}
\nek{\omb} {\om_2}
\nek{\omil}{\om_1^\rL}
\nek{\za}  {\zeta}
\nek{\ali} {\aleph_1}
\nek{\ald} {\aleph_2}
\nek{\alil} {\aleph_1^\rL}
\nek{\fs}[2]{{\mathbf\Sigma}^{#1}_{#2}}
\nek{\fp}[2]{{\mathbf\Pi}^{#1}_{#2}}
\nek{\fd}[2]{{\mathbf\Delta}^{#1}_{#2}}
\nek{\is}[2]{\varSigma^{#1}_{#2}}
\nek{\ip}[2]{\varPi^{#1}_{#2}}
\nek{\isp}[2]{(\varSigma+\varPi)^{#1}_{#2}}
\nek{\id}[2]{\varDelta^{#1}_{#2}}
\nek{\lsp}[2]
{\cL\hspace*{-0.4ex}(\varSigma{+\hspace*{-0.15ex}}\varPi)^{#1}_{#2}}
\nek{\ls}[2]{\cL\hspace*{-0.4ex}\varSigma^{#1}_{#2}}
\nek{\lp}[2]{\cL\hspace*{-0.4ex}\varPi^{#1}_{#2}}
\nek{\ld}[2]{\cL\hspace*{-0.4ex}\varDelta^{#1}_{#2}}
\nek{\lsh}[1]{\cL\varSigma^{\hc}_{#1}}
\nek{\lph}[1]{\cL\varPi^{\hc}_{#1}}
\nek{\ldh}[1]{\cL\varDelta^{\hc}_{#1}}
\nek{\BBB}{} %{\hspace{0.0ex}}
\nek{\dA}{{\BBB{\mathbb A}\BBB}}
\nek{\dB}{{\BBB{\mathbb B}\BBB}}
\nek{\dC}{{\BBB{\mathbb C}\BBB}}
\nek{\dD}{{\BBB{\mathbb D}\BBB}}
\nek{\dJ}{{\BBB{\mathbb J}\BBB}}
\nek{\dP}{{\BBB{\mathbb P}\BBB}}
\nek{\dQ}{{\BBB{\mathbb Q}\BBB}}
\nek{\dR}{{\BBB{\mathbb R}\BBB}}
\nek{\gM}{{\BBB{\goth M}\BBB}}
\nek{\gN}{{\BBB{\goth N}\BBB}}
\nek{\gK}{{\BBB{\goth K}\BBB}}
\nek{\cM}{{\BBB{\skri M}\BBB}}
\nek{\cP}{{\BBB{\skri P}\BBB}}
\nek{\cS}{\skri S}
\nek{\pu}  {\varnothing}
\nek{\sq}  {\subseteq}
\nek{\qs}  {\supseteq}
\nek{\su}  {\subset}
\nek{\sneq}{\subsetneqq}
\nek{\eqv} {\mathbin{\,\Longleftrightarrow\,}}
\nek{\imp} {\mathbin{\,\Longrightarrow\,}}
\nek{\mpi} {\mathbin{\,\Longleftarrow\,}}
\nek{\lra} {\longrightarrow} 
\nek{\we}  {{\mathbin{\hspace*{0.2ex}^\smallfrown}}}
\nek{\sus} {{\exists\,}}
\nek{\kaz} {{\forall\,}}
\nek{\ti}  {\times}
\nek{\dm}  {$$}
\nek{\iy}  {\infty}
\nek{\abs} [1] {|#1|}
\nek{\abt} [1]
{\mathopen{{|}\hspace*{-0.2ex}{|}}#1\mathclose{{|}\hspace*{-0.2ex}{|}}}

\nek{\nin}{\notin}
\nek{\res} {{\hspace*{0.2ex}{\restriction}\hspace*{0.2ex}}}

\nek{\bang} [1] {\big\langle #1\big\rangle}
\nek{\ang} [1] {\langle #1\rangle}
\nek{\ans} [1] {\{\hspace{0.1ex}#1\hspace{0.1ex}\}}
\nek{\dd}[1]{$\mtho\hspace{0.0ex}{#1}$-\hspace{0.0ex}}
\nek{\itla}{\item\label}
\nek{\itlb} [1] {\itla{#1}\imas{#1}}
\nek{\itlc} [1] {\itla{#1}\imat{#1}}

\nek{\bez}{\smallsetminus}
\nek{\rit} [1] {{\it#1\/}}

\nek{\ens} [2] {\ans{{#1\hspace{0.5ex}{:}}\zz\hspace{0.5ex}#2}}
\nek{\zz} {\linebreak[0]} 
\nek{\goth}{\mathfrak}
\nek{\skri}{\mathscr}

\nek{\yo} {,\linebreak[0]}
\nek{\yi} {\hspace{\mathsurround},\linebreak[0]\hspace{\mathsurround}}
\nek{\yd} {\hspace{\mathsurround},\linebreak[0]\:}
\nek{\yt} {\hspace*{\mathsurround}\text{,}\linebreak[0]\;}

\nek{\viv} {\text{vice versa}}
\nek{\ie} {\text{\sl i.\resp e.}}
\nek{\eg} {\text{\sl e.\resp g.}}
\nek{\poo} {\text{w.\resp r.\resp t.}}
\nek{\noo} {\text{w.\resp l.\resp o.\resp g.}}

\nek{\iesp}{\hspace{0.3ex}}
\nek{\resp}{\hspace{0.25ex}}

\nek{\itsep}{\itemsep=0.25ex plus 0.1ex minus 0.1ex}

\nek{\tenu}[1]{

\itsep}

\nek{\tenv}[1]{
\def\theenumii{#1}
\itsep}

\nek{\cenu}{ 
}

\nek{\Aenu}{\tenu{(\Alph{enumi})}}
\nek{\aenu}{\tenu{(\alph{enumi})}}
\nek{\aenr}{\tenu{{\rm(\alph{enumi})}}}
\nek{\aenup}{\tenu{(\alph{enumi}$'$)}}
\nek{\nenu}{\tenu{{\rm(\arabic{enumi})}}}
\nek{\nenui}{\tenv{{\rm(\arabic{enumii})}}}

\nek{\aenri}{\tenv{{\rm(\alph{enumii})}}}

\nek{\fenu}{\itsep\tenu{{\mtho$(\fnsymbol{enumi})$}}}
\nek{\fenv}{\tenv{$(\fnsymbol{enumii})$}}
\nek{\renu}{\itsep\tenu{{\rm(\roman{enumi})}}\itsep}
\nek{\Renu}{\itsep\tenu{{\rm(\Roman{enumi})}}\itsep}

\nek{\atc} {\addtocounter{enumi}1}
\nek{\atm} {\addtocounter{enumi}{-1}}

\nek{\stk} [2] {\ang{#1\hspace{0.3ex};\hspace{0.1ex}#2}}
\nek{\sis} [2] {\ang{#1}_{#2}}
\nek{\sid} [3] {\bang{#1}{}^{#2}_{#3}}

\renek{\cM} {\goth M}
\nek{\cN} {\goth N}

\nek{\lam} [1]
{\label{#1}\hspace*{-3pt}\imaq{#1}%
}%
\nek{\las} [1]
{\label{#1}\imae{#1}}%

\nek{\jmar}[1]{\marginpar[%\vspace{-1ex}%
\flushright\footnotesize%
$\mtho\longrightarrow$\\%
\vspace{-1ex}{#1}\vspace*{1ex}]%
{%\vspace{-1ex}%
\flushleft\footnotesize%
$\mtho\longleftarrow$\\%
\vspace{-1ex}{#1}\vspace*{1ex}}%
}%
\nek{\pws} [1] {\cP(#1)}
\nek{\lla} {\,\land\,}

\nek{\snos} [1] {\,\footnote{\ #1}}
\nek{\snom}   {\,\footnotemark}
\nek{\snot} [1] {\footnotetext{\ #1}}

\nek{\ap} {{\hspace*{0.2ex}\boldsymbol\cdot\hspace*{0.2ex}}} 

\nek{\onto} {\overset{{\text{onto}}}{\longrightarrow}}

\nek{\La} {\Lambda}
\nek{\alo} {{\aleph_0}}

\nek{\leqv} {\,\eqv\,}

\nek{\dop} [1] {{#1}{}^{\complement}}

\nek{\ba}{\beta}

\nek{\qand} {\quad\text{and}\quad}

\nek{\vt}{\vartheta}

\nek{\nse} {\om^{<\om}}
\nek{\bse} {2^{<\om}}

\nek{\namx} [1] 
{\overset{\text{\mtho$\hspace*{0.2ex}_\text{\Large\bf.}$}}{#1}}
\nek{\namy} [1] 
{\overset{\text{\mtho$\hspace*{0.5ex}_\text{\Large\bf.}$}}{#1}}
\nek{\namz} [1] 
{\overset{\text{\mtho$\hspace*{0.9ex}_\text{\Large\bf.}$}}{#1}}

\nek{\cD} {\mathscr D}

\nek{\np}{\newpage}

\nek{\sqf} {\sq^{\text{\tt fin}}}
\nek{\sqfd} {\sq^{\text{\tt fd}}}

\nek{\ka}{\kappa}

\nek{\pwor} [1] {{\text{\rm PWO}}(#1)}
\nek{\sepa} [1] {{\text{\ubf Sep}}(#1)}
\nek{\redu} [1] {{\text{\ubf Red}}(#1)}

\nek{\pet} {\text{\ubf PT}}
\nek{\ret} [1] {\res_{\hspace*{-0.05ex}#1}}

\renek{\leq} {\leqslant}
\renek{\geq} {\geqslant}

\nek{\limp} {\mathrel{\,\imp\,}}

\nek{\spe} [1] {\text{\bf SC}(#1)}
\nek{\spf} [1] {\text{\bf SC}(#1)}
\nek{\spg} [1] {\text{\bf SC}^{<\om}(#1)}

\nek{\vys} [1] {\text{\tt hgt}(#1)} 

\nek{\cle} {\preccurlyeq}
\nek{\cls} {\prec}

\renek{\refname} {{\large\bf References}}

\nek{\dphi} {{\mathbb\Phi}}

\nek{\dU}{\mathbb U}
\nek{\uu}{^{\text{\tt fu}}}

\nek{\Tv}{\overrightarrow{T}}
\nek{\Sv}{\overrightarrow{S}}

\nek{\jLa} {\boldsymbol\La}
\nek{\jhi} {\boldsymbol\chi}
\nek{\jfi} {\boldsymbol\vpi}
\nek{\jsi} {\boldsymbol\psi}
\nek{\jta} {\boldsymbol\tau}
\nek{\jpi} {\boldsymbol\pi}
\nek{\jsg} {\boldsymbol\sg}
\nek{\jrho} {\boldsymbol\rho}
\nek{\jro} {\jrho}
\nek{\ju} {\boldsymbol u}
\nek{\jv} {\boldsymbol v}
\nek{\jw} {\boldsymbol w}

\nek{\dn} {2^\om}
\nek{\bn} {\om^\om}

\nek{\dox}{{\namy{\boldsymbol x}}}
\nek{\doy}{{\namy{\boldsymbol y}}}
\nek{\doa}{{\namy {a}}}
\nek{\dob}{{\namy {b}}}
\nek{\doH}{{\namy {H}}}

\nek{\rc} {\mathbf c}
\nek{\rd} {\mathbf d}
\nek{\rpi} {\dox}

\nek{\kn} [1] {\underline{#1}}

\nek{\dplom} {\dd{\dP\lom}}
\nek{\dplo} {\dd\plo}

\nek{\plo} {\dP\lom}
\nek{\uplo} {(\dP\cup\dU)\lom}

\nek{\cll} {\mathrel{{\cls}\hspace*{-0.9ex}{\cls}}}

\nek{\vpj} [2] {\vpi^{#1}_{#2}}
\nek{\Phj} [2] {\Phi^{#1}_{#2}}

\nek{\vpt} [3] {\vpi^{#1}_{#2#3}}

\nek{\od} {\text{\rm OD}}

\nek{\Eo} {\mathrel{{\text{\sf E}}_0}}

\nek{\zf} {\text{\ubf ZF}}

\nek{\uf} [2] {{\boldsymbol Q}^\dphi_{#1}(#2)}
\nek{\ufi} [1] {{\boldsymbol Q}^\dphi_{#1}}

\nek{\qf} [4] {{\boldsymbol Q}^\dphi_{#1#2,#3}(#4)}

\nek{\tf} [3] {{\boldsymbol Q}^\dphi_{#1#2}(#3)}
\nek{\tfi} [2] {{\boldsymbol Q}^\dphi_{#1#2}}

\nek{\tx} [2] {{\boldsymbol T}^\dphi_{#1}(#2)}
\nek{\txi} [1] {{\boldsymbol T}^\dphi_{#1}}

\nek{\qq} [3] {{\boldsymbol T}^\dphi_{#1#2,#3}}

\nek{\ty} [3] {{\boldsymbol T}^\dphi_{#1#2}(#3)}
\nek{\qy} [4] {{\boldsymbol T}^\dphi_{#1#2,#3}(#4)}

\nek{\lel} {\leq_{\rL}}

\nek{\kc} [2] {C_{#1#2}}
\nek{\xc} [2] {C'_{#1#2}}
\nek{\kcc} [2] {C'_{#1#2}}
\nek{\kd} [2] {D_{#1#2}}
\nek{\kcp} [3] {C_{#1#2}^{#3}}
\nek{\kr} [3] {\jrho_{#1#2}^{#3}}
\nek{\kR} [2] {R_{#1#2}}

\nek{\jcp} [4] {C_{#1#2}^{#3#4}}
\nek{\jr} [4] {\jrho_{#1#2}^{#3#4}}

\nek{\roo} [1] {\text{\tt stem}(#1)}

\nek{\cI} {\mathcal I}
\nek{\fso} {\text{\bf FSS}}
\nek{\fss} [1] {\text{\bf FSS}(#1)}
\nek{\fsd} [1] {\text{\bf FSS}\lom(#1)}

\nek{\ptf} {\text{\ubf {PTF}}}

\nek{\jp} {{\mathbb p}}
\nek{\bbu} {\mathbb u}
\nek{\bbp} {{\jta}}
\nek{\bbt} {{\boldsymbol\tau}}

\nek{\bbq} {\mathbb q}
\nek{\bbpu} {\bbp{\lor}\bbu}
\nek{\bbpr} {\jpi{\lor}\jqo}
\nek{\jpq} {\jpi\kw\jqo}

\nek{\mus} {multi\-sys\-tem}
\nek{\Mus} {Multisystem}
\nek{\muq} {multisequence}
\nek{\Muq} {Multisequence}
\nek{\ms} [1] {\text{\ubf MS}(#1)}

\nek{\Muf} {Multiforcing}
\nek{\muf} {multiforcing}
\nek{\Mut} {Multitree}
\nek{\mut} {multitree}
\nek{\aac} {{\sc aac}}
\nek{\maa} {multi-{\sc aac}}
\nek{\maac} {multi almost antichain}
\nek{\md} {\text{\ubf MT}}
\nek{\mt} [1] {\text{\ubf MT}(#1)}
\nek{\mtr} [1] {\text{\ubf MT}_{#1}}
\nek{\mfr} [1] {\text\mf_{#1}}
\nek{\mf} {\text{\ubf sMF}}
\nek{\omim} {\omi^\cM}
\nek{\et}{\eta}             

\nek{\bas} {\underline s}
\nek{\bat} {\underline t}
\nek{\bah} {H}
\nek{\baH} {2^H}

\nek{\gh} {2^{H}} %{\mathfrak H}

\nek{\jvt} {\bbtheta}

\nek{\ivt} {\underline\vt^{-1}}

\nek{\wo} {\text{\ubf WO}}
\nek{\mto} {\mapsto}

\nek{\ahh} [1] {A_{#1}^\ast}
\nek{\bhh} [1] {B_{#1}^\ast}

\nek{\uH}{\mathbf h}

\nek{\bp} {\bar p}
\nek{\bq} {\bar q}

\nek{\mtp} [1] {\mt{#1}^+}

\nek{\Ord} {\text{\ubf Ord}}

\nek{\xh} {N}

\nek{\gfu} {\mathbb W}

\nek{\mor} [1] {\rL[G\res\Da_{#1}]}
\nek{\morp} [1] {\rL[G\res\Da'_{#1}]}

\nek{\tup}{\textup}

\nek{\mtho}{\mathsurround=0mm}
\nek{\msur}{\hspace*{-1\mathsurround}}
\nek{\dsur}{\hspace{-0.3\mathsurround}}
\nek{\hsur}{\hspace{-0.5\mathsurround}}
\nek{\noi}{\noindent}
\nek{\vim}{\vspace{-0.5mm}}
\nek{\vom}{\vspace{1mm}}
\nek{\vtm}{\vspace{2mm}}

\nek{\fP} {\pmb{\mathbb P}}
\nek{\fQ} {\pmb{\mathbb Q}}
\nek{\tSg} [1] {\widetilde\varSigma^1_{#1}}
\nek{\tPi} [1] {\widetilde\varPi^1_{#1}}
\nek{\tSP} {(\widetilde\varSigma+\widetilde\varPi)^1_1}

\nek{\bcu} {q_1}

\nek{\raw} [2] {#1{(\to#2)}}
\nek{\req} [2] {{#1}\ret{#2}}
\nek{\pes} {\pet}
\nek{\nq} [1] {\mathrel{{\sq}_{#1}}}
\nek{\leqs} {\leqslant}
\nek{\mtt} [2] {\tau^{#1}_{#2}}
\nek{\mth} [2] {h^{#1}_{#2}}
\nek{\ntt} [3] {T^{#1}_{#2}(#3)}
\nek{\snenu}{\tenu{{\rm(\thesubsection.\arabic{enumi})}}}
\nek{\sif} {{\ubf PTF}}
\nek{\sct} [2] {\text{\ubf CT}_{#2}(#1)}
\nek{\sco} [1] {\text{\ubf CT}(#1)}

\nek{\TS}{\textstyle}

\nek{\qh} [4] {h^{#1}(#2,#3,#4)}
\nek{\qT} [5] {T^{#1}_{#2#3,#4}(#5)}

\nek{\zh} [3] {h^{#1}(#2,#3)}
\nek{\zt} [3] {\tau^{#1}_{#2#3}}
\nek{\zT} [4] {T^{#1}_{#2#3}(#4)}
\nek{\zd} [3] {T^{#1}_{#2#3}}

\nek{\zpi} {\boldsymbol\pi}
\nek{\zro} {\boldsymbol\rho}
\nek{\zsg} {\boldsymbol\sigma}
\nek{\zp} {{\boldsymbol p}}
\nek{\zq} {{\boldsymbol q}}
\nek{\zr} {{\boldsymbol r}}

\nek{\ta} {\tau}

\nek{\Ta} {T^\ast}

\nek{\ah} {\ibb} %[1] {\mathbb0_{#1}}
\nek{\bh} {\dbb} % [1] {\mathbb1_{#1}}
\nek{\ch} {\idbb} %[1] {\mathbb2_{#1}}

\nek{\kk} [3] {K^{#3}_{#1#2}}
\nek{\kkc} [2] {\kk{#1}{#2}\rc}

\nek{\ki} [2] {\kk{}{#2}{#1}}
\nek{\kic} [1] {\ki\rc{#1}}

\nek{\kip} [3] 
{\ki{#1}{#2}{\hspace*{0.2ex}\uparrow\hspace*{0.2ex}}{#3}}
\nek{\kipc} [2] {\kip\rc{#1}{#2}}

\nek{\kkcp} [3] {\kkc{#1}{#2}{\uparrow}{#3}}

\nek{\bg} {\mathbf g}

\nek{\etz}{{\boldsymbol\za}}

\nek{\dpo} {\dP_{\text{\rm coh}}}

\nek{\rE} {\mathrel{\mathsf E}}

\nek{\hgh} [2] {\hc(#1,#2)}

\nek{\ibb} [1] {\mathbb1_{#1}}
\nek{\dbb} [1] {\mathbb2_{#1}}
\nek{\idbb} [1] {\mathbb{12}_{#1}}
\nek{\obb} [1] {\mathbb{0}_{#1}}

\nek{\idt} {\ans{1,2,12}}
\nek{\idq} {\ans{0,1,2,12}}

\nek{\xxh} [2] {X_{#1}(#2)}
\nek{\yyh} [2] {Y_{#1}(#2)}

\nek{\dxh} [1] {X_{#1}(\doH)}
\nek{\dyh} [1] {Y_{#1}(\doH)}

\nek{\uG} {{\underline G}}
\nek{\unH} {{\underline H}}
\nek{\ua} {{\underline a}}
\nek{\ux} {{\underline x}}
\nek{\uy} {{\underline y}}

\nek{\jqo} {\text{\bfit\textgreek{\qoppa}}}
\nek{\jQo} {\text{\bfit\textgreek{\Qoppa}}}

\nek{\qfd} [4] {{\boldsymbol Q}^\dphi_{#1#2,#3}(#4)}
\nek{\tfd} [4] {{\boldsymbol T}^\dphi_{#1#2,#3}(#4)}

\nek{\qfc} [3] {{\boldsymbol Q}^\dphi_{#1#2,#3}}
\nek{\tfc} [3] {{\boldsymbol T}^\dphi_{#1#2,#3}}

\nek{\dqf} [2] {{\dQ}^\dphi_{#1#2}}

\nek{\dpi} [2] {\dpj{#1}{#2}\jpi}
\nek{\dqo} [2] {\dpj{#1}{#2}\jqo} 
\nek{\dpq} [2] {\dpj{#1}{#2}{\jpq}}
\nek{\dro} [2] {\dpj{#1}{#2}\jro} 
\nek{\dsg} [2] {\dpj{#1}{#2}\jsg} 

\nek{\dpj} [3] {{#3}(#1,#2)}
\nek{\zdp} [3] {{#1}(#2,#3)}

\nek{\jf} {{\boldsymbol\mu}}
\nek{\bssq} {\mathrel{\boldsymbol\sqsubset}}

\nek{\ssq} {\sqsubset}
\nek{\ssm} [1] {\mathrel{{\bssq}\hspace*{-1.6ex}{\bssq}}_{#1}}

\nek{\ssa} [1] {\ssq_{#1}}
\nek{\ssb} [1] {\bssq_{#1}}
\nek{\ssc} [1] {\bssq^3_{#1}}
\nek{\sse} [3] {\bssq^{#1}_{#2#3}}

\nek{\gom} [1] {\boldsymbol h_{#1}}
\nek{\gon} [2] {\boldsymbol h_{#1#2}}
\nek{\obr} {^{-1}}

\nek{\pil} [1] {\jpi[#1]}
\nek{\tal} [1] {\jta[#1]}
\nek{\zal} [1] {\bz[#1]}
\nek{\bz} {\mathbf z}
\nek{\ja} {\mathbf a}
\nek{\jb} {\mathbf b}

\nek{\mmf} {\text{\ubf MMF}}

\nek{\mmja} [1] {\cM[#1]}
\nek{\mmpi} [1] {\jpi[#1]}

\nek{\abd} [2] {|#1|^{#2}}
\nek{\pro} [2] {(#1)_{#2}}

\nek{\ttt} {\mathbf s}

\nek{\ddi} [4] {\zD^{#4}_{#1#2}(#3)}

\nek{\mes} [1] 
{{}\mathord{{\res}\hspace*{-0.42em}{\res}}_{\hspace*{-0.2ex}#1}}

\nek{\mfl} [2] {\mf_{#1}\mes{{<}\hspace*{0.2ex}#2}}
\nek{\mfg} [2] {\mf_{#1}\mes{{\ge}\hspace*{0.2ex}#2}}
\nek{\mfe} [2] {\mf_{#1}\mes{\hspace*{0.2ex}#2}}

\nek{\mtl} [1] {\md{}\mes{{<}\hspace*{0.2ex}#1}}
\nek{\mtg} [1] {\md{}\mes{{\ge}\hspace*{0.2ex}#1}}
\nek{\mte} [1] {\md{}\mes{\hspace*{0.2ex}#1}}

\nek{\ntl} [2] {\mt{#1}\mes{{<}\hspace*{0.2ex}#2}}
\nek{\ntg} [2] {\mt{#1}\mes{{\ge}\hspace*{0.2ex}#2}}
\nek{\nte} [2] {\mt{#1}\mes{\hspace*{0.2ex}#2}}

\nek{\nfl} [1] {\mf\mes{{<}\hspace*{0.2ex}#1}}
\nek{\nfg} [1] {\mf\mes{{\ge}\hspace*{0.2ex}#1}}
\nek{\nfe} [1] {\mf\mes{\hspace*{0.2ex}#1}}

\nek{\npl} [1] {\mfp\mes{{<}\hspace*{0.2ex}#1}}
\nek{\npg} [1] {\mfp\mes{{\ge}\hspace*{0.2ex}#1}}
\nek{\npe} [1] {\mfp\mes{\hspace*{0.2ex}#1}}

\nek{\vmf}{\vv\mf}
\nek{\vmi}{\vmf_{\omi}}
\nek{\nvl} [1] {\vmf\mes{{<}\hspace*{0.2ex}#1}}
\nek{\nvg} [1] {\vmf\mes{{\ge}\hspace*{0.2ex}#1}}
\nek{\nve} [1] {\vmf\mes{\hspace*{0.2ex}#1}}
\nek{\ivl} [1] {\vmi\mes{{<}\hspace*{0.2ex}#1}}
\nek{\ivg} [1] {\vmi\mes{{\ge}\hspace*{0.2ex}#1}}
\nek{\ive} [1] {\vmi\mes{\hspace*{0.2ex}#1}}

\nek{\kvl} [2] {\vv{\mf}(#2)\mes{{<}\hspace*{0.2ex}#1}}
\nek{\kvg} [2] {\vv{\mf}(#2)\mes{{\ge}\hspace*{0.2ex}#1}}
\nek{\kve} [2] {\vv{\mf}(#2)\mes{\hspace*{0.2ex}#1}}

\nek{\prl} [2] {{#1}\mes{{<}\hspace*{0.2ex}#2}}
\nek{\prg} [2] {{#1}\mes{{\ge}\hspace*{0.2ex}#2}}
\nek{\pre} [2] {{#1}\mes{\hspace*{0.2ex}#2}}

\nek{\sbr} [1] {\lseil#1\rseil}

\nek{\vjpi} {\vv\jpi}
\nek{\vjqo} {\vv\jqo}
\nek{\vjro} {\vv\jro}
\nek{\vjsg} {\vv\jsg}

\nek{\nos} [2] {{#1}\lceil{#2}\rceil}
\nek{\nop} [2] {{#1}_{\hspace*{-0.3ex}+}\lceil{#2}\rceil}
\nek{\nom} [3] {{#1}\lceil{#3,#2}\rceil}

\nek{\nor} [2] {{#1}_{#2}}
\nek{\norl} [2] {{#1}_{<#2}}

\nek{\duz} [3] {{#1}^{\abc{#2}}_{#3}}
\nek{\dqu} [3] {\zD(#1,#2,#3)}

\nek{\ufo} {\text{\ubf un}}
\nek{\ups} {\Upsilon}
\nek{\jPi} {\nmi\Pi}

\nek{\vjPi} {\vv{\jPi}}
\nek{\mup} [2] 
{\mathord{\lceil#1,#2\rceil}\hspace*{-0.1ex}^+}

\nek{\len} [1] {\dom(#1)}

\nek{\dma} [2] {\dD_{#1}^{#2}}

\nek{\mfp} {\text{\ubf MF}}

\nek{\jpn}  [1] {\jPi_{#1}}           
\nek{\pilg} [1] {\jPi_{<#1}}
\nek{\pigg} [1] {\jPi_{\ge#1}}

\nek{\dcp} {\dC'}

\nek{\fpl}  [1] {\fP_{<#1}}
\nek{\fpg}  [1] {\fP_{\ge#1}}

\nek{\reb} [1] {\mes{\hspace*{0.2ex}#1}}
\nek{\abc} [1] {|#1|}

\nek{\mdi} {special}
\nek{\mre} {regular}
\nek{\evd} {absolutely incompatible}
\nek{\evdy} {absolutely incompatible}
\nek{\emd} {\evd}

\nek{\cwun} {\bigcup^{\text{\tt cw}}}
\nek{\kw} {\cup^{\text{\tt cw}}}
\nek{\bkw} {\bigcup^{\text{\tt cw}}}

\nek{\zD} {{\boldsymbol D}}
\nek{\cmp} {\cM^+}

\nek{\bo}{\mathbf 0}
\nek{\sko} [1] {\textup{(}#1\/\textup{)}}

\nek{\cL}{\mathscr L}

\nek{\xm} {\widetilde m}
\nek{\xk} {\widetilde k}

\renek{\cM}{M}

\nek{\sfo} {{\tt sforc}}
\nek{\sfoe} [1] 
{\mathop{\hspace*{0.4ex}\sfo_{#1}\hspace*{0.4ex}}}
\nek{\sfof}   
{\mathop{\hspace*{0.4ex}\sfo\hspace*{0.4ex}}}

\nek{\wfo} {{\tt wforc}}
\nek{\wfoe} [1] 
{\mathop{\hspace*{0.4ex}\wfo_{#1}\hspace*{0.4ex}}}
\nek{\wfof}   
{\mathop{\hspace*{0.4ex}\wfo\hspace*{0.4ex}}}

\nek{\fo} {{\tt forc}}
\nek{\fod} [2] 
{\mathop{\hspace*{0.4ex}\fo_{#1#2}\hspace*{0.4ex}}}
\nek{\foe} [1] 
{\mathop{\hspace*{0.4ex}\fo_{#1}\hspace*{0.4ex}}}
\nek{\fof}   
{\mathop{\hspace*{0.4ex}\fo\hspace*{0.4ex}}}
\nek{\foa} [1] 
{\mathop{\hspace*{0.4ex}\fo^\ast_{#1}\hspace*{0.4ex}}}

\nek{\zfm} {\ZFL^{\text{\ubf--}}}
\nek{\mo} {\models}
\nek{\mm} {\gM}
\nek{\nn} {\gN}
\nek{\wmf} [1] {\vmf[\prl\vjPi{#1}]}

\nek{\bi} [1] {[#1]}

\renek{\cM} {\gM}

\nek{\imar}[1]{\marginnote[%\vspace{-1ex}%
\flushright\footnotesize%
\vspace*{-5.5ex}
{\scriptsize\rm #1$\Rightarrow$}%
\vspace*{1ex}]%
{\vspace*{-4.5ex}%
\flushleft\footnotesize%
{\scriptsize$\Leftarrow$\rm #1}%
}%
}%

\nek{\imaq}[1]{\marginnote[%\vspace{-1ex}%
\flushright\footnotesize%
\vspace*{-5.5ex}
{\scriptsize\rm #1$\rightarrow$}%
\vspace*{1ex}]%
{\vspace*{-4.5ex}%
\flushleft\footnotesize%
{\scriptsize$\leftarrow$\rm #1}%
}%
}%

\nek{\imas}[1]{\marginnote[% 
\flushright\footnotesize%
\vspace*{-6.5ex}
{\scriptsize\rm #1$\rightarrow$}%
]%
{\vspace*{-5.5ex}%
\flushleft\footnotesize%
{\scriptsize$\leftarrow$\rm #1}%
}%
}%

\nek{\imat}[1]{\marginnote[% 
\flushright\footnotesize%
\vspace*{-5ex}
{\scriptsize\rm #1$\rightarrow$}%
]%
{\vspace*{-4ex}%
\flushleft\footnotesize%
{\scriptsize$\leftarrow$\rm #1}%
}%
}%

\vyk{
\nek{\imae}[1]{\marginnote[%\vspace{-1ex}%
\flushright\footnotesize\vspace*{-10ex}%
$\rightarrow${#1}]%
{%\vspace{-1ex}%
\flushleft\footnotesize\vspace{-11ex}%
$\leftarrow${#1}}
}%
}

\nek{\imae}[1]{\marginpar[%\vspace{-1ex}%
\flushright\footnotesize\vspace{-4ex}%
$\rightarrow${#1}\vspace*{1ex}]%
{%\vspace{-1ex}%
\flushleft\footnotesize\vspace{-4ex}%
{#1}$\leftarrow$\vspace*{1ex}}
}%
\nek{\kmar}[1]
{\marginnote
[{\scriptsize\rm#1$\Rightarrow$}]%
{{\scriptsize\rm$\Leftarrow$#1}}}%

\nek{\msp} [2] {\ang{#1,#2}}
\nek{\fos} [1] {\text{\ubf FORC}[#1]}                      

\nek{\rn} [1] {\dd{#1}real name}
\nek{\qn} {real name}

\nek{\aut} {\text{\rm PERM}}
\nek{\ima} [2] {{#1}\text{\rm\hspace*{0.3ex}''\hspace*{-0.1ex}}{#2}}

\nek{\tba} {\tilde\ba}

\nek{\zqo} {\zq_0}
\nek{\zpo} {\zp_0}
\nek{\rco} {\rc_0}
\nek{\gao} {\ga_0}
\nek{\gai} {\ga_1}
\nek{\gad} {\ga_2}
\nek{\gat} {\ga_3}
\nek{\nuo} {\nu_0}
\nek{\nui} {\nu_1}
\nek{\zao} {\zeta_0}

\nek{\ad} {a.\hspace*{0.15ex}d.}
\nek{\sad} {s.\hspace*{0.15ex}a.\hspace*{0.15ex}d.}
\renek{\evd}{\ad}

\nek{\hh} {\boldsymbol h\hspace*{-0.00ex}}

\nek{\reg} [1] {{\res}_{\ge#1}}

\nek{\ofo} {\mathrel{{\|}\hspace*{-0.5ex}{-}}}

\renek{\kmar} [1] {}
\renek{\imae} [1] {}
\renek{\imat} [1] {}
\renek{\imas} [1] {}
\renek{\imaq} [1] {}
\renek{\imar} [1] {}

\begin{document}

\selectlanguage{english}

%$\mathbb A\nmi A_{\mathbb A\nmi A}$   {\small$\mathbb A$} 

\title
{The full basis theorem does not imply analytic wellordering}

\author 
%, Lyubetsky]
{
%Ali~Enayat\thanks{address, email}  
%\and
Vladimir~Kanovei
\thanks
{IITP RAS and MIIT, Moscow, Russia, \ 
{\tt kanovei@googlemail.com} --- 
contact author.} 
\thanks
{Thankful the Department of Philosophy, Linguistics and Theory of
Science at the University of Gothenburg and 
the Erwin Schrodinger International Institute for 
Mathematics and Physics (ESI) at Vienna 
for their hospitality and support in resp.\ May 2015 and 
December 2016.
\vyk{
Partial support of University of Gothenburg in 2015 
the Erwin Schrodinger International Institute for 
Mathematics and Physics (ESI) at Vienna in 2016 
is acknowledged.} 
} 
\and
Vassily~Lyubetsky
\thanks{IITP RAS, Moscow, Russia, \ {\tt lyubetsk@iitp.ru}.}
\thanks{Supported in part by RNF Grant \#14-50-00150.}
}

\date 
{\today}

\maketitle

\begin{abstract}
We make use of a finite support product of $\omi$ clones 
of the Jensen minimal $\varPi^1_2$ singleton forcing to 
define a model in which every non-empty 
lightface analytically definable 
set of reals contains a lightface analytically definable real 
(the full basis theorem),  
but there is no lightface  analytically definable wellordering of the continuum.
\end{abstract}
\vspace*{-5ex}

\np

{\scriptsize
\def\contentsname{}
\tableofcontents
}

\np

\parf{Introduction}
\las{int}

The uniformization problem, introduced by Luzin 
\cite{lhad2,lhad}, 
as well as the related basis problem,
are well known in modern set theory. 
(See Moschovakis~\cite{mDST}, Kechris~\cite{dst}, 
Hauser and Schindler \cite{hashi} for 
both older and more recent studies.)
In particular, it is known that every non-empty 
$\is12$ set of reals contains a $\id12$ real, 
but on the other hand, 
it is consistent that there exists a non-empty $\ip12$ set 
of reals containing even no ordinal-definable real.

%\cite{uspl}

The negative part of this result was strengthened in 
\cite{kl22} to the 
effect that the counter-example set $X\sq\bn$ is a  
$\ip12$ \dd\Eo equivalence class (hence, a countable set),
see related discussions at the Mathoverflow 
exchange desk\snos
{\label{snos1}
A question about ordinal definable real numbers. 
Mathoverflow, March 09, 2010. 
{\tt http://mathoverflow.net/questions/17608}. 
}
and at FOM\snos
{\label{snos2}%
Ali Enayat. Ordinal definable numbers. FOM Jul 23, 2010.
{\tt http://cs.nyu.edu/pipermail/fom/2010-July/014944.html}}.
Recall that $\Eo$ is an equivalence relation on $\bn$ defined 
so that $x\Eo y$ iff $x(n)=y(n)$ for all but finite $n$.

As for the positive direction, the most transparent way to 
get a basis result is to make use of an analytically definable 
wellordering $<$ of the reals, 
which enables one to pick the \dd<least 
real in each non-empty set of reals. 
This leads to the question: 
is the existence of an analytically definable 
wellordering $<$ of the reals independent of the basis theorem. 
We answer it in the positive:

\vyk{
Ali Enayat (Footnote~\ref{snos2}) conjectured that 
Question~\ref{vo2} can be solved in the positive by the   
finite-support product $\plo$ 
of countably many copies of the Jensen ``minimal $\ip12$ real 
singleton forcing'' $\dP$ defined in \cite{jenmin}\snos 
{\rit{Jensen's forcing} below, for the sake of brevity --- 
on this forcing, see also 28A in \cite{jechmill}.}.
Enayat demonstrated in \cite{ena} 
that a symmetric part of the \dd\plo generic 
extension of $\rL$, the constructible universe, 
definitely yields a model of $\zf$ 
(not a model of $\ZFC$!) 
in which there is a Dedekind-finite infinite \od\ set of 
reals with no \od\ elements. 

Following the mentioned conjecture, we proved in \cite{klOD} 
that indeed it is true in a \dd\plo generic extension of $\rL$  
that the set of \dd\dP generic reals is a countable non-empty  
$\ip12$ set with no $\od$ elements.\snos  
{We also proved in \cite{klE} that the existence of a 
$\ip12$ \dd\Eo class with no $\od$ elements is consistent 
with $\ZFC$, using a \dd\Eo invariant version of the Jensen 
forcing.}  
Using a finite-support product $\prod_{\xi<\omi}{\dP_\xi}\lom$, 
where the forcing notions $\dP_\xi$ are pairwise 
different clones of Jensen's forcing $\dP$, 
we answer Question~\ref{vo1} in the positive.
}

\bte
\lam{mt}
In a suitable generic extension of\/ $\rL$, 
it is true that in which every non-empty 
lightface analytically definable 
set of reals contains a lightface analytically definable real 
(the full basis theorem),  
but there is no lightface analytically definable 
%and even projective, 
wellordering of the continuum.

More precisely, there is a cardinal-preserving generic 
extension\/ $\rL[X]$ of\/ $\rL$, such that\/ 
$X=\sis{x_{\xi k}}{\xi<\omil\land k<\om}$, where each\/ 
$x_{\xi k}$ is a real in\/ $\dn$, and in addition
\ben
\Renu
\itlb{mt1}
if\/ $m<\om$ then the submodel\/ $\rL[X_m]$ admits a\/ 
%``good''\/ 
$\id1{m+3}$ wellordering of the reals 
of length $\omi$, where\/ 
$X_m=\sis{x_{\xi k}}{\xi<\omil\land k<m}\;;$

\itlb{mt2}
if\/ $m<\om$ then\/ $\bn\cap\rL[X_m]$ is a\/ 
$\is1{m+3}$ set in\/ $\rL[x]\;;$

\itlb{mt3}
if\/ $m<\om$ then\/ $\rL[X_m]$ is an
elementary submodel of\/ $\rL[x]$ \poo\ 
all\/ $\is1{m+2}$ formulas with reals in\/ 
$\rL[X_m]$ as parameters$;$

\itlb{mt4}
it is true in\/ $\rL[X]$ that there is no 
lightface analytically definable wellordering 
of the reals$.$
\een
\ete

To see that the additional claims imply the main claim 
(the full basis theorem),  
let, in $\rL[X]$, $Z\sq\bn$ be a non-empty $\is1{m+2}$
set of reals. 
Then $Z'=Z\cap\rL[X_m]$ is a $\is1{m+3}$ set by \ref{mt2}, 
and $Z'\ne\pu$ by \ref{mt3}. 
It remains to pick the least real in $Z'$ in the sense of 
the lightface $\id1{m+3}$ wellordering given by \ref{mt1}.

\parf{Comments}
\las{com}

To prove the theorem, we define, in $\rL$, a system of
forcing notions $\dP_{\xi k}$, $\xi<\omi$ and $k<\om$,
whose finite-support product $\dP=\prod_{\xi,k}\dP_{\xi k}$
adds an array $X=\sis{x_{\xi k}}{\xi<\omi,k<\om}$ of reals
$x_{\xi k}$ to $\rL$, such that 
\ref{mt1}, \ref{mt2}, \ref{mt3}, \ref{mt4} hold in $\rL[X]$.
\vyk{Then it is true in $\rL[X]$ that every non-empty
lightface $\is1{m+2}$ set of reals contains a real in
$\rL[X_m]$ by \ref{mt3}, 
and hence contains a $\id1{m+3}$ real by \ref{mt1}.
Adding \ref{mt4} in the mix, we establish the main statement 
of the theorem.
}

Regarding the history of this research, in goes down to 
Jensen~\cite{jenmin}, where a forcing 
$\dJ=\bigcup_{\al<\omi}\dJ_\al$ is defined in $\rL$, 
the constructible universe, such that each $\dJ_\al$ is 
a countable set of perfect trees in $\bse,$ the 
canonical \dd\dJ generic real is a single \dd\dJ generic 
real in the extension, and `being a \dd\dJ generic real' 
is a $\ip12$ property, so as a result we get a $\ip12$ 
nonconstructible singleton in any \dd\dJ generic 
extension of $\rL$.
See 28A in \cite{jechmill} for a more modern exposition of 
Jensen's forcing.

A nonconstructible $\ip12$ 
singleton also was defined in \cite{jsad} by 
means of the almost-disjoint forcing, yet the construction in 
\cite{jenmin} has the advantage of 
\rit{minimality} of \dd\dJ generic reals 
and some other advantages 
(as well as some disadvantages). 

Jensen's forcing construction (including its iterations) 
was exploited by Abraham \cite{abr,abr2}, including a 
definable minimal collapsing real. 
Another modification of Jensen's forcing construction in 
\cite{jj} yields such a forcing notion in $\rL$ that 
any extension of $\rL$, containing two generic reals 
$x\ne y$, necessarily satisfies $\omil<\omi$.
See \cite{baga,baga+} on some other modifications in 
coding purposes. 

A different modification of Jensen's forcing construction 
was engineered in \cite{ian78} in order to define an 
extension of $\rL$ in which, for a given $n\ge2$, there 
is a nonconstructible $\ip1n$ singleton while all $\is12$ 
reals are constructible. 
(An abstract appeared in \cite{dan75}.)
The idea is to complicate the inductive construction of 
Jensen's sequence $\vec\dJ=\sis{\dJ_\al}{\al<\omi}$ in $\rL$ 
by the requirement that it intersects any set, 
of a certain definability level,  
dense in the collection of all possible countable 
initial steps of the construction. 
The same \rit{inner genericity} idea, with respect to the 
Jensen -- Johnsbraten forcing notion in \cite{jj}, 
was developed in \cite{jjc}.

Such an \rit{inner genericity} modification of the 
Jensen -- Solovay almost-disjoint forcing \cite{jsad} 
was developed in \cite{h74} towards some great results 
which unfortunately have never been published in a 
mathematical journal.  
Except for a one result, a model in which the set of 
all analytically definable reals is equal to the set 
of all constructible reals, independently obtained in 
\cite{ian2}. 
We employ the inner definable genericity idea 
here in such a way that if $m<\om$ 
then the \dd mtail 
$\sis{\dP_{\xi k}}{\xi<\omi\land k\ge m}$ of the 
forcing construction, bears an amount of inner definable 
genericity which strictly depends on $m$. 
(See Definition~\ref{blo}, where a key concept is 
introduced.)

Ali Enayat (Footnote~\ref{snos2}) conjectured that 
some definability questions can be solved by     
finite-support products of Jensen's \cite{jenmin} 
forcing $\dJ$.
Enayat demonstrated in \cite{ena} 
that a symmetric part of the \dd{\dJ^\om}generic 
extension of $\rL$  
definitely yields a model of $\zf$ 
(not a model of $\ZFC$!) 
in which there is a Dedekind-finite infinite $\ip12$ set of 
reals with no \od\ elements. 
Following the conjecture, we proved in \cite{klOD} that 
indeed it is true in a \dd{\dJ^\om}generic extension 
of $\rL$  that the set of \dd{\dJ}generic reals is 
a countable non-empty $\ip12$ set with no $\od$ elements.  
We also proved in \cite{klE} that the existence of a 
$\ip12$ \dd\Eo class with no $\od$ elements is consistent 
with $\ZFC$, using a \dd\Eo invariant version of Jensen's 
forcing.  
We further employed another finite-support product of Jensen's 
forcing to define a generic extension of $\rL$ where 
there is a $\ip12$ set $P\sq \bn\ti\bn$ which has countable 
cross-sections $P_x=\ens{y}{\ang{x,y}\in P}$ and is 
non-uniformizable by any projective set 
\cite{KLunif}.

\back
The idea of making use of a suitable finite-support product 
of Jensen-like forcing notions in order to obtain a model, 
in which the full basis theorem holds 
but there is no lightface analytically definable 
wellordering of the continuum, 
was communicated to  an author of this paper (VK)  
by Ali Enayat in 2015, and 
we thank Ali Enayat for fruitful discussions and helpful 
ideas. 
\eack

\parf{The structure of the paper}
\las{lay}

The general organization of the paper is as follows. 
{\ubf Chapter~\ref{sek1}} contains a general formalism 
related to forcing by perfect trees and finite-support 
products, convenient for our goals. 
Following Jensen~\cite{jenmin}, 
we consider forcing notions of the form 
$\dP=\bigcup_{\al<\la}\dP_\al$, where $\la<\omi$ and 
each $\dP_\al$ is a countable set of perfect trees 
in $\bse.$ 
Each term $\dP_\al$ has to satisfy some routine 
conditions of \rit{refinement} with respect to the 
previous terms, in particular, to make sure that 
each $\dP_\al$ remains pre-dense at further steps. 
Also, each $\dP_\al$ has to \rit{lock} some dense sets 
in $\bigcup_{\xi<\al}\dP_\xi$ so that they 
remain pre-dense at further steps as well.   
And this procedure has to be extended from single 
forcing notions to their finite-support products. 
These issues are dealt with in {\ubf Chapter~\ref{sek2}}.

Then we consider real names with respect to 
finite-support products of perfect-tree forcing 
notions in {\ubf Chapter~\ref{sek3}}. 
Here the key issue is to make sure that if $\dP$ 
is a factor in a product forcing considered then 
there is no other \dd\dP generic real in the whole 
product extension except for the obvious one.

In {\ubf Chapter~\ref{sek4}}  
we define the forcing notion 
$\fP=\prod_{\xi<\omi,k<\om}\dP_{\xi k}$ 
to prove the main theorem, in the form of a limit 
of a certain increasing sequence of countable 
products of countable perfect-tree forcing notions.
Quite a complicated construction of this sequence 
in $\rL$ involves ideas related 
to diamond-style constructions, 
as well as to some sort of definable genericity, 
as explained above.

The forcing $\dP$ is not analytically definable; 
basically, each \dd kth layer 
$\sis{\dP_{\xi k}}{\xi<\omi}$ belongs to $\id1{k+4}$. 
But it is a key property that the \dd\fP forcing 
relation restricted to $\is1n$ formulas is essentially 
$\is1n$. 
We prove this in {\ubf Chapter~\ref{afn}}, with the help of 
an auxiliary forcing notion $\fof$.
We also establish the invariance of $\fof$ with respect to 
countable-support permutations of $\omi\ti\om$.

We finally prove Theorem~\ref{mt} in 
{\ubf Chapter~\ref{them}}, on the base of the results 
obtained in two previous chapters.

\vyk{
Using an appropriate generic extension of a 
submodel of the same model, similar to  
models considered in Harrington's unpublished notes 
\cite{h74}, we also prove 

\bte
\lam{Tsep}
In a suitable generic extension of\/ $\rL$, 
it is true  that there is a pair of 
disjoint lightface\/ $\ip13$ sets\/ $X,Y\sq\dn$, not 
separable by disjoint\/ $\fs13$ sets, and hence\/ 
$\fp13$ Separation and\/ 
$\ip13$ Separation fail.
\ete

This result was first proved by Harrington in \cite{h74} 
on the base of almost disjoint forcing of Jensen -- 
Solovay \cite{jsad}, 
and in this form 
has never been published, but was mentioned, \eg, in 
\cite[5B.3]{mDST} and \cite[page 230]{hin}.
A complicated alternative proof of Theorem~\ref{Tsep} 
can be obtained 
%in \cite{k83} 
with the help of \rit{countable-support} products and 
iterations of Jensen's forcing studied earlier in 
\cite{abr,k79,k83}. 
The \rit{finite-support} approach which we pursue here yields 
a significantly more compact proof, which still uses some 
basic constructions from \cite{h74}. 
As far as Theorem~\ref{Tun} is concerned, countable-support 
products and iterations hardly can lead to the 
countable-section non-uniformization results.

We recall that $\fp13$ Separation \rit{holds} in $\rL$. 
Thus Theorem~\ref{Tsep} in fact shows that the $\fp13$ 
Separation principle is destroyed in an appropriate generic 
extension of $\rL$.  
It would be interesting to find a generic extension in which, 
the other way around, the $\fs13$ Separation 
(false in $\rL$) holds.
This can be a difficult problem.
At least, the model used to prove Theorem~\ref{Tsep} does
not help: we prove (Theorem~\ref{Tses} below) that any
pair of disjoint $\fs13$ sets, non-separable by
disjoint $\fp13$ sets in $\rL$, remains $\fs13$ and
non-separable by disjoint $\fp13$ sets  in the extension.
}

\gla{Basic constructions}
\las{sek1}

We begin with some basic things: 
perfect trees in the Cantor space $\dn$, 
perfect tree forcing notions 
(those which consist of perfect trees), 
their finite-support products, and a 
splitting construction of perfect trees.

\parf{Perfect trees} 
\las{tre1}

Let $\bse$ be the set of all \rit{strings} (finite sequences) 
\index{strings}%
of numbers $0,1$.
If $t\in\bse$ and $i=0,1$ then 
$t\we k$ is the extension of $t$ by $k$. 
\index{zzt^k@$t\we k$}%
If $s,t\in\bse$ then $s\sq t$ means that $t$ extends $s$, while 
$s\su t$ means proper extension. 
\index{zzstsu@$s\su t$}%
\index{zzstsq@$s\sq t$}%
If $s\in\bse$ then $\lh s$ is the length of $s$,  
\index{zzlhs@$\lh s$}%
and $2^n=\ens{s\in\bse}{\lh s=n}$ (strings of length $n$).%
\index{zz2^n@$2^n$}%

A set $T\sq\bse$ is a \rit{tree} iff 
\index{tree}%
%it is an initial segment, that is, 
for any strings $s\su t$ in $\bse$, if $t\in T$ then $s\in T$.
Every non-empty tree $T\sq\bse$ contains the empty 
string $\La$. 
\index{strings!empty string, $\La$}%
\index{zzLa@$\La$}%
If $T\sq\bse$ is a tree and $s\in T$ then put 
$T\ret s=\ens{t\in T}{s\sq t\lor t\sq s}$. 
\index{zzTis@$T\ret s$}%
%this is a tree as well.

Let $\pet$ be the set of all \rit{perfect} trees 
\index{zzPT@$\pet$}%
$\pu\ne T\sq \bse$. 
\imar{pet}%
Thus a non-empty tree $T\sq\bse$ belongs to $\pet$ iff 
it has no endpoints and no isolated branches. 
Then there is a largest string $s\in T$ such that 
$T=T\ret s$; it is denoted by $s=\roo T$   
(the {\it stem\/} of $T$); 
\index{stem, $\roo T$}%
we have $s\we 1\in T$ and $s\we 0\in T$ in this case.

\bdf
[perfect sets]
\lam{body}
If $T\in\pet$ then   
$ 
[T]=\ens{a\in\dn}{\kaz n\,(a\res n\in T)} 
$  
\index{zzT[]@$[T]$}% 
is the set of all \rit{paths through $T$}, 
a perfect set in $\dn$. 
Conversely if $X\sq\dn$ is a perfect set then 
$\tre X=\ens{a\res n}{a\in X\land n<\om}\in\pet$
\index{tree!$\tre X$}%
\index{zztreeX@$\tre X$}%
and $[\tre X]=X$.

Trees $T,S\in\pet$ are \rit{almost disjoint}, \rit{\ad{}} 
for brevity, iff 
\kmar{evd}
\index{tree!almost disjoint}%
\index{tree!almost disjoint@\ad\ trees}%
\index{almost disjoint@\ad\ (trees)}%
the intersection $S\cap T$ is finite; 
this is equivalent to just $[S]\cap[T]=\pu$.
\edf

The \rit{simple splitting\/} of a tree $T\in\pet$ consists 
\index{splitting}%
of smaller trees
$$
\raw T0=T\ret{\roo T\we 0}
\quad\text{and}\quad
\raw T1=T\ret{\roo T\we 1} 
$$
\index{zzT->i@$\raw T i$}%
in $\pet$, so that $[\raw Ti]=\ens{x\in[T]}{x(h)=i}$, 
where $h=\lh{\roo T}$.
We let
$$
\raw Tu= \raw{\raw{\raw{\raw T{u(0)}}{u(1)}}{u(2)}\dots}{u(n-1)}
$$
\index{zzT->u@$\raw T u$}%
for each string $u\in\bse\yd \lh u=n$; 
and separately $\raw T\La=T$.

\ble
\lam{splem}
Suppose that\/ $T\in\pet$. 
%\lam{clops}
Then$:$
\ben
\renu
\itlb{splem1}
if\/ $u\in\bse$ then there is a string\/ $s\in\bse$ such that\/ 
$\raw Tu=\req Ts\;;$  

\itlb{splem2}
if\/ $s\in\bse$ then there is a string\/ $u\in\bse$ 
such that\/ $\req Ts=\raw Tu\;;$  

\itlb{splem3}
if\/ $\pu\ne U\sq[T]$ 
is a (relatively) open subset of\/ $[T]$, or at least\/ 
$U$ has a non-empty interior in\/ $[T]$, then
there is a string\/ $s\in T$ such that\/ $\req Ts\sq U$.
\qed 
\een
\ele

If $T\in\pet$ and $a\in\dn$ then the intersection 
$\raw Ta=\bigcap_{n<\om}\raw T{a\res n}=\ans{\gom T(a)}$ 
is a singleton, and the map $\gom T$ is a 
\index{zzhT@$\gom T$}%
\index{canonical homeomorphism!zzhT@$\gom T$}%
\rit{canonical homeomorphism} from $\dn$ onto $[T]$.  
Accordingly if $S,T\in\pet$ then the map 
$\gon ST(x)=\gom T({\gom S}\obr(x))$ is a 
\rit{canonical homeomorphism} from $[S]$ onto $[T]$.
\index{zzhST@$\gon ST$}%
\index{canonical homeomorphism!zzhST@$\gon ST$}%

\parf{Perfect tree forcing notions} 
\las{tre2}

A {\ubf perfect-tree forcing notion} is any non-empty set 
\index{perfect-tree forcing, $\ptf$}
$\dP\sq\pet$ such that if $s\in T\in\dP$ then $T\ret s\in \dP$, 
or equivalently, by Lemma~\ref{splem}, 
if $u\in \bse$ then $\raw Tu\in \dP$.
Let $\ptf$ be the set of all such forcing notions  
\index{zzptf@$\ptf$}
$\dP\sq\pet$.

\bpri
\lam{cloL}
If $s\in\bse$ then the 
tree $\bi s=\ens{t\in\bse}{s\sq t\lor t\sq s}$ 
\index{tree!@$\bi s$}%
\index{zzs[]@$\bi s$}%
%\index{0Is@$I_s$}% 
belongs to $\pet$.
%and  $T[s]=\raw{(\bse)}s=\req{(\bse)}s\,,\:\kaz s$. 
%
%Generally if\/ $T\in\pet$ and\/ $2^n\sq T$ then\/ 
%$\raw Ts=\req Ts$ for all\/ $s\in 2^n$. 
%
The set $\dpo=\ens{\bi s}{s\in\bse}$ of all such trees 
(the Cohen forcing) 
\index{Cohen forcing, $\dpo$} 
\index{zzPcoh@$\dpo$} 
is a regular perfect-tree forcing notion.
\epri

\ble
\lam{regfn11}
Let\/ $\dP\in\ptf$.
If\/ $T\in\dP$ and a set\/ $X\sq[T]$ is (relatively) open
(resp., clopen) in\/ $[T]$, then there is a countable
(resp., finite) set\/ $\cS$ of pairwise \evd\ trees\/ $S\in\dP$, 
satisfying\/ $\bigcup_{S\in\cS}[S]=X$.\qed
\ele

\ble
\label{suz}
\ben
\renu
\itlb{suz1}
If\/ 
%\imar{suz}
$s\in T\in\dP\in\ptf$ then\/ $T\ret s\in\dP$.

\itlb{suz2}
If\/ $\dP,\dP'\in\ptf$, $T\in\dP$, $T'\in\dP'$, then 
there are trees\/ $S\in\dP$, $S'\in\dP'$ such that\/ 
$S\sq T$, $S'\sq T'$, and\/ $[S]\cap[S']=\pu$.
\een
\ele
\bpf
\ref{suz1} use Lemma~\ref{splem}.
\ref{suz2} 
If $T=T'$ then let $S=\raw T0$, 
$S'=\raw T1$.
If say $T\not\sq T'$ then let $s\in T\bez T'$, 
$S=T\ret{s}$, and simply $S'=T'$.
\epf

\bdf
\lam{ptf2}
A set $\dA\sq\pet$ is an \rit{antichain} iff any trees 
\index{antichain}
$T\ne T'$ in $\dA$ are \evd, that is, 
$[T]\cap[T']=\pu$.
A forcing notion $\dP\in\ptf$ is:  
\bde
\item[\rit{small}\rm,]  
if it is countable; 
\index{perfect-tree forcing, $\ptf$!small}% 
\index{small}%  
  
\item[\rit{special}\rm,]  
if there is 
\index{perfect-tree forcing, $\ptf$!special}% 
\index{special}%  
an antichain $\dA\sq\dP$ 
such that $\dP=\ens{T\ret s}{s\in T\in \dA}$ --- 
note that $\dA$ 
\index{perfect-tree forcing, $\ptf$!base, $\baz\dP$}% 
\index{base, $\baz\dP$}% 
\index{zzbazP@$\baz\dP$}% 
is unique if exists; we write 
$\dA=\baz\dP$ (the \rit{base} of $\dP$);
  
\item[\rit{regular}\rm,] 
if for any $S,T\in\dP$, the intersection  
\index{perfect-tree forcing, $\ptf$!regular}% 
\index{regular}% 
$[S]\cap [T]$ is clopen in $[S]$ or clopen in $[T]$ 
(or clopen in both $[S]$ and $[T]$).
\qed
\ede
\eDf

\ble
\lam{regfn}
Let\/ $\dP\in\ptf$.
If\/ $\dP$ is special and\/
$S,T\in\dP$ are not \evd, then they are comparable$:$
$S\sq T$ or\/ $T\sq S$. 

If\/ $\dP$ is special then\/ $\dP$ is regular.
If\/ $\dP$ is regular, then
\ben
\renu
\itlb{regfn3}
if\/ $S,T\in\dP$ are not
%\imar{regfn3}
\evd, then they are compatible in\/ $\dP$, that is,
there is a tree\/ $R\in\dP$ such that\/ $R\sq S\cap T$.

\itlb{regfn4}
if\/
$S_1,\dots,S_k\in\dP$ then there is a finite set of 
pairwise \evd\ trees $R_1,\dots,R_n\in \dP$ such that 
$[S_1]\cap\ldots\cap[S_k]=[R_1]\cup\ldots\cup [R_n]$. 

\itlb{regfn5}
if\/ $\cS_1,\dots,\cS_k$ are
finite collections of trees in\/ $\dP$ then there
is a finite set of trees $R_1,\dots,R_n\in \dP$ such that 
${\bigcup_{S\in\cS_1}[S]}\cap\ldots\cap
{\bigcup_{S\in\cS_k}[S]}=
[R_1]\cup\ldots\cup [R_n]$,
and for any\/ $\cS_i$ and\/ $R_j$,
there is\/ $S\in\cS_i$ such that\/ $R_j\sq S$. 
\een
\ele
\bpf
\ref{regfn5}
Apply \ref{regfn4} to every set of the form
$[S_1]\cap\ldots\cap[S_k]$, where $S_i\in\cS_i$, $\kaz i$,
then gather all
trees $R_i$ obtained in one finite set.
\epf

\bre
\lam{tfn}
Any set $\dP\in\ptf$ can be considered as a forcing notion 
(if $T\sq T'$ then $T$ is a stronger condition); 
then $\dP$ adds a real $x\in\dn$.
\ere

\ble
\lam{gx}
If a set $G\sq\dP$ is generic over a ground set universe\/ 
$\rV$ \sko{resp., over a transitive model, \eg\ $\rL$}
then 
\ben
\renu
\itlb{gx1} 
the intersection\/ $\bigcap_{T\in G}[T]$ contains a 
single real\/ $x=x[G]\in\dn$, and 

\itlb{gx2} 
this real\/ $x$ is\/ 
\dd\dP{\ubf generic}, 
\index{real!Pgeneric@\dd\dP generic}%
\index{Pgeneric@\dd\dP generic}%
\index{generic!Pgeneric@\dd\dP generic real}%
in the sense that if\/ $D\sq\dP$ is 
dense in\/ $\dP$ and belongs to\/ $\rV$ 
\sko{resp., to the ground model} 
then\/ $x\in \bigcup_{T\in D}[T]$.\qed
\een
\ele

As usual, a set $D\sq\dP$ is:
\bit
\item[$-$]
\rit{open} in $\dP$, if for any trees $T\sq S$ in $\dP$, 
\index{set of trees!open}%
\index{open}%
$T\in D\imp S\in D$;

\item[$-$]
\rit{dense} in $\dP$, if for any $T\in \dP$  
\index{set of trees!dense}%
\index{dense}%
there is $S\in D$, $S\sq T$;

\item[$-$]
\rit{pre-dense} in $\dP$, if the set 
\index{set of trees!pre-dense}%
\index{pre-dense}%
$D'=\ens{T\in \dP}{\sus S\in D\,(T\sq S)}$ 
is dense in $\dP$.
\eit

\parf{Splitting construction} 
\las{spe} 

%This section starts the proof of Theorem~\ref{xist}. 
We proceed with an important splitting/fusion 
construction of perfect trees by means of infinite 
splitting systems of such trees.      

\bdf
\lam{sped}
%If $\dP\in\ptf$ then let $\fss \dP$ be the set of all 
Let $\fso$ be the set of all 
\rit{finite splitting systems}, that is,  
\index{finite splitting system, $\fso$}% 
\index{zzfss@$\fso$}% 
systems of the form 
$\vpi=\sis{T_s}{s\in2^{\le n}}$, where 
$n=\vys\vpi<\om$ (the height of $\vpi$), 
\index{finite splitting system, $\fso$!height, $\vys\vpi$}% 
\index{finite splitting system, $\fso$!empty system, $\jLa$}% 
each value $T_s=T^\vpi_s=\vpi(s)$ is a tree in $\pet$, and
\ben
\fenu
\itlb{spe2}
if $s\in 2^{<n}$ and $i=0,1$ (so $s\we i\in2^{\le n}$)
then $T_{s\we i}\sq \raw{T_s}i$  --- 
%\imar{spe2}
it easily follows that 
$[T_{s\we0}]\cap [T_{s\we1}]=\pu$.
%\qed
\een
We add the \rit{empty system} $\jLa$ to $\fso$, with $\vys\jLa=-1$. 
\edf

A tree $T$ \rit{occurs in\/ $\vpi\in\fso$} if 
\index{occurs}% 
\index{finite splitting system, $\fso$!tree occurs in}% 
$T=\vpi(s)$ for some $s\in 2^{\le\vys\vpi}$.
If all trees occurring in $\vpi$ belong to some $\dP\in\ptf$ then
say that $\vpi$ is a finite splitting system
\rit{over $\dP$}, symbolically $\vpi\in\fss \dP$. 
\index{finite splitting system, $\fso$!over $\dP$, $\fss \dP$}% 
\index{zzfssP@$\fss \dP$}% 

Let $\vpi,\psi$ be systems in $\fso$.
Say that 
%\bit
%\item[$-$]
$\vpi$ \rit{extends} $\psi$, symbolically $\psi\cle\vpi$, if 
\index{finite splitting system, $\fso$!extension, $\psi\cle\vpi$}% 
$n=\vys\psi\le\vys\vpi$ and $\psi(s)=\vpi(s)$ for 
all $s\in2^{\le n}$, and \rit{properly extends}, $\psi\cls\vpi$,
\index{finite splitting system, $\fso$!extension proper, $\psi\cls\vpi$}%
if in fact 
$\vys\psi<\vys\vpi$ strictly.
\vyk{
\item[$-$]
\rit{properly extends} $\psi$, 
symbolically $\psi\cls\vpi$, if in 
addition $\vys\psi<\vys\vpi$;

\item[$-$]
\rit{reduces} $\psi$, if $n=\vys\psi=\vys\vpi$, 
$\vpi(s)\sq\psi(s)$ for all $s\in 2^{\vys\vpi},$ and 
$\vpi(s)=\psi(s)$ for all $s\in 2^{<\vys\vpi}$.
\eit
In other words, reduction allows to shrink trees in the top 
layer of the system, but keeps intact those in the lower 
layers.
}

Each system $\vpi\in\fss\dP$ with $\vys\vpi=0$ consists
essentially of a single tree $T^\vpi_\La\in\dP$. 
%, where $\La$ is the empty string.  
The next lemma provides systems of arbitrary height.

\ble
\lam{n+1}
Assume that\/ $\dP\in\ptf$.
If\/ $n\ge1$ and\/ $\psi=\sis{T_s}{s\in2^{\le n}}\in\fss\dP$ 
then there is a system\/ 
$\vpi=\sis{T_s}{s\in2^{\le n+1}}\in\fss\dP$
which properly extends\/ $\psi$.
\ele
\bpf
If $s\in2^{n}$ and $i=0,1$ then let 
$T_{s\we i}=\raw{T_s}{i}$. 
\epf

The next well-known lemma belongs to the type of 
\rit{splitting/fusion} lemmas widely used in connection
with the perfect set forcing and some similar forcings.

\ble
\lam{infty}
Let\/ $\dP\in\ptf$.
Then there is an\/ \dd\cls in\-creas\-ing sequence\/ 
$\sis{\vpi_n}{n<\om}$ of systems in\/ $\fss\dP$. 
And if\/   $\sis{\vpi_n}{n<\om}$ is such then$:$
\ben
\renu
\itlb{infty1}
the limit system\/ 
$\vpi=\bigcup_n\vpi_n=\sis{T_s}{s\in\bse}$
satisfies\/ \ref{spe2} of Definition~\ref{sped}
%Definition~\ref{sped} 
on the whole domain of strings\/ $s\in\bse\,;$ 

\itlb{infty2}
$T=\bigcap_n\bigcup_{s\in2^n}T_s$ is 
a perfect tree in\/ $\pet$ 
%{\rm(yet not necessarily in $\dP$)}, 
and\/ $[T]=\bigcap_n\bigcup_{s\in2^n}[T_s]\,;$

%If\/ $s\in\bse$ then\/
%$T'_s=\bigcap_{n\ge \lh s}\bigcup_{t\in2^n,s\sq t}T_t\in\pet$,
%$T'_s=T\cap T_s$, 
%$[T'_s]=\bigcap_{n\ge \lh s}\bigcup_{t\in2^n,s\sq t}[T_s]$,
%and\/ $T'_s=T\ret u$, where\/ $u=\roo{T_s}$.
\itlb{infty3}  
if\/ $u\in\bse$ then\/
$\raw Tu=T\cap T_u=\bigcap_{n\ge \lh u}\bigcup_{s\in2^n,u\sq s}T_s$.
\qed
\een
\ele

\parf{\Muf s and \mut s} 
\las{muft}
             
We'll systematically make use of finite support products 
of perfect tree forcings in this paper. 
The following definitions introduce suitable notation.

Call a {\ubf\muf} 
\index{multiforcing}  
any map $\jpi:\abc\jpi\to\ptf$, where 
$\abc\jpi=\dom\jpi\sq\omi\ti\om$. 
%is at most countable. 
Thus each set $\dpi\xi k$, 
\kmar{dpi xi k}
$\ang{\xi,k}\in\abc\jpi$, 
%$=\jpi(\xi,k)$ 
is a perfect tree forcing notion. 
Such a $\jpi$ is:
\bit
\item[$-$]
\rit{small}, if both $\abc\jpi$ and each 
\index{multiforcing!small}%  
\index{small}%  
forcing $\dpi\xi k$, $\ang{\xi,k}\in\abc\jpi$, are countable;

\item[$-$]
\rit{\mdi}, 
\index{multiforcing!special}%  
\index{special}%  
if each $\dpi\xi k$ is special in the sense of 
Definition~\ref{ptf2};

\item[$-$] 
\rit{\mre}, if each $\dpi\xi k$ is regular, 
\index{multiforcing!regular}%  
\index{regular}%  
in the sense of Definition~\ref{ptf2}.
\eit
Let $\mfp$ be the set of all \muf s. 
\kmar{mfp}

Let a {\ubf\mut} be any map  
\index{multitree}% 
$\zp:\abc\zp\to\pet$, such that $\abc\zp=\dom\zp\sq\omi\ti\om$ 
is finite and each value 
$\zd \zp\xi k=\zp(\xi,k)$ is a tree in $\pet$.
\index{tree!Tpxik@$\zd \zp\xi k$}% 
\index{zzTpxik@$\zd \zp\xi k$}% 
\index{multitree!zzTpxik@$\zd \zp\xi k$}% 
\kmar{zd zp xi k}
In this case we define a cofinite-dimensional
perfect cube in $2^{\omi\ti\om}$
$$
\bay{rcl}
[\zp]
\index{zzp[]@$[\zp]$}%
\index{multitree!zzp[]@$[\zp]$}%
&=&\ens{x\in2^{\omi\ti\om}}
{\kaz \ang{\xi,k}\in\abs\zp\,(x(\xi,k)\in [\zd\zp\xi k])}=
\\[1.5ex] 
&=&\ens{x\in2^{\omi\ti\om}}
{\kaz \ang{\xi,k}\in\abs\zp\;\kaz m\,
(x(\xi,k)\res m\in \zd\zp\xi k)}\,. 
\eay
$$
Let $\md$ be the set of all \mut s. 
\index{multitree!$\md$}% 
\index{zzMT@$\md$}% 
\kmar{md, leq}
% of size $\vt$.
We order $\md$ componentwise: $\zq\leq\zp$ 
($\zq$ is stronger) 
iff $\abc\zp\sq\abc\zq$ and 
$\zd\zq\xi k\sq\zd\zp\xi k$ for all $\ang{\xi,k}\in\abc\zp$;
this is equivalent to $[\zq]\sq[\zp]$, so that stronger
\mut s correspond to smaller cubes. 
The weakest \mut\ $\jLa\in\md$ is just the empty map; 
\index{multitree!empty multitree, $\jLa$}% 
\kmar{jLa}
$\abc\jLa=\pu$ and $[\jLa]=2^{\omi\ti\om}$. 

\Mut s $\zp,\zq$ are \rit{somewhere almost disjoint}, 
or \rit{\sad}, 
if, for at least one pair of indices
\kmar{sad}
$\ang{\xi,k}\in\abc\zp\cap\abc\zq$,    
\index{sad@\sad, somewhere \ad}%
\index{somewhere ad@\sad, somewhere \ad}%
\index{somewhere almost disjoint, \sad}%
\index{multitree!somewhere almost disjoint, \sad}%
the trees $\zd \zp\xi k$, $\zd \zq\xi k$ are \evd,
that is,  $[\zd \zp\xi k]\cap [\zd \zq\xi k]=\pu$, 
or equivalently, $\zd \zp\xi k\cap \zd \zq\xi k$ is 
finite.

\bcor
[of Lemma~\ref{regfn}\ref{regfn3}]
\lam{sad}
If\/ $\jpi$ is a regular \muf\ and \mut s\/ 
$\zp,\zq\in\mt\jpi$ are not \sad, then\/ $\zp,\zq$ 
are compatible in\/ $\mt\jpi$, so that there is a 
\mut\ $\zr\in\mt\jpi$ with\/ 
$\zr\leq\zp$, $\zr\leq\zq$.\qed
\ecor

If $\jpi$ is a \muf\ then
a \rit{\dd\jpi\mut} is any \mut\ $\zp$ with 
\index{multitree!pimultitree@\dd\jpi\mut}% 
$\abc\zp\sq\abc\jpi$ and  
$\zd \zp\xi k\in\dpi\xi k$ for all 
$\ang{\xi,k}\in\abc\zp$.
Let $\mt\jpi$ 
\kmar{mt jpi}
be the set  of all \dd\jpi\mut s;
\index{multitree!$\mt\jpi$}% 
\index{zzMTpi@$\mt\jpi$}% 
it is equal to the finite support product 
$\prod_{\ang{\xi,k}\in\abc\jpi}\dpi\xi k$.

The following is similar to Lemma~\ref{regfn}\ref{regfn5}.

\ble
\lam{regfm}
If a \muf\/ $\jpi$ is regular, $\xi\sq\abs\jpi$
is finite, and\/ $U_1,\dots,U_k$ are
finite collections of \mut s in\/ $\mt\jpi$
with\/ $\abs\zp=\xi$ for all\/ $\zp\in\bigcup_iU_i$, 
then there is a finite set of \mut s\/
$\ju_1,\dots,\ju_n\in\mt\jpi$ such that\/
$\abs{\ju_j}=\xi$, $\kaz j$,\pagebreak[0]  
$$
\textstyle
{\bigcup_{\zp\in U_1}[\zp]}\cap\ldots\cap
{\bigcup_{\zp\in U_k}[\zp]}=
[\ju_1]\cup\ldots\cup [\ju_n],
$$
and for any\/ $U_i$ and\/ $\ju_j$, there is\/
$\zp\in U_i$ such that\/ $[\ju_j]\sq [\zp]$. 
\qed
\ele

We consider sets of the form $\mt\jpi$ in the 
role of {\ubf product forcing notions}. 
A set $D\sq\mt\jpi$ is:
\bit
\item[$-$]
\rit{open} in $\mt\jpi$, if for any $\zp\leq\zq$ in $\mt\jpi$, 
\index{set of multitrees!open}%
\index{open}%
$\zq\in D\imp\zp\in D$;

\item[$-$]
\rit{dense} in $\mt\jpi$, if for any $\zp\in \mt\jpi$, 
\index{set of multitrees!dense}%
\index{dense}%
there is $\zq\in D$, $\zq\leq\zp$;

\item[$-$]
\rit{pre-dense} in $\mt\jpi$, if the set 
\index{set of multitrees!pre-dense}%
\index{pre-dense}%
$D'=\ens{\zp\in \mt\jpi}{\sus\zq\in D\,(\zp\leq\zq)}$ 
is dense in $\mt\jpi$.
\eit

\bre
\lam{adds}
As a forcing notion, each 
$\mt\jpi$ adds an array   
$\sis{x_{\xi k}}{\ang{\xi,k}\in\abc\jpi}$ of reals, 
where each real $x_{\xi k}\in\dn$ 
is a \dd{\dpi\xi k}generic real. 
Namely if a set $G\sq\mt\jpi$ is generic over the ground 
set universe $\rV$ then each factor 
$$
G(\xi,k)=\ens{\zd \zp\xi k}{\zp\in G\land\ang{\xi,k}\in\abc\zp}
\sq  \dpi\xi k
$$ 
(where $\ang{\xi,k}\in\abc\jpi$) 
is accordingly a set \dd{\dpi\xi k}generic over $\rV$, the real 
$x_{\xi k}=x_{\xi k}[G]=x[G(\xi,k)]\in\dn$ is the only 
\index{zzxxikG@$x_{\xi k}[G]$}
real satisfying $x_{\xi k}\in\bigcap_{T\in G(\xi,k)}[T]$, and 
$x_{\xi k}$ is \dd{\dpi\xi k}generic over $\rV$ as in 
Lemma~\ref{gx}.
\ere

The reals of the form $x_{\xi k}[G]$ will be called 
\rit{principal generic reals\/} in $\rV[G]$. 
\index{reals!principal generic reals, $x_{\xi k}[G]$}%
\index{principal generic reals, $x_{\xi k}[G]$}%

\bdf
%[componentwise union]
\lam{kw}
A \rit{componentwise union} of \muf s 
$\jpi,\jqo$ is a \muf\ $\jpi\kw\jqo$ satisfying 
\kmar{kw}
$\abc{(\jpi\kw\jqo)}=\abc\jpi\cup\abc\jqo$ and 
\index{componentwise union!Ucw@$\jpi\kw\jqo$}%
\index{zzUcw@$\jpi\kw\jqo$}
$$
(\jpi\kw\jqo)(\xi,k)=
\left\{
\bay{rcl}
\jpi(\xi,k),&\text{whenever}&\ang{\xi,k}\in\abc\jpi\bez\abc\jqo\\[1ex]
\jqo(\xi,k),&\text{whenever}&\ang{\xi,k}\in\abc\jqo\bez\abc\jpi\\[1ex]
\jpi(\xi,k)\cup\jqo(\xi,k),
&\text{whenever}&\ang{\xi,k}\in\abc\jpi\cap\abc\jqo 
\eay
\right.
$$ 
\vyk{
resp., $\jqo(\xi,k)$,
resp., $\jpi(\xi,k)\cup\jqo(\xi,k)$,
in cases $\ang{\xi,k}\in\abc\jpi\bez\abc\jqo$, 
resp., $\ang{\xi,k}\in\abc\jqo\bez\abc\jpi$,
resp., $\ang{\xi,k}\in\abc\jpi\cap\abc\jqo$. 
}%
Similarly, if $\vjpi=\sis{\nor\jpi\al}{\al<\la}$ is a sequence 
of \muf s then 
define a \muf\ 
$\jpi=\bkw\vjpi=\bkw_{\al<\la}\nor\jpi\al$ so that 
\kmar{bkw}
\index{componentwise union!Ucwb@$\bkw\vjpi=\bkw_{\al<\la}\nor\jpi\al$}%
\index{zzUcwb@$\bkw_{\al<\la}\nor\jpi\al$}
$\abc\jpi=\bigcup_{\al<\la}\abc{\nor\jpi\al}$ and if 
$\ang{\xi,k}\in\abc\jpi$ then 
$\dpj\xi k{\jpi}=
\bigcup_{\al<\la,\:\ang{\xi,k}\in\abc{\nor\jpi\al}}
\dpj\xi k{\nor\jpi\al}$.
\edf

\parf{\Mus s} 
\las{mus}

The next definition introduces 
\rit{\mus s}, a multi version of the 
splitting/fusion technique of Section~\ref{spe}, 
whose intention is to define suitable \muf s, 
as will be shown in Section~\ref{jex} below.

\bdf
\lam{muss}
A {\ubf\mus} 
\index{multisystem}%
%of size $\vt=\siz\jfi$ 
is any map $\jfi:\abc\jfi\to\fso$, 
such that 
$\abc\jfi\sq\omi\ti\om\ti\om$ is finite.
This amounts to 
\ben
\nenu
\item 
the map
$\qh\jfi\xi km=\vys{\jfi(\xi,k,m)}:\abc\jfi\to\om$,
\index{multisystem!hfixikm@$\qh\jfi\xi km$}% 
\index{zzhfixikm@$\qh\jfi\xi km$}% 
%such that the support 
%$\ens{\ang{\xi,k,m}}{\qh\jfi\xi km\ne -1}=\abs\jfi$
%is finite, 
and 

\item
the finite collection of trees
$\qT\jfi\xi km s=\jfi(\xi,k,m)(s)$,
\index{multisystem!Tfixikms@$\qT\jfi\xi km s$}% 
\index{zzTfixikms@$\qT\jfi\xi km s$}% 
where $\ang{\xi,k,m}\in\abc\jfi$
and $s\in 2^{\le \qh\jfi\xi km}$, such that if
$\ang{\xi,k,m}\in\abs\jfi$ then
$\jfi(\xi,k,m)=\sis{\qT\jfi\xi km s}{s\in 2^{\le \qh\jfi\xi km}}$ 
is a finite splitting system in $\fso$.
\een
If $\jpi$ is a \muf, $\abc\jfi\sq{(\abc\jpi)}\ti\om$, and
$\jfi(\xi,k,m)\in\fss{\dpi\xi k}$ for all 
$\ang{\xi,k,m}\in\abc\jfi$
(or equivalently, 
$\qT\jfi\xi km s\in \dpi\xi k$ whenever
$\ang{\xi,k,m}\in\dm\jfi$ and $s\in 2^{\le \qh\jfi\xi km}$),
then say that $\jfi$ is a \dd\jpi{\ubf\mus},  
$\jfi\in\ms\jpi$.
\index{multisystem!pimultisystem@\dd\jpi\mus}% \edf
\index{multisystem!MSpi@$\ms\jpi$}% \edf
\index{zzMSpi@$\ms\jpi$}% 
\edf

Let $\jfi,\jsi$ be multisystems.  
Say that  
$\jfi$ \rit{extends} $\jsi$, 
symbolically $\jsi\cle\jfi$, 
\index{multisystem!extension, $\jsi\cle\jfi$}% \edf
if $\abc\jsi\sq\abc\jfi$, and, 
for every $\ang{\xi,k,m}\in\abc\jsi$,
$\jfi(\xi,k,m)$ extends $\jsi(\xi,k,m)$, that is, 
$\qh\jfi\xi km\ge \qh\jsi\xi km$ and  
$\qT\jfi\xi km s=\qT\jsi\xi km{s}$ for all 
$s\in 2^{\le\qh\jsi\xi km}$.

It will be demonstrated in Section~\ref{jex} that a 
suitably increasing infinite sequence 
$\jfi_0\cle\jfi_1\cle\jfi_2\cle\dots$ of \mus s in 
some $\ms\jpi$ leads to a ``limit'' \muf\ $\jqo$ 
with $\abc\jqo=\bigcup_n\abc{\jfi_n}$, such that each 
factor $\dqo\xi k$, $\ang{\xi,k}\in\abc\jpi$, is filled 
in by trees $Q_{\xi k,m}$, $m<\om$, in such a way, that 
the \dd{(\xi,k,m)}components of the systems $\jfi_n$ 
are responsible for the construction of the 
tree $Q_{\xi k,m}$. 

The next lemma introduces different ways to extend 
a given \mus. 

Say that a \mus\ $\vpi$ is \rit{2wise disjoint} if
\index{multisystem!2wise disjoint}%  
$[\qT\vpi\xi km{s}]\cap [\qT\vpi\et \ell n{t}]=\pu$ 
for all triples $\ang{\xi,k,m}\ne\ang{\et,\ell,n}$
in $\abc\vpi$ 
and all $s\in 2^{\qh\jfi\xi km}$ and
$t\in 2^{\qh\jfi\et\ell n}$.  

\ble
\lam{nwm}
Let\/ $\jpi$ be a \muf\ and\/ $\jfi\in\ms\jpi$. 
\ben
\renu
\itlb{nwm1}
If\/ $\ang{\xi,k,m}\in\abc\jfi$ and\/ $h=\qh\jfi\xi km$ then
%\imar{nwm1}
the extension\/ $\jsi$ of\/ $\jfi$ by\/
$\qh{\jsi}\xi km=h+1$ and\/
$\qT{\jsi}\xi km{s\we i}=\raw{\qT{\jfi}\xi km{s}}i$
for all\/ $s\in2^h$ and\/ $i=0,1$,
belongs to\/
$\ms\jpi$ and\/ $\jfi\cle\jsi$.

\itlb{nwm-1}
If\/ $\ang{\xi,k,m}\nin\abc\jfi$ then
%\imar{nwm-1}
the extension\/ $\jsi$ of\/ $\jfi$ by\/
$\abc{\jsi}=\abc\jfi\cup\ans{\ang{\xi,k,m}}$, 
$\qh{\jsi}\xi km=0$ and\/
$\qT{\jsi}\xi km{\La}=T$, where\/ $T\in\dpi\xi k$
and\/ $\La$ is the empty string, 
belongs to\/ $\ms\jpi$ and\/ $\jfi\cle\jsi$.

\itlb{nwm3} 
If\/ $\ang{\xi,k,m}\in\abc\jfi$ and a set\/ $D\sq\dpi\xi k$ 
is open dense in\/ $\dpi\xi k$ then there is a \mus\  
%\imar{nwm3}
$\jsi\in\mt\jpi$ such that\/ 
$\abc{\jsi}=\abc\jfi$, $\jfi\cle\jsi$, and\/ 
$\qT\jsi\xi km s\in D$ whenever $s\in 2^{\qh\jsi\xi km}\;.$

\itla{nwm2} 
There is a 2wise disjoint 
\imar{nwm2}
$\jsi\in\mt\jpi$ such that\/ 
$\abc{\jsi}=\abc\jfi$ and\/ $\jfi\cle\jsi\;.$
\een
\ele

\bpf
To prove \ref{nwm3} first use \ref{nwm1} to get a \mus\ 
$\jsi\in\ms\jpi$ with $\jfi\cle\jsi$ and 
$\qh{\jsi}\xi km=h+1$, where $h=\qh\jfi\xi km$.
Then replace each tree $\qT\jsi\xi km s=\jsi(\xi,k,m)(s)$, 
$s\in 2^{h+1}$, with a suitable tree $T'\in D$, 
$T'\sq \qT\jsi\xi km s$.

To prove \ref{nwm2} first apply \ref{nwm1} to get a \mus\ 
$\jsi\in\ms\jpi$ with $\jfi\cle\jsi$, $\abc\jsi=\abc\jfi$, 
and $\qh{\jsi}\xi km=\qh\jfi\xi km+1$ for all 
$\ang{\xi,k,m}\in\abs\jfi$.
Now if $\ang{\xi,k,m}\ne\ang{\et,\ell,n}$ are triples 
in $\abc\vpi$ and $s\in 2^{\qh\jfi\xi km+1}$, 
$t\in 2^{\qh\jfi\et\ell n+1}$, then, by Lemma~\ref{suz}\ref{suz2}, 
there are trees $S\in\dpi\xi k$ and $T\in\dpi\et\ell$ 
satisfying $[S]\cap[T]=\pu$ and $S\sq \qT\jsi\xi km s$, 
$T\sq \qT\jsi\et\ell nt$.
Replace the trees $\qT\jsi\xi km s$, 
$T\sq \qT\jsi\et\ell nt$ with resp.\ $S$, $T$.
Iterate this shrinking construction for all triples 
$\ang{\xi,k,m}\ne\ang{\et,\ell,n}$ and strings $s,t$ as above.
\epf

\gla{Refinements}
\las{sek2}

Here we consider \rit{refinements} of perfect tree forcings 
and \muf s, the key technical tool of definition of various 
forcing notions in this paper.

\parf{Refining perfect tree forcings} 
\las{em}  

If $T\in\pet$ (a perfect tree) and $D\sq\pet$ then 
$X\sqf\bigcup D$ will mean 
\index{zzsqf@$\sqf$}%
that there is a finite set $D'\sq D$ such that 
$T\sq\bigcup D'$, or equivalently $[T]\sq\bigcup_{S\in D'}[S]$.

\bdf
\lam{fm}
Let $\dP,\dQ\in\ptf$ be perfect tree forcing 
notions. 
Say that $\dQ$ is a \rit{refinement} of $\dP$ 
%\index{perfect tree forcing, $\ptf$!refinement, $\dP\ssq\dQ$}%
\index{perfect-tree forcing, $\ptf$!refinement, $\dP\ssq\dQ$}
\index{zzPssqQ@$\dP\ssq\dQ$}%
(symbolically $\dP\ssq\dQ$) if 
\ben
\nenu
\itlb{fm1} 
the set 
$\dQ$ is dense in $\dP\cup\dQ$: 
%\imar{fm1}
if $T\in\dP$ then $\sus Q\in\dQ\,(Q\sq T)$;

%\itla{fm2}
%if trees $S,T\in\dP$ are compatible in $\dP\cup\dQ$ 
%then $S,T$ are compatible  in $\dP$;

\itlb{fm3}
if $Q\in\dQ$ 
%\imar{fm3}
then $Q\sqf\bigcup\dP$;

\itlb{fm4}
if $Q\in\dQ$ and $T\in\dP$ then $[Q]\cap[T]$ is 
%\imar{fm4}
clopen in $[Q]$ and $T\not\sq Q$.\qed 
%--- in particular, it follows that $\dP\cap\dQ=\pu$.
\een
\eDf

\vyk{
\ble
\lam{ffin}
If\/ $\dP\ssq\dQ$ and\/ $S\in\dP$, $T\in\dQ$, 
then there are finite sets\/ $A,B\sq\dQ$ such that\/ 
$[T]\cap[S]=\bigcup_{Q\in A}[Q]$  
and\/ $[T]\bez[S]=\bigcup_{Q\in B}[Q]$ .
\ele
\bpf
By \ref{fm4} $[T]\cap[S]$ is clopen in $[T]$, therefore 
there is a finite set $\sg\sq T$ such that 
$[T]\cap[S]=\bigcup_{s\in\sg}T\ret s$. 
But all trees $T\ret s$ belong to $\dQ\in\ptf$. 
\epf
}

\vyk{
\ble
\lam{mea}
If\/ $\dP\ssq\dQ$ and\/ $S\in\dP$, $T\in\dQ$, then\/ 
$[S]\cap[T]$ is meager in\/ $[S]$. 
Therefore\/ $\dP\cap\dQ=\pu$.
\ele
\bpf
Otherwise there is a string $u\in S$ such that 
$S\ret u\sq[T]\cap[S]$. 
But $S\ret u\in\dP$, which contradicts to \ref{fm4}. 
\epf
}

\ble
\label{pqr}
\ben
\renu
\itlb{pqr0}
If\/ $\dP\ssq\dQ$ and\/ $S\in\dP$, $T\in\dQ$, then\/ 
$[S]\cap[T]$ is meager in\/ $[S]$, 
therefore\/ $\dP\cap\dQ=\pu$ and\/ $\dQ$ is 
open dense in\/ $\dP\cup\dQ\;;$.

\itlb{pqr1}
if\/ $\dP\ssq\dQ\ssq\dR$ then\/ $\dP\ssq\dR$,
thus\/ $\ssq$ is a strict partial order$;$  
%\imar{pqr}

\itlb{pqr2}
if\/ $\sis{\dP_\al}{\al<\la}$  
%\imar{pqr2}
is a\/ \dd\ssq increasing sequence in\/ \ptf\ and\/  
$0<\mu<\la$ then\/
$\dP=\bigcup_{\al<\mu}\dP_\al\ssq
\dQ=\bigcup_{\mu\le\al<\la}\dP_\al\;;$

\itlb{pqr3}
if\/ $\sis{\dP_\al}{\al<\la}$  
%\imar{pqr3}
is a\/ \dd\ssq increasing sequence in\/ \ptf\ and each\/  
$\dP_\al$ is special then\/ $\dP=\bigcup_{\al<\la}\dP_\al$ 
is a regular forcing in\/ $\ptf\;;$

\itlb{pqr4}
in\/ \ref{pqr3}, each\/ $\dP_\ga$ is pre-dense in\/ 
$\dP=\bigcup_{\al<\la}\dP_\al$.
\een                                             
\ele
\bpf
\ref{pqr0}
Otherwise there is a string $u\in S$ such that 
$S\ret u\sq[T]\cap[S]$. 
But $S\ret u\in\dP$, which contradicts to 
\ref{fm}\ref{fm4}. 

\ref{pqr1}, \ref{pqr2}
Make use of \ref{pqr0} to establish \ref{fm}\ref{fm4}.

\ref{pqr3}
To check the regularity let 
$S\in\dP_\al$, $T\in\dP_\ba$, $\al\le\ba$. 
If $\al=\ba$ then, as $\dP_\al$ is special, the trees 
$S,T$ are either \evd\ or \dd\sq comparable by
Lemma~\ref{regfn}.
If $\al<\ba$ then $[S]\cap[T]$ is clopen in $[T]$
by \ref{fm}\ref{fm4}. 

\ref{pqr4}
Let $S\in\dP_\al$, $\al\ne\ga$. 
If $\al<\ga$ then by \ref{fm}\ref{fm1} there is a tree 
$T\in\dP_\ga$, $T\sq S$. 
Now let $\ga<\al$. 
Then $S\sqf\bigcup\dP_\ga$ by \ref{fm}\ref{fm3}, in particular, 
there is a tree $T\in\dP_\ga$ such that $[S]\cap[T]\ne\pu$. 
However $[S]\cap[T]$ is clopen in $[S]$ by \ref{fm}\ref{fm4}. 
Therefore $S\ret u\sq T$ for a string $u\in S$.
Finally $S\ret u\in \dP_\al$ since $\dP_\al\in\ptf$.
\epf

\vyk{
\ble
\lam{resp}
Suppose that\/ $\dP\ssq\dQ\ssq\dR$ are special forcings 
in\/ $\ptf$ 
{\rm(see Definition~\ref{ptf2})} 
and\/ $T\in\dP$. 
Then the sets\/ 
\bce
$ 
\dP'=\ens{P\in\dP}{P\sq T}$, 
$\dQ'=\ens{Q\in\dQ}{Q\sq T}$,
$\dR'=\ens{R\in\dR}{R\sq T}$  
\ece
are special forcings 
in\/ $\ptf$ as well, and still\/ $\dP'\ssq\dQ'\ssq\dR'$.
\ele
\bpf
Let $\dA=\baz\dP$, the base of $\dP$ (Definition~\ref{ptf2}). 
Let $\dA'$ consist of all trees $P\ret s$, where $s\in P\in \dA$ 
and $P\ret s\sq T$, 
and there is no strictly shorter string $t\su s$ with 
$P\ret t\sq T$. 
Then $\dA'$ is an antichain and the base of $\dP'$. 

The sets $\dQ',\dR'$ are special forcings by the same argument. 

To prove $\dP'\ssq\dQ'$ check \ref{fm}\ref{fm3}.
Let $Q\in\dQ'$. 
Then $Q\sqf\bigcup\dP$, in other words, there is a finite 
set $P'\sq\dP$ such that $Q\sq\bigcup P'$. 
As $\dP$ is special, any set $[T]\cap[S]$, $S\in P'$, 
is clopen in $[S]$, therefore $[T]\cap[S]$ is a finite 
%disjoint 
union of sets of the form $S\ret u\in\dP$, $u\in S$.
\epf
}

\vyk{
\ble
\lam{expq}
Suppose that\/ $\la<\omi$, and\/ 
$\sis{\dP_\al}{\al<\la}$ is an\/ 
\dd\ssq increasing sequence of special forcings 
in\/ $\ptf$. 
Then each\/ $\dP_\ga$ is pre-dense in\/ 
$\bigcup_{\al<\la}\dP_\al$.
\ele
\bpf
Let $S\in\dP_\al$, $\al\ne\ga$. 
If $\al<\ga$ then by \ref{fm}\ref{fm1} there is a tree 
$T\in\dP_\ga$, $T\sq S$. 
Now let $\ga<\al$. 
Then $S\sqf\bigcup\dP_\ga$ by \ref{fm3}, in particular, 
there is a tree $T\in\dP_\ga$ such that $[S]\cap[T]\ne\pu$. 
However $[S]\cap[T]$ is clopen in $[S]$ by \ref{fm4}. 
Therefore $S\ret u\sq T$ for a string $u\in S$.
Finally $S\ret u\in \dP_\al$ since $\dP_\al\in\ptf$.
\epf
}

Note that if $\dP,\dQ\in\ptf$ and $\dP\ssq\dQ$ then a 
dense set $D\sq \dP$ is not necessarily dense or even 
pre-dense in $\dP\cup\dQ$. 
Yet there is a special type of refinement which preserves 
at least pre-density.
We modify the relation $\ssq$ as follows. 

\bdf
\lam{ssqm}
Let $\dP,\dQ\in\ptf$ and 
%satisfy $\dP\ssq\dQ$. Let 
$D\sq \dP$. 
% (usually a dense or pre-dense set).  
Say that $\dQ$ \rit{locks\/ $D$ over\/ $\dP$}, 
\index{refinement!locks}%
\index{locks}%
\index{refinement1D@refinement, $\dP\ssa D\dQ$}
\index{perfect-tree forcing, $\ptf$!refinement1D@refinement, 
$\dP\ssa D\dQ$}%
\index{zzPssq1DQ@$\dP\ssa D\dQ$}%
symbolically $\dP\ssa D\dQ$, if $\dP\ssq\dQ$ holds and  
\imar{ssa D}
every tree $S\in\dQ$ satisfies $S\sqf \bigcup D$.
Then simply $\dP\ssq\dQ$ is equivalent to $\dP\ssa\dP\dQ$. 
\edf

As we'll see now, a locked set has to be pre-dense both before 
and after the refinement. 
The additional importance of locking refinements lies in fact that, 
once established, it preserves under further simple 
refinements, that is, $\ssa D$ is transitive in 
a combination with $\ssq$ in the sense of \ref{pqm1} 
of the following lemma:

\ble
\label{pqm}
\ben
\renu
\itlb{pqm0}
If\/ $\dP\ssa D\dQ$ then\/ $D$ is pre-dense in\/ $\dP\cup\dQ$, 
and if in addition\/ $\dP$ is regular then\/ 
$D$ is pre-dense in\/ $\dP$ as well$;$
\imar{pqm}

\itlb{pqm1}
if\/ $\dP\ssa D\dQ\ssq\dR$ 
{\rm(note: the second $\ssq$ is not $\ssa D$!)} 
then\/ $\dP\ssa D\dR\;;$ 
%\imar{pqm1}

\itlb{pqm2}
if\/ $\sis{\dP_\al}{\al<\la}$  
%\imar{pqm2}
is a\/ \dd\ssq increasing sequence in\/ \ptf,  
$0<\mu<\la$, and\/ 
$\dP=\bigcup_{\al<\mu}\dP_\al\ssa D\dP_\mu$, then\/
$\dP\ssa D\dQ=\bigcup_{\mu\le\al<\la}\dP_\al\;.$
\een                                             
\ele
\bpf
\ref{pqm0}
To see that $D$ is pre-dense in $\dP\cup\dQ$, 
let $T_0\in\dP\cup\dQ$. 
By \ref{fm}\ref{fm1}, there is a tree $T\in\dQ$, $T\sq T_0$. 
Then $T\sqf \bigcup D$, in particular, there is a tree 
$S\in D$ with $X=[S]\cap[T]\ne\pu$.  
However $X$ is clopen in $[T]$ by \ref{fm}\ref{fm4}.
Therefore, by Lemma~\ref{regfn11},
there is a tree $T'\in\dQ$ 
with $[T']\sq X$, thus $T'\sq S\in D$ and $T'\sq T\sq T_0$. 
We conclude that $T_0$ is compatible with $S\in D$ in $\dP\cup\dQ$.

To see that $D$ is pre-dense in $\dP$ (assuming $\dP$ is regular), 
let $S_0\in\dP$. 
It follows from the above that $S_0$ is compatible with some 
$S\in D$, hence, $S$ and $S_0$ are not absolutely incompatible. 
It remains to use Lemma~\ref{regfn}\ref{regfn3}.  

To prove \ref{pqm1} on the top of Lemma~\ref{pqr}\ref{pqr1}, 
let $R\in\dR$. 
Then $R\sqf\bigcup\dQ$, but each $T\in\dQ$ satisfies 
$T\sqf\bigcup D$.
The same for \ref{pqm2}.
\epf

The existence of \dd{\ssa D}refinements will be established below.

\parf{Refining \muf s} 
\las{emf1}

Let $\jpi,\jqo$ be \muf s. 
% of sizes resp.\ $\vt\le\za$. 
Say that $\jqo$ is an \rit{refinement} of $\jpi$, 
\index{multiforcing!refinement, $\jpi\bssq\jqo$}%
\index{zzpissqqo@$\jpi\bssq\jqo$}%
\kmar{bssq}
symbolically $\jpi\bssq\jqo$, if $\abc\jpi\sq\abc\jqo$ 
and 
$\dpi\xi k\ssq\dqo\xi k$ whenever $\ang{\xi,k}\in\abc\jpi$.

\vyk{
\bcor
[of Lemma~\ref{pqr}]
\label{pqrC}
\ben
\renu
\itla{pqrC1}
if\/ $\jpi\bssq\jqo\bssq\jro$ then\/ $\jpi\bssq\jro\;;$ 
\imar{pqrC}

\itlb{pqrC2}
if\/ $\sis{\nor\jpi\al}{\al<\la}$ is a\/ \dd\bssq increasing 
sequence of\/ \muf s then the componentwise union\/ 
$\bkw_{\al<\la}\nor\jpi\al$ is 
a \mre\ \muf, and if in addition\/ $0<\mu<\la$, then\/
$\bkw_{\al<\mu}\nor\jpi\al\bssq
\bkw_{\mu\le\al<\la}\nor\jqo\ga$.\qed
\een                                             
\ecor
}

\bcor
[of Lemma~\ref{pqr}]
\lam{pqrC}
If\/ $\jpi\bssq\jqo\bssq\jro$ then\/ $\jpi\bssq\jro$. 

If\/ $\jpi\bssq\jqo$ then the \muf\/ $\mt\jqo$ 
is open dense in\/ $\mt{\jpi\kw\jqo}$. 
\qed
\ecor

Our next goal is to introduce a version of 
Definition~\ref{ssqm} suitable for \muf s; 
we expect an appropriate version of 
Lemma~\ref{pqm} to hold. 

First of all, we accomodate the definition of the 
relation $\sqf$ in Section~\ref{em} for \mut s. 
Namely if $\ju$ is a \mut\ and $\zD$ a collection of \mut s, 
then $\ju\sqf\bigvee \zD$ 
\kmar{sqf}%
will mean that there is a finite 
\index{zzsqf@$\sqf$}%
\index{zzsqfV@$\sqf\bigvee$}%
set $\zD'\sq \zD$ satisfying 
1) $\abc \jv=\abc \ju$ for all $\jv\in \zD'$, and 
2) $[\ju]\sq\bigcup_{\jv\in \zD'}[\jv]$.  

\bdf
\lam{ssl}
Let\/ $\jpi,\jqo$ be \muf s, and\/ $\jpi\bssq\jqo$. 
Say that\/ 
\rit{$\jqo$ locks a set\/ $\zD\sq\mt\jpi$ over\/ $\jpi$}, 
\index{locks}% 
symbolically\/ $\jpi\ssb\zD\jqo$
\kmar{ssb zD}
%(the upper index {\sf mf} accounts for the `miltiforcing version'), 
\index{multiforcing!refinement2D@refinement, $\jpi\ssb\zD\jqo$}%
\index{refinement2D@refinement, $\jpi\ssb\zD\jqo$}%
\index{zzpissq2Dqo@$\jpi\ssb\zD\jqo$}%
if the following condition holds:
\ben
\fenu
\itlb{ssl*}
if\/ 
%\imar{ssl*}
$\zp\in\mt\jpi$, $\ju\in\mt\jqo$, $\abc\ju\sq\abc\jpi$, 
$\abc\ju\cap\abc\zp=\pu$, 
then there is\/ $\zq\in\mt\jpi$  
such that\/ $\zq\leq\zp$, still\/ $\abc{\zq}\cap\abc\ju=\pu$,  
%$\pro\zr\jqo=\pro\ju\jqo$, 
and\/ $\ju\sqf \bigvee \duz \zD\ju\zq$, where\/ 
$$
\duz \zD\ju\zq=\ens{\ju'\in\mt\jpi}
{\abc{\ju'}=\abc\ju\:\text{ and }\:\ju'\cup\zq\in \zD}\,. 
\index{zzDuq@$\duz D\ju\zq$}%
\eqno\qed
$$
\een
\eDf

Note that if $\zp,\ju,\zD,\zq$ are as indicated then still 
$\ju\cup\zq\sqf \bigvee\zD$ holds via the finite set
$\zD'=\ens{\ju'\cup\zq}{\ju'\in\duz \zD\ju\zq}\sq\zD$. 
Anyway the definition of $\ssb\zD$ in \ref{ssl}
looks somewhat different and more complex 
then the definition of $\ssa D$ in \ref{ssqm}, 
which reflects the fact 
that finite-support products of
forcing notions in \ptf\ behave 
differently (and in more complex way) 
than single perfect-tree forcings.
Accordingly, the next lemma, similar to Lemma~\ref{pqm}, is 
way harder to prove.

\ble
\lam{pqn}
Let\/ $\jpi,\jqo,\jsg$ be \muf s and\/
%, $\jpi$ is regular, 
$\zD\sq\mt\jpi$.
Then$:$
\ben
\renu
\itlb{pqn0}
if\/ $\jpi\ssb \zD\jqo$ then\/ $\zD$ is 
%\imar{pqn0}
dense in\/ $\mt\jpi$ and
pre-dense in\/ $\mt{\jpq}\;;$

\itlb{pqn0a}
if\/ $\jpi\ssb \zD\jqo$ and\/ $\zD\sq\zD'\sq\mt\jpi$   
%\imar{pqn0a}
then\/ $\jpi\ssb{\zD'}\jqo\;;$

\itlb{pqn0b}
if\/ $\jpi$ is regular, 
$\jpi\ssb{\zD_i}\jqo$ for\/ $i=1,\dots,n$, all sets\/ 
$\zD_i\sq\mt\jpi$ 
are open dense in\/ $\mt\jpi$, and\/ $\zD=\bigcap_i\zD_i$,   
%\imar{pqn0b}
then\/ $\jpi\ssb{\zD}\jqo\;;$

\itlb{pqn1}
if\/ ${\zD}$ is open dense in\/ $\mt\jpi$ 
and\/ $\jpi\ssb {\zD}\jqo\bssq\jsg$ then\/ $\jpi\ssb {\zD}\jsg\;;$ 
%\imar{pqn1}

\itlb{pqn2}
if\/ $\sis{\jpi_\al}{\al<\la}$  
%\imar{pqn2}
is a\/ \dd\bssq increasing sequence in\/ \mfp, 
$0<\mu<\la$, $\jpi=\bkw_{\al<\mu}\jpi_\al$, 
${\zD}$ is open dense in\/ $\mt\jpi$, 
and\/ $\jpi\ssb {\zD}\jpi_\mu$, then\/
$\jpi\ssb {\zD}\jqo=\bkw_{\mu\le\al<\la}\jpi_\al$.
\een                                             
\ele
\bpf
\ref{pqn0}
To check that ${\zD}$ is pre-dense in $\mt\jpq$, 
let $\zr\in \mt\jpq$. 
Due to the product character of $\mt\jpq$, we can assume 
that $\abc\zr\sq\abc\jpi$. 
%; otherwise just restrict $\zr$ to $\abc\zr\cap\abc\jpi$. 
Let  
$$
X=\ens{\ang{\xi,k}\in\abc\zr}{\zd\zr \xi k\in\mt\jqo}\,, 
\quad 
Y=\ens{\ang{\xi,k}\in\abc\zr}{\zd\zr \xi k\in\mt\jpi}\,. 
$$
Then $\zr=\ju\cup \zp$, where 
$\ju=\zr\res X\in\mt\jqo$, $\zp=\zr\res Y\in\mt\jpi$. 
As $\jqo$ locks ${\zD}$, there is a \mut\ $\zq\in\mt\jpi$  
such that $\zq\leq\zp$, $\abc{\zq}\cap\abc\ju=\pu$,  
and $\ju\sqf \bigcup \duz {\zD}\ju\zq$.
By an easy argument, 
there is a \mut\ $\ju'\in \duz {\zD}\ju\zq$ compatible 
with $\ju$ in $\mt\jqo$; let $\jw\in\mt\jqo$, $\jw\leq\ju$, 
$\jw\leq\ju'$, $\abc\jw=\abc{\ju'}=\abc\ju$. 
Then the \mut\ $\zr'=\jw\cup\zq\in\mt\bbpr$ 
satisfies $\zr'\leq\zr$ and $\zr'\leq\ju'\cup\zq\in {\zD}$. 

To check that ${\zD}$ is dense in $\mt\jpi$, 
suppose that $\zp\in\mt\jpi$. 
Let $\ju=\jLa$ (the empty \mut) 
in \ref{ssl*} of Definition~\ref{ssl}, 
so that $\abc\ju=\pu$ and $\duz {\zD}\ju\zq={\zD}$.\vom

\ref{pqn0a} is obvious. 
To prove \ref{pqn0b}, let $\zp\in\mt\jpi$, 
$\ju\in\mt\jqo$, $\abc\ju\sq\abc\jpi$, 
$\abc\ju\cap\abc\zp=\pu$. 
Iterating \ref{ssl*} for $\zD_i$, $i=1,\dots,n$, we find  
a \mut\ $\zq\in\mt\jpi$  
such that\/ $\zq\leq\zp$, $\abc{\zq}\cap\abc\ju=\pu$,  
%$\pro\zr\jqo=\pro\ju\jqo$, 
and\/ $\ju\sqf \bigvee \duz{(\zD_i)}\ju\zq$ for all $i$, where\/ 
$$
\duz{(\zD_i)}\ju\zq=\ens{\ju'\in\mt\jpi}
{\abc{\ju'}=\abc\ju\:\text{ and }\:\ju'\cup\zq\in \zD_i}\,. 
$$
Thus there are finite sets $U_i\sq \duz{(\zD_i)}\ju\zq$ 
such that $[\ju]\sq\bigcup_{\jv\in U_i}[\jv]$ for all $i$. 
Using the regularity assumption and Lemma~\ref{regfm},
we refine \mut s in $\bigcup_iU_i$,
getting a finite set $W\sq \mt\jpi$  
such that still $\abc{\jw}=\abc\ju$ for all $\jw\in W$, 
$\bigcap_i\bigcup_{\jv\in U_i}[\jv]=\bigcup_{\jw\in W}[\jw]$, 
and if $i=1,\dots,n$ and $\jw\in W$ then   
$[\jw]\sq [\jv]$ for some $\jv\in U_i$ --- therefore
$\jw\cup\zq\in \zD_i$.
We conclude that if $\jw\in W$ then $\jw\cup\zq\in \zD$,
hence $\jw\in\duz{\zD}\ju\zq$.
Thus $W\sq\duz{\zD}\ju\zq$. 
However $[\ju]\sq\bigcup_{\jw\in W}[\jw]$ by the choice of $W$.
We conclude that $\ju\sqf \bigvee \duz{\zD}\ju\zq$, as required.
\vom

\ref{pqn1}
It follows from Corollary~\ref{pqrC} that $\jpi\bssq\jsg$, 
hence it remains to check that $\jsg$ locks ${\zD}$ over $\jpi$.
Assume that $\ju\in\mt\jsg$,  
%$\abc\ju=\ans{\ang{\xi_1,k_1},\dots,\ang{\xi_n,k_n}}\sq\abc\jpi$, 
$\abc\ju\sq\abc\jpi$, 
$\zp\in\mt\jpi$, $\abc\ju\cap\abc\zp=\pu$.
As $\jqo\bssq\jsg$, there is a finite set $U\sq\mt\jqo$ such that 
$\abc\jv=\abc\ju$ for all $\jv\in U$, and 
$[\ju]\sq\bigcup_{\jv\in U}[\jv]$.
As $\jpi\ssb\zD\jqo$, 
by iterated application of \ref{ssl*} of Definition~\ref{ssl}, 
we get a \mut\ $\zq\in\mt\jpi$  
such that $\zq\leq\zp$, still $\abc{\zq}\cap\abc{\ju}=\pu$,  
and if $\jv\in U$ then $\jv\sqf\bigvee\duz {\zD}{\ju}\zq$, where\/ 
$$
\duz{\zD}{\ju}\zq=\ens{\jv'\in\mt\jpi}
{\abc{\jv'}=\abc{\jv}=\abc{\ju}\land{\jv'}\cup\zq\in {\zD}}\,.
$$
Note finally that $\ju\sqf\bigvee U$ by construction, hence 
$\ju\sqf\bigvee\duz {\zD}{\ju}\zq$ as well.\vom

\vyk{
Then each tree $T_i=\zd\ju{\xi_i}{k_i}$,
$i=1,\dots,n$, belongs to $\jsg(\xi_i,k_i)$, 
and hence $T_i\sqf\bigcup\jqo(\xi_i,k_i)$. 
It follows that there are finite sets 
$A_i\sq \jqo(\xi_i,k_i)$ satisfying 
$T_i\sq\bigcup A_i$ ($i=1,\dots,n$).
Let $U\sq\mt{\jqo}$ be the (finite) set of all 
\mut s $\jv$ such that 
$\abc\jv=\abc{\ju}$ and $\zd{\jv}{\xi_i}{k_i}\in A_i$ 
for all $i$.
As $\jpi\ssb\zD\jqo$, 
by iterated application of \ref{ssl*} of Definition~\ref{ssl}, 
we get a \mut\ $\zq\in\mt\jpi$  
such that $\zq\leq\zp$, still $\abc{\zq}\cap\abc{\ju}=\pu$,  
and if $\jv\in U$ then $\jv\sqf\bigvee\duz {\zD}{\ju}\zq$, where\/ 
$$
\duz{\zD}{\ju}\zq=\ens{\jv'\in\mt\jpi}
{\abc{\jv'}=\abc{\jv}=\abc{\ju}\land{\jv'}\cup\zq\in {\zD}}\,.
$$
Note finally that $\ju\sqf\bigvee U$ by construction, hence 
$\ju\sqf\bigvee\duz {\zD}{\ju}\zq$ as well.\vom
}

\ref{pqn2}
We have to check that $\jqo$ locks ${\zD}$ over $\jpi$.
Let $\ju\in\mt\jqo$, 
$\abc\ju
%=\ans{\ang{\xi_1,k_1},\dots,\ang{\xi_n,k_n}}
\sq\abc\jpi$, 
$\zp\in\mt\jpi$, $\abc\ju\cap\abc\zp=\pu$.
As above, there is a finite set $U\sq\mt{\jpi_\mu}$ such that 
$\abc\jv=\abc\ju$ for all $\jv\in U$ and 
$[\ju]\sq\bigcup_{\jv\in U}[\jv]$.
\vyk{
Then each term $T_i=\zd\ju{\xi_i}{k_i}$,
$i=1,\dots,n$, belongs to some $\jpi_\al(\xi_i,k_i)$, 
$\mu\le\al<\la$, 
and hence $T_i\sqf\bigcup\jpi_0(\xi_i,k_i)$ and there are 
finite sets 
$A_i\sq\jpi_0(\xi_i,k_i)$ satisfying 
$T_i\sq\bigcup A_i$.
}%
And so on as in the proof of \ref{pqn1}.
\epf

\parf{Generic refinement of a \muf} 
\las{jex}

Here we introduce a construction, 
due to Jensen in its original form, which implies
the existence of refinements 
of forcings and \muf s, of types $\ssa{D}$ and $\ssb{\zD}$. 

\vyk{
From now on and until Section~\ref{emf4} we stick to objects 
$\la<\omi$, a \dd\bssq increasing sequence 
$\sis{\nor\jpi\al}{\al<\la}$ of small \mdi\ \muf s, the 
componentwise union $\jpi=\bkw_{\al<\la}\nor\jpi\al$, 
a set $Z$, at most countable and satisfying 
$\abc\jpi\sq Z\sq\omi\ti\om$, and another countable set $\cM$, 
as in Theorem~\ref{xist}. 
{\ubf We first consider the case $Z=\abc\jpi$.} 
}

\bdf
\lam{dPhi}
1.  
Suppose that $\jpi$ is a small \muf, and $\cM\in\hc$ is 
any set.
(Recall that $\hc$ = all hereditarily countable sets.) 
This is the input.\vom 

2. 
The set $\cmp$ of all sets $X\in\hc$, \dd\in definable 
in $\hc$ by formulas with sets in $\cM$ as parameters, is 
still countable. 
Therefore there exists a \dd\cle increasing sequence  
$\sis{\jfi(j)}{j<\om}$ of \mus s $\jfi(j)\in\ms\jpi$, 
\dd{\cmp}{\it generic\/} in the sense that it intersects 
any set $\Da\sq\ms\jpi$, $\Da\in{\cmp}$, 
dense in $\ms\jpi$. 
(The density means that for any $\jsi\in\ms\jpi$ there is 
a \mus\ $\jfi\in\Da$ with $\jsi\cle\jfi$.)
 
Let us fix any such a \dd{\cmp}generic  sequence 
$\dphi=\sis{\jfi(j)}{j<\om}$. 
\vom

3. 
Suppose that $\ang{\xi,k}\in\abc\jpi$ and $m<\om$. 
In particular, the sequence $\dphi$ intersects 
every (dense by Lemma~\ref{nwm}\ref{nwm1},\ref{nwm-1}) set 
of the form   
$$
\Da_{\xi km h}=\ens{\jfi\in\ms\jpi}
{\qh\jfi\xi km\ge h}\in{\cmp}\,,\quad\text{where }h<\om\,. 
$$
Hence a tree $\tfd\xi km{s} \in\dpi\xi k$ 
can be associated to any  
$s\in\bse,$ such that, for all $j$, if 
$\ang{\xi,k,m}\in\abc{\jfi(j)}$ and $\lh s\le\qh{\jfi(j)}\xi km$ 
then $\ntt{\jfi(j)}{\xi k,m}s = \tfd\xi km{s}$.\vom

4.
Then it follows from Lemma~\ref{infty} that each set 
$ 
\TS
\qfc\xi km
= 
\bigcap_{h}\bigcup_{s\in2^h}\tfd\xi kms
%\eqno(1)
$ 
is a tree in $\pes$ (not necessarily in $\dpi\xi k$), 
\index{tree!QFixikm@$\qfc\xi km$}%
\index{zzQFixikm@$\qfc\xi km$}%
as well as the trees
%$\uf{\xi m}s:=\raw{\ufi{\xi m}}s$, 
%and still by Lemma~\ref{fus}, 
$$
\TS
\qfd\xi kms 
=\bigcap_{n\ge \lh s}\bigcup_{t\in2^n,\:s\sq t}\tfd\xi kmt 
\,,
\index{tree!QFixikms@$\qfd\xi kms$}%
\index{zzQFixikms@$\qfd\xi kms$}%
%\eqno(2)
$$
and obviously $\qfc{\xi} km=\qfd\xi km\La$. 
Let $\dqf\xi k=\ens{\qfd\xi k ms}{m<\om\land s\in\bse}$.  
\index{zzQQFixik@$\dqf\xi k$}%
\vom 
 
5. 
If $\ang{\xi,k}\in\abc\jpi$ then let 
$\jqo(\xi,k)=\dqf\xi k=\ens{\qfd\xi k ms}{m<\om\land s\in\bse}$.\vom 
 
6.  
Finally if $\jqo=\jqo[\dphi]$ is obtained this way from an 
\dd{\cmp}generic sequence $\dphi$ of \mus s in $\ms\jpi$, 
then $\jqo$ is called 
an \dd\cM\rit{generic refinement of $\jpi$}. 
\index{multiforcing!refinement!generic}%
\index{refinement!generic}%
%, symbolically $\jpi\ssm\cM\jqo$.
\edf

\bpro
[by the countability of $\cmp$]
\lam{gee} 
If\/ 
%$\la<\omi$, $\sis{\nor\jpi\al}{\al<\la}$ 
%is a \dd\bssq increasing sequence of small \mdi\ \muf s, 
$\jpi$ is a small \muf\ 
%=\bkw_{\al<\la}\nor\jpi\al$, 
and\/ $\cM\in\hc$ 
%is a countable set, 
then there is an\/ 
\dd\cM generic refinement\/ $\jqo$ of\/ $\jpi$.
\qed
\epro

\bte
\lam{dj}
If\/ $\jpi$ is a small \muf, 
a set\/ $\cM\in\hc$ contains\/ $\jpi$, $\abc\jpi\sq\cM$,
and\/ $\jqo$ is an\/ \dd\cM generic refinement of\/ $\jpi$, 
then$:$
\ben
\renu
\itlb{dj1}
$\jqo$ is a small \mdi\ \muf, $\abc\jqo=\abc\jpi$, 
and\/ $\jpi\bssq\jqo\;;$% 

\itlb{dj3}
%\label{disj2+}
if\/ $\ang{\xi,k}\in\abc\jpi$ and a set\/ $D\in\cM$, 
%\imar{disj2+}
$D\sq\dpi\xi k$ is pre-dense in\/ $\dpi\xi k$ then\/ 
$\dpi\xi k\ssa D\dqo\xi k\;;$

\aenr
\atm\atm

\itlb{dj4}
%\label{disj3}
if\/ $\ang{\xi,k}\in\abc\jpi$, $m<\om$, and\/ 
%\imar{disj3}
$s\in\bse$ then\/ $\qfd\xi kms=\raw{\qfc\xi km}s\;;$ 

\itlb{dj5}
%\label{disj3+}
if\/ $\ang{\xi,k}\in\abc\jpi$, $m<\om$, and\/ 
%\imar{disj3+}
$s\in\bse$ then\/ $\qfd{\xi}k{m}s\sq\tfd{\xi}k{m}s\;;$ 

\itlb{dj6}
if\/ $\ang{\xi,k}\in\abc\jpi$, $m<\om$, 
and strings\/ $t'\ne t$ in\/ $\bse$ 
%\label{disj4}
%\imar{disj4}
are\/ \dd\sq incomparable then\/ 
$[\qfd\xi km{t'}]\cap[\qfd\xi kmt]=[\tfd\xi k m{t'}]
\cap[\tfd\xi kmt]=\pu\;;$ 

\itlb{dj7}
%\label{uu21}
if\/ $\ang{\xi,k,m}\ne\ang{\et,\ell,n}$ then\/ 
%\imar{uu21}
$[\qfc\xi km]\cap[\qfc\et \ell n]=\pu\;;$ 

\itlb{dj8}
%\label{uu22}
if\/ $\ang{\xi,k}\in\abc\jpi$, $S\in\dqo\xi k$ and\/ $T\in\dpi\xi k$ 
%\imar{uu22}
then\/ $[S]\cap[T]$ is clopen in\/ $[S]$ 
and\/ $T\not\sq S$,  
in particular,\/ $\dpi\xi k\cap \dqo\xi k=\pu\;;$  

\itlb{dj9}
%\label{uu23}
if\/ $\ang{\xi,k}\in\abc\jpi$ then 
%\imar{uu23}
the set\/ $\dqo\xi k$ is open dense in\/ $\dqo\xi k\cup\dpi\xi k$. 
\een
If in addition\/ $\jpi=\bkw_{\al<\la}\nor\jpi\al$, where\/
$\la<\omi$, $\sis{\nor\jpi\al}{\al<\la}$ 
is a\/ \dd\bssq increasing sequence of small \mdi\ \muf s, 
and\/ $\cM$ contains\/ $\sis{\nor\jpi\al}{\al<\la}$ 
and all\/ $\al<\la$, then
\ben
\renu
\atc\atc
\itlb{dj2}
if\/ $\al<\la$ then\/ $\nor\jpi\al\bssq\jqo$. 
%\imar{dj2}
%$\jpi\bssq\jqo$ by Corollary~\ref{pqrC}$;$
\een
\ete  

\bpf
Let $\jqo=\jqo[\dphi]$ be obtained from an 
\dd{\cmp}generic sequence $\dphi$ of \mus s in $\ms\jpi$,
as above. 
We argue in the notation of Definition~\ref{dPhi}.

If $\ang{\xi,k}\in\abc\jpi$ and $m<\om$ then by construction 
the system of trees 
$\tfd\xi km{s} \in\dpi\xi k$, $s\in\bse,$
satisfies \ref{sped}\ref{spe2} on the whole domain $s\in\bse.$ 
This leads to \ref{dj4}, \ref{dj5} 
(essentially corollaries of Lemma~\ref{infty}) 
and \ref{dj6}.\vom

To prove \ref{dj7} note that the set $\Da$ 
of all 2wise disjoint \mus s $\jfi$ such that 
$\abc\jfi$ contains both $\ang{\xi,k,m}$ and $\ang{\et,\ell,n}$,
is dense in $\ms\jpi$ by Lemma~\ref{nwm}, and obviously 
$\Da\in{\cmp}.$
Therefore there is $j<\om$ such that $\jfi(j)\in \Da$. 
Let $h=\qh{\jfi(j)}\xi km$ and $h'=\qh{\jfi(j)}\et\ell n$.
Then the sets 
$$
\textstyle
A=\bigcup_{s\in2^h}[\ntt{\jfi(j)}{\xi k,m}s]
=\bigcup_{s\in2^h}[\tfd\xi km{s}],\;
B=\bigcup_{t\in2^{h'}}[\ntt{\jfi(j)}{\xi \ell,n}t]
=\bigcup_{t\in2^{h'}}[\tfd\xi \ell n{t}]
$$ 
are disjoint as $\jfi(j)\in \Da$.
However $[\qfc\xi km] \sq A$ and $[\qfc\et \ell n] \sq B$.

\ref{dj1} 
It follows that the sets $\dqo\xi k=\dqf\xi k$ are 
special \ptf s (Definition~\ref{ptf2}), and hence 
$\jqo$ is a small \mdi\ \muf, as in \ref{dj1}, and 
$\abc\jqo=\abc\jpi$. 
\vom 

%The following lemma presents rather simple consequences 
%of genericity of the background sequence of \mus s 
%$\dphi=\sis{\vpi(j)}{j<\om}$ in Definition~\ref{dPhi}.

\ref{dj8}
To prove the clopenness claim, note that 
the set $\Da$ of all \mus s $\jfi\in\ms\jpi$ such that 
$\ang{\xi,k,m}\in\abc\vpi$ and if $s\in2^h,$ where 
$h=\qh{\jfi}\xi km$, then either 
$\qT\jfi\xi kms\sq T$ or $[\qT\jfi\xi kms]\cap [T]=\pu$, 
is dense. 
To prove $T\not\sq S$,  
the set $\Da'$ of all \mus s $\jfi\in\ms\jpi$ such that 
$\ang{\xi,k,m}\in\abc\vpi$ and 
$T\not\sq \bigcup_{s\in2^h}\qT\jfi\xi kms$, 
where $h=\qh{\jfi}\xi km$, is dense.
Note that $\Da,\Da'\in{\cmp}$ and argue as above.\vom

\ref{dj9}
\rit{Density}. 
If $T\in\dpi\xi k$ then 
the set $\Da(T)$ of all \mus s   
$\jfi\in\ms\jpi$,
such that $\qT\jfi\xi km\La=T$ for some $m$, 
is dense in $\ms\jpi$ by Lemma~\ref{nwm}\ref{nwm-1}, 
therefore $\vpi(j)\in \Da(T)$ for some $j$. 
Then $\tfd\xi km\La=T$ for some $ m<\om$. 
However $\qfd\xi mk\La\sq \tfd\xi km\La$. 
\rit{Openness}. 
Suppose that $S\in\dqo\xi k$, $T\in\dqo\xi k\cup\dpi\xi k$, 
$T\sq S$. 
Then $T\nin\dpi\xi k$ by \ref{dj8}. 
Therefore $T\in\dqo\xi k$.\vom
\vyk{
It remains to establish $\jpi\bssq\jqo$ in \ref{disj1}, that 
is, $\dpi\xi k\ssq\dqo\xi k$ for all $\ang{\xi,k}\in\abc\jpi$.
We check items \ref{fm1}, \ref{fm3}, \ref{fm4} 
in Section~\ref{em}. 
\ref{fm1} follows from \ref{uu23}.
\ref{fm4} follows from \ref{uu22}.
}%

\ref{dj1}, continuation.
To establish $\jpi\bssq\jqo$, 
let $\ang{\xi,k}\in\abc{\jpi}$.
We have to prove that $\jpi(\xi,k)\ssq\jqo(\xi,k)$. 
This comes down to conditions \ref{fm1}, \ref{fm3}, \ref{fm4} 
of Definition~\ref{fm}, of which \ref{fm1} follows from 
\ref{dj9} and \ref{fm4} from \ref{dj8}, and \ref{fm3} 
is obvious since 
$\qfd\xi kms\sq \tfd{\xi}k{m}s \in \jpi(\xi,k)$ 
for all $m$.\vom

\ref{dj3} 
As $\jpi\bssq\jqo$ has been checked, it remains to prove 
$\qfc\xi km\sqf\bigcup D$ for all $m$. 
It follows from the pre-density of $D$ that the set 
$D'=\ens{T\in\jpi(\xi,k)}{\sus S\in D(T\sq S)}$ 
is open dense in $\jpi(\xi,k)$, and still $D'\in{\cmp}$.
%
% in $\jpi(\xi,k)$, hence we can 
%argue as in the proof of \ref{dj2}, with $D$ instead of 
%$\nor\jpi\al(\xi,k)$.
%
%By definition $\sis{\nor\jpi\ga(\xi,k)}{\al\le\ga<\la}$ is 
%a \dd\bssq increasing sequence of \ptf s.
%Thus $\nor\jpi\al(\xi,k)$ is pre-dense in 
%$\jpi(\xi,k)=\bigcup_{\al\le\ga<\la}\nor\jpi\ga(\xi,k)$ 
%by Lemma~\ref{pqr}\ref{pqr4}. 
%Then $D=\ens{T\in\jpi(\xi,k)}
%{\sus S\in \nor\jpi\al(\xi,k)(T\sq S)}$ 
%is open dense in $\jpi(\xi,k)$.
%It follows that 
Then the set $\Da\in{\cmp}$ of all \mus s 
$\jfi\in\ms\jpi$ such that $\ang{\xi,k,m}\in\abc\jfi$ and 
$\qT\jfi\xi kms\in D$ for all $s\in\qh\jfi\xi km$, 
is dense in $\ms\jpi$ by Lemma~\ref{nwm}\ref{nwm3}. 
Thus $\jfi(j)\in \Da$ for some $j$, which witnesses   
$\qfc\xi km\sqf\bigcup D$.\vom

\ref{dj2} 
We have to prove that $\nor\jpi\al(\xi,k)\ssq\jqo(\xi,k)$
whenever $\ang{\xi,k}\in \abc{\nor\jpi\al}$. 
And as $\jpi(\xi,k)\ssq\jqo(\xi,k)$ has been checked, 
it suffices to prove that 
$\qfc\xi km\sqf\bigcup \nor\jpi\al(\xi,k)$ for all $m$. 
However $D=\nor\jpi\al(\xi,k)$ is pre-dense in 
$\jpi(\xi,k)$ by Lemma~\ref{pqr}\ref{pqr4}, and still 
$D\in{\cmp}$, hence we can refer to 
\ref{dj3}.
\epf

\bcor
\lam{xistt}
In the assumptions of Proposition~\ref{gee}, if\/ 
$\abc\jpi\sq Z\sq\omi\ti\om$ and\/ $Z$ is at most countable 
then there is a small \mdi\ \muf\/ $\jqo$ 
such that\/ $\abc\jqo=Z$ 
%$\nor\jpi\al\bssq\jqo$ for all\/ $\al<\la$, 
and\/ $\jpi\bssq\jqo$.
\ecor
\bpf
If $\abc\jpi=Z$ then let $\cM$ be any countable set containing 
%$\sis{\nor\jpi\al}{\al<\la}$ and satisfying $\la\sq\cM$, 
$\jpi$, pick $\jqo$ by Proposition~\ref{gee}, 
and apply Theorem~\ref{dj}. 
If $\abc\jpi\sneq Z$ then we trivially extend the 
construction by $\dqo\xi k=\dpo$ (see Example~\ref{cloL}) 
for all $\ang{\xi,k}\in Z\bez\abc\jpi$.
\epf

\vyk{
\bte
\lam{dj}
If\/ $\la<\omi$, $\sis{\nor\jpi\al}{\al<\la}$ 
is a \dd\bssq increasing sequence of small \mdi\ \muf s, 
$\jpi=\bkw_{\al<\la}\nor\jpi\al$, 
$\cM\in\hc$ is a countable set  
%
%In the assumptions of Proposition~\ref{gee}, if\/ $\cM$ 
containing\/ $\sis{\nor\jpi\al}{\al<\la}$ 
and all ordinals\/ $\al<\la$, 
and\/ $\jqo$ is an\/ \dd\cM generic refinement of\/ $\jpi$, 
then$:$
\ben
\renu
\itla{dj1}
$\jqo$ is a small \mdi\ \muf\ and\/ $\abc\jqo=\abc\jpi\;;$% 

\itla{dj2}
if\/ $\al<\la$ then\/ $\nor\jpi\al\bssq\jqo$, hence, 
%\imar{disj2}
$\jpi\bssq\jqo$ by Corollary~\ref{pqrC}$;$

\itla{dj3}
%\label{disj2+}
if\/ $\ang{\xi,k}\in\abc\jpi$ and a set\/ $D\in\cM$, 
%\imar{disj2+}
$D\sq\dpi\xi k$ is pre-dense in\/ $\dpi\xi k$ then\/ 
$\dpi\xi k\ssa D\dqo\xi k\;;$

\aenr
\atm\atm\atm

\itla{dj4}
%\label{disj3}
if\/ $\ang{\xi,k}\in\abc\jpi$, $m<\om$, and\/ 
%\imar{disj3}
$s\in\bse$ then\/ $\qfd\xi kms=\raw{\qfc\xi km}s\;;$ 

\itla{dj5}
%\label{disj3+}
if\/ $\ang{\xi,k}\in\abc\jpi$, $m<\om$, and\/ 
%\imar{disj3+}
$s\in\bse$ then\/ $\qfd{\xi}k{m}s\sq\tfd{\xi}k{m}s\;;$ 

\itla{dj6}
If\/ $\ang{\xi,k}\in\abc\jpi$, $m<\om$, 
and strings\/ $t'\ne t$ in\/ $\bse$ 
%\label{disj4}
%\imar{disj4}
are\/ \dd\sq incomparable then\/ 
$[\qfd\xi km{t'}]\cap[\qfd\xi kmt]=[\tfd\xi k m{t'}]
\cap[\tfd\xi kmt]=\pu\;;$ 

\itla{dj7}
%\label{uu21}
If\/ $\ang{\xi,k,m}\ne\ang{\et,\ell,n}$ then\/ 
%\imar{uu21}
$[\qfc\xi km]\cap[\qfc\et \ell n]=\pu\;;$ 

\itla{dj8}
%\label{uu22}
if\/ $\ang{\xi,k}\in\abc\jpi$, $S\in\dqo\xi k$ and\/ $T\in\dpi\xi k$ 
%\imar{uu22}
then\/ $[S]\cap[T]$ is clopen in\/ $[S]$ 
and\/ $T\not\sq S$,  
in particular,\/ $\dpi\xi k\cap \dqo\xi k=\pu\;;$  

\itla{dj9}
%\label{uu23}
if\/ $\ang{\xi,k}\in\abc\jpi$ then 
%\imar{uu23}
the set\/ $\dqo\xi k$ is open dense in\/ $\dqo\xi k\cup\dpi\xi k$. 
\een
\ete  

\bpf
Let $\jqo=\jqo[\dphi]$ be obtained from an 
\dd{\cmp}generic sequence $\dphi$ of \mus s in $\ms\jpi$,
as above. 
We argue in the notation of Definition~\ref{dPhi}.

If $\ang{\xi,k}\in\abc\jpi$ and $m<\om$ then by construction 
the system of trees 
$\tfd\xi km{s} \in\dpi\xi k$, $s\in\bse,$
satisfies \ref{sped}\ref{spe2} on the whole domain $s\in\bse.$ 
This leads to \ref{dj4}, \ref{dj5} 
(essentially corollaries of Lemma~\ref{infty}) 
and \ref{dj6}.\vom

To prove \ref{dj7} note that the set $\Da$ 
of all 2wise disjoint \mus s $\jfi$ such that 
$\abc\jfi$ contains both $\ang{\xi,k,m}$ and $\ang{\et,\ell,n}$,
is dense in $\ms\jpi$ by Lemma~\ref{nwm}, and obviously 
$\Da\in{\cmp}.$
Therefore there is $j<\om$ such that $\jfi(j)\in \Da$. 
Let $h=\qh{\jfi(j)}\xi km$ and $h'=\qh{\jfi(j)}\et\ell n$.
Then the sets 
$$
\textstyle
A=\bigcup_{s\in2^h}[\ntt{\jfi(j)}{\xi k,m}s]
=\bigcup_{s\in2^h}[\tfd\xi km{s}],\;
B=\bigcup_{t\in2^{h'}}[\ntt{\jfi(j)}{\xi \ell,n}t]
=\bigcup_{t\in2^{h'}}[\tfd\xi \ell n{t}]
$$ 
are disjoint as $\jfi(j)\in \Da$.
However $[\qfc\xi km] \sq A$ and $[\qfc\et \ell n] \sq B$.

It follows that the sets $\dqo\xi k=\dqf\xi k$ are 
special \ptf s (Definition~\ref{ptf2}), and hence 
$\jqo$ is a small \mdi\ \muf, as in \ref{dj1}.\vom 

%The following lemma presents rather simple consequences 
%of genericity of the background sequence of \mus s 
%$\dphi=\sis{\vpi(j)}{j<\om}$ in Definition~\ref{dPhi}.

\ref{dj8}
To prove the clopenness claim, note that 
the set $\Da$ of all \mus s $\jfi\in\ms\jpi$ such that 
$\ang{\xi,k,m}\in\abc\vpi$ and if $s\in2^h,$ where 
$h=\qh{\jfi}\xi km$, then either 
$\qT\jfi\xi kms\sq T$ or $[\qT\jfi\xi kms]\cap [T]=\pu$, 
is dense. 
To prove $T\not\sq S$,  
the set $\Da'$ of all \mus s $\jfi\in\ms\jpi$ such that 
$\ang{\xi,k,m}\in\abc\vpi$ and 
$T\not\sq \bigcup_{s\in2^h}\qT\jfi\xi kms$, 
where $h=\qh{\jfi}\xi km$, is dense.
Note that $\Da,\Da'\in{\cmp}$ and argue as above.\vom

\ref{dj9}
\rit{Density}. 
If $T\in\dpi\xi k$ then 
the set $\Da(T)$ of all \mus s   
$\jfi\in\ms\jpi$,
such that $\qT\jfi\xi km\La=T$ for some $m$, 
is dense in $\ms\jpi$ by Lemma~\ref{nwm}\ref{nwm-1}, 
therefore $\vpi(j)\in \Da(T)$ for some $j$. 
Then $\tfd\xi km\La=T$ for some $ m<\om$. 
However $\qfd\xi mk\La\sq \tfd\xi km\La$. 
\rit{Openness}. 
Suppose that $S\in\dqo\xi k$, $T\in\dqo\xi k\cup\dpi\xi k$, 
$T\sq S$. 
Then $T\nin\dpi\xi k$ by \ref{dj8}. 
Therefore $T\in\dqo\xi k$.\vom
\vyk{
It remains to establish $\jpi\bssq\jqo$ in \ref{disj1}, that 
is, $\dpi\xi k\ssq\dqo\xi k$ for all $\ang{\xi,k}\in\abc\jpi$.
We check items \ref{fm1}, \ref{fm3}, \ref{fm4} 
in Section~\ref{em}. 
\ref{fm1} follows from \ref{uu23}.
\ref{fm4} follows from \ref{uu22}.
}%

\ref{dj2}. 
Suppose that $\ang{\xi,k}\in\abc{\nor\jpi\al}$.
We have to prove that $\nor\jpi\al(\xi,k)\ssq\jqo(\xi,k)$. 
This comes down to conditions \ref{fm1}, \ref{fm3}, \ref{fm4} 
of Definition~\ref{fm}, of which \ref{fm1} follows from 
\ref{dj9} and \ref{fm4} from \ref{dj8}. 
It remains to prove \ref{fm3}, 
that is, if $m<\om$ then $\qfc\xi km\sqf\bigcup\nor\jpi\al(\xi,k)$. 
By definition $\sis{\nor\jpi\ga(\xi,k)}{\al\le\ga<\la}$ is 
a \dd\ssq increasing sequence of \ptf s.
Thus $\nor\jpi\al(\xi,k)$ is pre-dense in 
$\jpi(\xi,k)=\bigcup_{\al\le\ga<\la}\nor\jpi\ga(\xi,k)$ 
by Lemma~\ref{pqr}\ref{pqr4}. 
Then 
$D=\ens{T\in\jpi(\xi,k)}
{\sus S\in \nor\jpi\al(\xi,k)(T\sq S)}$ 
is open dense in $\jpi(\xi,k)$.
It follows that the set $\Da\in{\cmp}$ of all \mus s 
$\jfi\in\ms\jpi$ such that $\ang{\xi,k,m}\in\abc\jfi$ and 
$\qT\jfi\xi kms\in D$ for all $s\in\qh\jfi\xi km$, 
is dense in $\ms\jpi$ by Lemma~\ref{nwm}\ref{nwm3}. 
Thus $\jfi(j)\in \Da$ for some $j$, which witnesses   
$\qfc\xi km\sqf\bigcup\nor\jpi\al(\xi,k)$.\vom

Finally, \ref{dj3}. 
The set $D$ is is pre-dense in $\jpi(\xi,k)$, hence we can 
argue as in the proof of \ref{dj2}, with $D$ instead of 
$\nor\jpi\al(\xi,k)$.
\epf

\bcor
\lam{xistt}
In the assumptions of Proposition~\ref{gee}, if\/ 
$\abc\jpi\sq Z\sq\omi\ti\om$ and\/ $Z$ is at most countable 
then there is a small \mdi\ \muf\/ $\jqo$ 
such that\/ $\abc\jqo=Z$ 
%$\nor\jpi\al\bssq\jqo$ for all\/ $\al<\la$, 
and\/ $\jpi\bssq\jqo$.
\ecor
\bpf
If $\abc\jpi=Z$ then let $\cM$ be any countable set containing 
%$\sis{\nor\jpi\al}{\al<\la}$ and satisfying $\la\sq\cM$, 
$\jpi$, pick $\jqo$ by Proposition~\ref{gee}, 
and apply Theorem~\ref{dj}. 
If $\abc\jpi\sneq Z$ then we trivially extend the 
construction by $\dqo\xi k=\dpo$ (see Example~\ref{cloL}) 
for all $\ang{\xi,k}\in Z\bez\abc\jpi$.
\epf 
}

\bcor
\lam{xistC}
Suppose that\/ $\la<\omi$, and\/ 
$\sis{\dP_\al}{\al<\la}$ is an\/ 
\dd\ssq increasing sequence of countable special forcings 
in\/ $\ptf$. 
Then there is a countable special forcing\/ $\dQ\in\ptf$ 
such that\/ $\dP_\al\ssq\dQ$ for each\/ $\al<\la\;.$ 
\ecor
\bpf
If $\al<\la$ then let a \muf\ $\nor\jpi\al$ be defined 
by $\abc{\nor\jpi\al}=\ans{\ang{0,0}}$ and 
by $\dpj00{\nor\jpi\al}=\dP_\al$. 
%Let $\cM$ be a countable transitive model of $\zfcm$ which 
%contains the sequence $\sis{\dP_\al}{\al<\la}$ and is such 
%that both $\la$ and all sets $\dP_\al$ 
%are at most countable in $\cM$. 
By Proposition~\ref{gee} 
and Theorem~\ref{dj} there is a \muf\ $\jqo$ satisfying  
$\abc\jqo=\ans{\ang{0,0}}$ and   
${\nor\jpi\al}\bssq\jqo$, $\kaz\al$. 
Let $\dQ=\dqo00$.
\epf

%\qeD{Theorem~\ref{xist}}

\parf{Preservation of density} 
\las{pres}

This Section proves a special consequence of 
\dd{\cmp}genericity of \muf\ refinements, 
the relation $\ssb{}$ of Definition~\ref{ssl} between
a \muf\ and its refinement.

\bte
\lam{uu4}
In the assumptions of Theorem~\ref{dj}, if\/ 
%$\jqo$ is an\/ \dd\cM generic refinement of\/ $\jpi$,
%$\jpi\ssm\cM\jqo$ 
${\zD}\in\cmp$, ${\zD}\sq\mt\jpi$, and\/ 
${\zD}$ is open dense in $\mt\jpi$,
then $\jpi\ssb {\zD}\jqo$. 
\ete
\bpf
We suppose that $\jqo=\jqo[\dphi]$ is obtained from an
increasing 
\dd{\cmp}generic sequence $\dphi$ of \mus s in $\ms\jpi$,
as in Definition~\ref{dPhi}, and argue in the notation of 
\ref{dPhi}.

Suppose that $\zp\in\mt\jpi$, $\ju\in\mt\jqo$, 
%$\abc\ju\sq\abc\jpi$, 
$\abc\ju\cap\abc\zp=\pu$, as in \ref{ssl*} of 
Definition~\ref{ssl}; the extra condition 
$\abc\ju\sq\abc\jpi$ holds automatically as we still have 
$\abc\jqo=\abc\jpi$.
Let $X=\abc\ju$, $Y=\abc\jpi\bez X.$
If $\ang{\xi,k}\in X$ then 
%by definition we have
$\zd\ju\xi k=\qfd\xi k{m_{\xi k}}{s_{\xi k}}$, 
where $m_{\xi k}<\om$ and $s_{\xi k}\in\bse$. 
By obvious reasons we can assume that $s_{\xi k}=\La$, 
hence $\zd\ju\xi k=\qfc\xi k{m_{\xi k}}$,
for all $\ang{\xi,k}\in X$.

Consider the set $\Da$ of all \mus s   
$\jfi\in \ms\jpi$ such that there is a number $H>0$
and a \mut\ $\zq\in \ms\jpi$ 
satisfying \ref{uu42}, \ref{uu41}, \ref{uu41+}, \ref{uu43} below.
\ben
\nenu
\itlb{uu42}
%$\zp$ occurs in $\jfi$, 
$\abc\zq\cap X=\pu$ and  $\zq\leq\zp$; 
%(note that $\ju\res X'=\pro\ju\jpi\in\mt\jpi$);% 
%\imar{uu42}

\itlb{uu41}
if $\ang{\xi ,k}\in X$ then
%\imar{uu41}
$\ang{\xi,k,m_{\xi k}}\in\abc\jfi$;

\itlb{uu41+}
if $\ang{\xi,k,m}\in\abc\jfi$ then
%\imar{uu41+}
$\qh\jfi\xi km=H$.
\een
To formulate the last requirement, we need one more definition. 
Suppose that $\tau=\sis{t_{\xi k}}{\ang{\xi,k}\in X}$ is a 
system of strings $\tau(\xi,k)=t_{\xi k}\in 2^H$, symbolically 
$\tau\in(2^H){}^X.$ 
Define a \mut\ $\ttt(\jfi,\tau)\in\mt\jpi$ so that 
$\abc{\ttt(\jfi,\tau)}=X$ and 
$\zd{\ttt(\jfi,\tau)}\xi k=\qT\jfi\xi k{m_{\xi k}}{t_{\xi k}}$ 
for all $\ang{\xi,k}\in X$.
Note that $\abc{\ttt(\jfi,\tau)}=\abc\ju$, 
and hence the \mut\ 
$\ttt(\jfi,\tau)\cup\zq$ belongs to $\mt\jpi$ as well.\snos
{Here, if $\zp,\zq$ are \mut s satisfying
$\abc\zp\cap\abc\zq=\pu$ (disjoint domains),
then $\zp\cup\zq$, a \rit{disjoint union}, is a \mut\ such that 
\index{multitree!disjoint union, $\zp\cup\zq$}% 
\index{zzpuq@$\zp\cup\zq$}% 
$\abc{\zp\cup\zq}=\abc\zp\cup\abc\zq$ and 
$\zd{\zp\cup\zq}\xi k=\zd{\zp}\xi k$ whenever $\ang{\xi,k}\in\abc\zp$ 
but
$\zd{\zp\cup\zq}\xi k=\zd{\zq}\xi k$ whenever $\ang{\xi,k}\in\abc\zq$. 
} 
Now goes the last condition.
\ben
\nenu
\atc\atc\atc
\itlb{uu43}
If $\tau\in(2^H){}^X$  
%\imar{uu43}
then $\ttt(\jfi,\tau)\cup \zq\in \zD$.
\een

\ble
\lam{uu4L}
The set\/ $\Da$ is dense in\/ $\ms\jpi$.
\ele 
\bpf[Lemma]
Suppose that $\jsi\in \ms\jpi$; we have to find 
a \mus\ $\jfi\in \ms\jpi$ with $\jsi\cle\jfi$. 
First of all, by Lemma~\ref{nwm}\ref{nwm1}\ref{nwm-1}
we can assume that 
\ben
\aenu 
\itlb{uu4a}
if $\ang{\xi,k}\in X$ then 
%\imar{uu4a}
$\ang{\xi,k,m_{\xi k}}\in\abc\jsi$;

\itlb{uu4b}
there is a number $g>0$
such that $\qh\jsi\xi km=g$ for all $\ang{\xi,k,m}\in\abc\jsi$.
\een
Let $H=g+1$. 
Define 
%a \mus\ 
$\jhi\in \ms\jpi$ so that 
$\abc\jhi=\abc\jsi$, and 
$\qh\jhi\xi km=H$,  
$\qT\jhi\xi km {s\we i}=\raw{\qT\jsi\xi km s}i$ 
for all $\ang{\xi,k,m}\in\abc\jsi$ and $s\we i\in2^H$; 
then $\jsi\cle\jhi$.

It follows from the open density of $\zD$ that there is a \mut\ 
$\zq\in\mt\jpi$ satisfying \ref{uu42}, and a \mus\
$\jfi\in \ms\jpi$ satisfying \ref{uu43} 
and such that still $\abc{\jfi}=\abc\jsi$ 
and $\qh{\jfi}\xi km=H$ for all $\ang{\xi,k,m}\in\abc\jsi$, 
and in addition
\ben
\aenu
\atc\atc
\itlb{hip1}
if $\ang{\xi,k}\in X$ and $s\in2^H$ then 
%\imar{hip1}
$\qT{\jfi}\xi k{m_{\xi k}}s\sq\qT{\jhi}\xi k{m_{\xi k}}s$; 

\itlb{hip2}
$\qT{\jfi}\xi k{m}s=\qT{\jhi}\xi kms$ for all 
%\imar{hip2}
applicable $\xi,k,m,s$ not covered by \ref{hip1}. 
\vyk{
\itla{hip3}
if $\tau=\sis{t_{\xi k}}{\ang{\xi,k}\in X}$ is a 
system of strings $t_{\xi k}\in 2^H$ 
\imar{hip3}
then $\ttt(\jfi,\tau)\cup \zp\in D$.
}%
\een
Namely to achieve \ref{uu43} for one particular  
$\tau\in(2^H){}^X,$ 
consider the \mut\ 
$\zr=\ttt(\jhi,\tau)\cup \zp$. 
%where temporarily $\zp= \ju\res X'\in \mt\jpi$.
There is a \mut\ $\zr'\in \zD$, $\zr'\leq\zr$. 
Let a new \mus\ $\jhi'$ be obtained from $\jhi$ by the reassignment 
$\qT{\jhi'}\xi k{m_{\xi k}}{\tau(\xi, k)}=\zd{\zr'}\xi k$ for 
all $\ang{\xi,k}\in X$.
To get the input for the next step, 
let $\zp'=\zr'\res Y,$\snos
{Here $\zr'\res Y$ is the plain restriction of the function 
$\zr':\abc{\zr'}\to\pet$ to the set $\abc{\zr'}\cap Y$.}  
\index{multitree!restriction, $\zp\res X$}% 
\index{zzpIX@$\zp\res X$}% 
so that 
$\zr'=\ttt(\jhi',\tau)\cup \zp'\in \zD$.

Now consider another $\tau'\in(2^H){}^X$ and the \mut\ 
$\zr'=\ttt(\jhi',\tau')\cup \zp'$. 
There is $\zr''\in \zD$, $\zr''\leq\zr'$. 
Define $\jhi''$ from $\jhi'$ by the reassignment 
$\qT{\jhi''}\xi k{m_{\xi k}}{\tau'(\xi,k)}=\zd{\zr''}\xi k$ for 
all $\ang{\xi,k}\in X$. 
Let $\zp''=\zr''\res Y,$ so that 
$\zr''=\ttt(\jhi'',\tau')\cup \zp''\in \zD$.

And so on.   
The final \mus\ and \mut\ of this construction will be $\jfi$ 
and $\zq$ satisfying \ref{uu42}, \ref{uu41}, \ref{uu41+}, 
\ref{uu43}.
Note that $\jsi\cle\jfi$, 
as we only amend 
the \dd Hth level of $\jhi$ absent in $\jsi$. 
\epF{Lemma}  

Note that $\Da$ is defined in $\hc$ using sets 
$\zD$, $\jpi$, $\zp$, $X$, and the map 
$\ang{\xi,k}\to m_{\xi k}:X\to\om$ as parameters. 
Now, $\zD$, $\jpi$ belong to $\cmp$ straightforwardly, 
$X$ belongs to $\cmp$ since it is a finite subset of a 
set $\abc\jpi\sq\cM$, and $\zp$ belongs to $\cmp$  by 
similar reasons.
It follows that $\Da$ belongs to $\cmp$ as well. 

Therefore, by the lemma and the choice of $\dphi$,  
%we can assume that $\Da$ belongs to the list 
%of \rit{relevant} sets in Definition~\ref{dPhi}(i). 
%Therefore 
there is an index $j$ 
such that the \mus\ $\jfi(j)$ belongs 
to $\Da$, which is witnessed by a number $H>0$ and a \mut\ 
$\zq\in\mt\jpi$  
satisfying \ref{uu42}, \ref{uu41}, \ref{uu41+}, \ref{uu43} 
for $\jfi(j)$ instead of $\jfi$.
To prove that 
$\ju\sqf\bigvee \duz \zD\ju\zq$,
note that   
the \mut s $\ttt(\jfi(j),\tau)\cup \zq$, $\tau\in(2^H){}^X$, 
belong to $\zD$ by \ref{uu43}, and easily 
$[\ju]\sq\bigcup_{\tau\in(2^H){}^X}[\ttt(\jfi(j),\tau)]$.
\epf

\gla{Structure of real names}
\las{sek3}

Here we turn to some details of the structure of reals in 
models of \dd{\mt\jpi}generic type, $\jpi$ being a \muf. 
We are going focus on \rit{non-principal} reals, \ie, 
those different from the principal generic reals 
$x_{\xi k}[G]$ (Remark~\ref{adds}).
We'll work towards the goal of making every non-principal 
real to be non-generic with respect to each of 
the factor forcing notions $\jpi(\xi,k)$.

\parf{Real names} 
\las{rn}

Our next goal is to 
introduce a suitable notation related to names of reals 
in the context of forcing notions of the form  
$\mt\jpi$.

\bdf
\lam{rk'}
A \rit{real name} 
\index{real name}%
is any set $\rc\sq\md\ti(\om\ti\om)$ such that the 
sets $\kkc ni=\ens{\zp\in\md}{\ang{\zp,n,i}\in\rc}$ 
%be a system $\rc=\sis{\kkc ni}{n,i<\om}$ of sets  
\imar{kkc ni}%
\index{zzKcni@$\kkc ni$}%
satisfy the following: 
\ben
\fenu
\itlb{rk1}
if $n,k,\ell<\om$, $k\ne\ell$, and $\zp\in \kkc nk$, 
$\zq\in \kkc n\ell$, then 
$\zp,\zq$ are \emd.
\een
A \qn\ $\rc$ 
\index{real name!small}%
%\index{real name!dyadic}%
is \rit{small} if each set $\kkc ni$ is at most countable ---  
then the sets 
$\dom\rc=\bigcup_{n,i}\kkc ni\sq\md$ and  
$\abs\rc=\bigcup_{n,i}\bigcup_{\zp\in\kkc ni}\abc\zp
\sq\omi\ti\om$, 
\index{zzcII@$\abs\rc$}%
and $\rc$ itself, are countable, too.
\edf

\bdf
\lam{preta}
Let $\rc$ be a \qn\ and
$G\sq\md$ a pairwise compatible set.
Define the
\rit{evaluation} $\rc[G]\in\bn$ so that
\index{real name!evaluation, $\rc[G]$}%
$\rc[G](n)=i$ iff:
\index{zzcG@$\rc[G]$}%
\bit
\item[$-$]
either $\sus \zp\in G\:\sus\zq\in \kkc ni\:(\zp\leq\zq)$
(recall that $\zp\leq\zq$ means $\zp$ is stronger), 

\item[$-$]
or just $i=0$ and
$\neg\;\sus \zp\in G\:\sus\zq\in\bigcup_i\kkc ni\:(\zp\leq\zq)$
(default case).
\qed
\eit
\eDf

%If $\rc$ is a dyadic name then clearly $\rc[G]\in\dn$.

\bdf
\lam{rk''}
Let $\jpi$ be a \muf. 
A \qn\ $\rc$ is said to be a
\rit{\rn\jpi} 
\index{real name!pirealname@\rn\jpi}%
if, in addition to \ref{rk1} above, the following 
condition holds:
\ben
\fenu
\atc
\itlb{rk2}
each set $\kic n=\bigcup_i\kkc ni$ is 
\imar{kic n}%
\index{zzKcn@$\kic n$}%
\rit{pre-dense for\/ $\mt\jpi$}, in the sense that 
the set 
$\kipc n\jpi=\ens{\zp\in\mt\jpi}{\sus\zq\in \kic n\,(\zp\leq\zq)}$ 
\imar{kipc n jpi}%
is dense (then obviously open dense) in $\mt\jpi$.\qed 
\een
\eDf

Generally speaking, we do not assume that 
$\kic n\sq\mt\jpi$. 
However if, in addition to \ref{rk1}, \ref{rk2}  
above, $\kic n\sq\mt\jpi$ holds for all $n$, 
%(and then $\kipc n\jpi=\kic n$)
then say 
that $\rc$ is a \rit{true \rn\jpi}.%  
\index{real name!true@true \rn\jpi} 
\index{real name!pirealname@\rn\jpi!true}%
Then each set $\kic n=\bigcup_i\kkc ni$  
is a pre-dense subset of $\mt\jpi$.

\bre
\lam{ggem}
Let $\jpi$ be a \muf, $\rc$ be a \rn\jpi, and   
a set $G\sq\mt\jpi$ be \dd{\mt\jpi}generic  over 
the collection of all  sets $\kipc n\jpi$ as in
\ref{rk2}
(All of $\kipc n\jpi$ are dense by the choice of $\rc$.)
Then the ``or'' case in Definition~\ref{preta}
never happens as we have $G\cap(\kipc n\jpi)\ne\pu$
by the choice of $G$.
\ere

\bre
\lam{ai}
If $\jpi$ is a \rit{\mre} \muf\ then the notions of being \evd\ 
and being incompatible in $\mt\jpi$ are equivalent by 
Lemma~\ref{regfn}\ref{regfn3}, so that 
a true \rn\jpi\ is the same as a \dd{\mt\jpi}name 
for an element of 
$\bn$ in the general theory of (unramified) forcing.
\ere

\bpri
\lam{proj1}
If $\xi<\omi$, $k<\om$, then $\rpi_{\xi k}$
\index{real name!x.xik@$\rpi_{\xi k}$}%
\index{zzx.xik@$\rpi_{\xi k}$}%
is a \qn\ such that if $i=0,1$ then 
%each 
the set $\jcp ni\xi k=\kk ni{\rpi_{\xi k}}$ 
consists of a lone \mut\ $\zr=\zr^{\xi k}_{ni}$ with 
$\abc{\zr}=\ans{\ang{\xi,k}}$ and 
$\zd\zr \xi k =\ens{t\in\bse}{\lh t\le n\lor t(n)=i}$, 
and if $i\ge2$ then $\jcp ni\xi k=\pu$.
\epri

\bre
\lam{proj2}
If $\jpi\in\md$ and $\ang{\xi,k}\in\abc\jpi$ 
then $\rpi_{\xi k}$ is a \rn\jpi\ of the real 
$x_{\xi k}=x_{\xi k}[G]\in\dn$, 
the $(\xi,k)$th term of a \dd{\mt\jpi}generic sequence 
$\sis{x_{\xi k}[G]}{\ang{\xi,k}\in\abc\jpi}$. 
That is, if $G\sq\mt\jpi$ is generic then 
the real $x_{\xi k}[G]$ defined by \ref{adds} 
coincides with the real $\rpi_{\xi k}[G]$
defined by \ref{preta}.
\ere

\parf{Direct forcing} 
\las{df}

The following definition of the {\ubf direct forcing} 
relation is not 
explicitly associated with 
any concrete forcing notion, but in fact the direct forcing 
relation (in all three instances) is compatible with any 
forcing notion of the form $\mt\jpi$.

Let $\rc$ be a real name. 
Let us say that a multitree $\zp$:
\bit
\item
\rit{directly forces\/ $\rc(n)=i$}, 
\index{directly forces}%
where $n,i<\om$, iff there is a 
\mut\ $\zq\in\kkc ni$ 
such that $\zp\leq \zq$ (meaning: $\zp$ is stronger); 

\item
\rit{directly forces\/ $s\su\rc$},  
where $s\in\nse,$ iff for all $n<\lh s$, $\zp$ 
directly forces $\rc(n)=i$, where $i=s(n)$; 

\vyk{
\item
\rit{directly forces\/ $\rd\ne\rc$}, iff there are strings 
$s,t\in\bse,$ incomparable in $\bse$ and such that  
$\zp$ directly forces $s\su\rc$ 
and $t\su\rd$; 
}

\item
\rit{directly forces\/ $\rc\nin[T]$},  
where $T\in\pet$, iff there is a string $s\in\nse\bez T$ 
such that $\zp$ directly forces $s\su \rc$. 
\eit

\ble
\lam{avo}
If\/ $\jpi$ is a \muf\ and\/ $\rc$ 
is a\/ \rn\jpi, $\zp\in\mt\jpi$, 
$\ang{\xi,k}\in\abc\jpi$, $T\in\pet$, $n<\om$,
then 
\ben
\renu
\itlb{avo1}
there is a number\/ $i<\om$ and  
a \mut\/ $\zq\in\mt\jpi\yd \zq\leq\zp$, 
which directly forces\/ $\rc(n)=i\;;$  

\itlb{avo2}
there is  
a \mut\ $\zq\in\mt\jpi\yd \zq\leq\zp$, 
which directly forces\/ $\rc\nin[\raw T0]$ 
or directly forces\/ $\rc\nin[\raw T1]$. 
\een
\ele

Note that if $T\in\dpi\xi k$ 
then the trees $\raw Ti$, $i=0,1$ belong to 
$\dpi\xi k$.

\bpf
\ref{avo1}
Use the density of sets $\kipc n\jpi$ by \ref{rk''}\ref{rk2} above.

\ref{avo2}
Let $r=\roo{T}$, $n=\lh r$. 
By \ref{avo1}, there is a \mut\ $\zq\in\mt\jpi\yd \zp'\leq\zp$, 
and, for any $m\le n$, --- a number $i_m=0,1$, such that $\zq$ 
directly forces $\rc(m)=i_m$, $\kaz m<n$. 
Let $s\in 2^{n+1}$ be defined by $s(m)=i_m$ for every $m\le n$.
Then $\zq$ directly forces $s\su \rc$. 
On the other hand, $s$ cannot belong to both $\raw T0$ and 
$\raw T1$. 
\epf

\parf{Locking real names} 
\las{lorn}

The next definition extends Definition~\ref{ssl}  
to real names. 

\bdf
\lam{ssl+}
Assume that\/ $\jpi,\jqo$ are \muf s, 
$\rc$ is a real name, and\/ $\jpi\bssq\jqo$. 
Say that\/ 
\rit{$\jqo$ locks\/ $\rc$ over\/ $\jpi$}, 
\index{locks}% 
symbolically\/ $\jpi\ssb\rc\jqo$, 
\kmar{ssb rc}
\index{multiforcing!refinement2c@refinement, $\jpi\ssb\rc\jqo$}%
\index{refinement2c@refinement, $\jpi\ssb\rc\jqo$}%
\index{zzpissq2cqo@$\jpi\ssb\rc\jqo$}%
if $\jqo$ locks, over $\jpi$, each set $\kipc n\jpi$ defined in 
\ref{rk''}\ref{rk2}, in the sense of Definition~\ref{ssl}.
\edf

\bcor
[of Theorem~\ref{uu4}]
\lam{uu4c}
In the assumptions of Theorem~\ref{dj}, if\/ 
${\rc}\in\cmp$ and\/ $\rc$ is a\/ \rn\jpi\ 
then $\jpi\ssb {\rc}\jqo$. 
\ecor
\bpf
Note that each set $\kipc n\jpi$ belongs to $\cmp$
(as so do $\rc$ and $\jpi$) and is dense in $\mt\jpi$, 
so it remains to apply 
Theorem~\ref{uu4}.
\epf

\ble
%[of Lemma~\ref{pqn}]
\lam{pqs}
Let\/ $\jpi,\jqo,\jsg$ be \muf s and\/ $\rc$ be a real name.
Then
\ben
\renu
\itlb{pqs1}
if\/ $\jpi\ssb \rc\jqo$ then\/ $\rc$ is a\/
\rn\jpi\ and   a\/
\rn{(\jpq)}$;$

\itlb{pqs2}
if\/ $\jpi\ssb {\rc}\jqo\bssq\jsg$ then\/ 
$\jpi\ssb {\rc}\jsg\;;$ 
%\imar{pqs2}

\itlb{pqs3}
if\/ $\sis{\jpi_\al}{\al<\la}$  
%\imar{pqs3}
is a\/ \dd\bssq increasing sequence in\/ \mfp, 
$0<\mu<\la$, $\jpi=\bkw_{\al<\mu}\jpi_\al$, 
%${\zD}$ is open dense in\/ $\mt\jpi$, 
and\/ $\jpi\ssb {\rc}\jpi_\mu$, then\/
$\jpi\ssb {\rc}\jqo=\bkw_{\mu\le\al<\la}\jpi_\al\;.$ 
\een
\ele
\bpf
\ref{pqs1} 
By definition, we have $\jpi\ssb{\kipc n\jpi}\jqo$ 
for each $n$, 
therefore $\kipc n\jpi$ is dense in $\mt\jpi$  
(then obviously open dense) 
and pre-dense in $\mt\jpq$ by Lemma~\ref{pqn}\ref{pqn0}.
It follows that $\kipc n{(\jpq)}$ is dense in $\mt\jpq$.

To check \ref{pqs2}, \ref{pqs2} apply 
Lemma~\ref{pqn}\ref{pqn1},\ref{pqn2}.
\epf

\parf{Non-principal names and avoiding refinements} 
\las{npna}

Let $\jpi$ be a \muf. 
Then $\mt\jpi$ adds a collection of reals $x_{\xi k}$, 
$\ang{\xi,k}\in\abc\jpi$, 
where each \rit{principal real} 
$x_{\xi k}=x_{\xi k}[G]$ is \dd{\jpi(\xi,k)}generic 
over the ground set universe. 
Obviously many more reals are added, and given a 
\rn\jpi\ $\rc$, one can elaborate different requirements
for a condition $\zp\in\mt\jpi$ to force that $\rc$ is a name 
of a real of the form $x_{\xi k}$ or to force the opposite. 
But we are mostly interested in simple conditions  
related to the ``opposite'' part. 
The next definition provides such a condition.

\bdf
\lam{npri}
Let $\jpi$ be a \muf, $\ang{\xi,k}\in\abc\jpi$. 
A \qn\  $\rc$ is 
\rit{non-principal over\/ $\jpi$ at\/ $\xi,k$} 
if the following set 
is open dense in\/ $\mt\jpi$:
$$
\ddi\xi k\rc\jpi
\kmar{ddi xi k rc jpi}
=\ens{\zp\in\mt\jpi}{\zp\,\text{ directly forces }\,
\rc\nin[\zd\zp\xi k]}\,.\eqno\qed
$$
\eDf

We'll show below (Theorem~\ref{npn}\ref{npn1}) that the   
non-principality implies $\rc$ being {\ubf not} a name 
of the real $x_{\xi k}[\uG]$. 
And further, the avoidance condition in the next definition 
will be shown to imply $\rc$ being a name 
of a non-generic real.

\bdf
\lam{avod}
Let $\jpi,\jqo$ be \muf s, $\jpi\bssq\jqo$, 
$\ang{\xi,k}\in\abc\jpi$. 
Say that\/ 
\rit{$\jqo$ avoids a\/ \qn\ $\rc$ over\/ $\jpi$ 
at\/ $\xi,k$}, 
\index{locks}% 
in symbol\/ $\jpi\sse\rc\xi k\jqo$, 
\kmar{sse rc xi k}
\index{avoids}%
\index{zzpissqcqo@$\jpi\sse\rc\xi k\jqo$}%
\index{multiforcing!refinementc@refinement, $\jpi\sse\rc\xi k\jqo$}%
\index{refinementc@refinement, $\jpi\sse\rc\xi k\jqo$}%
if for each 
%$\ang{\xi,k}\in\abc\jpi$ and 
$Q\in\dqo\xi k$, $\jqo$ locks, over\/ $\jpi$, the set
$$
\dqu\rc Q\jpi
\kmar{dqu rc Q jpi}
=
\ens{\zr\in\mt\jpi} 
{\zr\,\text{\rm\ directly forces }\, \rc\nin[Q]}\,, 
%\eqno\qed
$$ 
in the sense of Definition~\ref{ssl} ---  
formally $\jpi\ssb{\dqu\rc Q\jpi}\jqo$.
\edf

%The following lemma is a simple corollary of Lemma~\ref{pqn}.

\ble
\lam{pqo}
Assume that\/ $\jpi,\jqo,\jsg$ are \muf s, 
$\ang{\xi,k}\in\abc\jpi$, and\/ 
$\rc$ is a\/ \rn\jpi.
Then$:$
\ben
\renu
\itlb{pqo0}
if\/ $\jpi\sse\rc\xi k\jqo$ and\/ $Q\in\dqo\xi k$ then\/ 
the set\/ $\dqu\rc Q\jpi$ is open dense in\/ $\mt\jpi$ 
and pre-dense in\/ $\mt\jpq\;;$ 
%\imar{pqo0}

\itlb{pqo1}
if\/ $\jpi\sse\rc\xi k\jqo\bssq\jsg$ then\/ 
$\jpi\sse\rc\xi k\jsg\;;$ 
%\imar{pqo1}

\itlb{pqo2}
if\/ $\sis{\jpi_\al}{\al<\la}$  
%\imar{pqo2}
is a\/ \dd\bssq increasing sequence in\/ \mfp, 
$0<\mu<\la$, $\jpi=\bkw_{\al<\mu}\jpi_\al$, 
and\/ $\jpi\sse\rc\xi k\jpi_\mu$, then\/
$\jpi\sse\rc\xi k\jqo=\bkw_{\mu\le\al<\la}\jpi_\al$. 
\een                                             
\ele
\bpf
\ref{pqo0} 
Apply Lemma~\ref{pqn}\ref{pqn0}.

\ref{pqo1} 
Let $\ang{\xi,k}\in\abc\jpi$ and $S\in\dsg\xi k$. 
Then, as $\jqo\bssq\jsg$, there is a finite set 
$\ans{Q_1,\dots,Q_m}\sq \dqo\xi k$ such that 
$S\sq Q_1\cup\dots\cup Q_m$.
We have 
$\jpi\ssb{\dqu\rc {Q_i}\jpi}\jqo$ for all $i$ 
since $\jpi\sse\rc\xi k\jqo$, therefore  
$\jpi\ssb{\dqu\rc {Q_i}\jpi}\jsg$, $\kaz i$, 
by Lemma~\ref{pqn}\ref{pqn1}.
Note that $\bigcap_{i}\dqu\rc {Q_i}\jpi\sq \dqu\rc {S}\jpi$ 
since $S\sq \bigcup_iQ_i$. 
We conclude that $\jpi\ssb{\dqu\rc{S}\jpi}\jsg$ 
by Lemma~\ref{pqn}\ref{pqn0a},\ref{pqn0b}.
% 
%$\jqo$ locks each set of the form $\dqu\rc {Q_i}\jpi$ over\/ $\jpi$. 
%Then 
%$\dqu\rc {Q_i}\jpi$ is open dense in $\mt\jpi$ by \ref{pqo0} 
%for all $i$.
%Lemma~\ref{pqn}\ref{pqn0}, therefore obviously open dense. 
%It follows by Lemma~\ref{pqn}\ref{pqn1} that $\jqo$ locks   
%$\dqu\rc {Q_i}\jpi$ over\/ $\jpi$.
Therefore $\jpi\sse\rc\xi k\jsg$. 
\epf

\parf{Generic refinements avoid non-principal names} 
\las{emf2}

The following theorem says that generic refinements as 
in Section~\ref{jex} avoid nonprincipal names. 
It resembles Theorem~\ref{uu4} to some extent, yet the 
latter is not directly applicable here as both 
the \mut\ $Q$ and  
the set $\dqu\rc Q\jpi$ depend on $\jqo$, and hence the 
sets $\dqu\rc Q\jpi$ do not necessarily belong to $\cmp$. 
However the proof will be based on rather similar
arguments. 

\bte
\lam{K}
In the assumptions of Theorem~\ref{dj}, if\/ 
$\ang{\et,K}\in\abc\jpi\sq\cM$ and\/
%$\jpi\ssm\cM\jqo$ and\/ 
$\rc\in\cM$ is a\/ \rn\jpi\ non-principal 
over\/ $\jpi$ at\/ $\et,K$ then\/ $\jpi\sse\rc\et K\jqo$.  
\ete   
\bpf
Assume that $\jqo=\jqo[\dphi]$ is obtained from an 
\dd{\cmp}generic sequence $\dphi$ of \mus s in $\ms\jpi$,
as in Definition~\ref{dPhi}. 
We stick to the notation of \ref{dPhi}.

Let 
%$\ang{\et,K}\in\abc\jpi$ and 
$Q\in\dqo\et K$; 
we have to prove that $\jqo$ locks the set
$\dqu\rc Q\jpi$ over\/ $\jpi$. 
By construction $Q=\qfd\et K{\xm}{s_0}\sq \qfc\et K{\xm}$ 
for some ${\xm}<\om$; 
it can be assumed that $Q=\qfc\et K{\xm}$. 
Following the proof of Theorem~\ref{uu4}, 
we suppose that $\zp\in\mt\jpi$, $\ju\in\mt\jqo$, 
$\abc\ju\cap\abc\zp=\pu$,  
define $X=\abc\ju$, $Y=\abc\jpi\bez X,$ and 
assume that 
$\zd\ju\xi k=\qfc\xi k{m_{\xi k}}$, 
where $m_{\xi k}<\om$, 
for each $\ang{\xi,k}\in X$. 
%We put $h(\xi,k)=\lh{s_{\xi k}}$ and  
%$h=\tmax_{\ang{\xi,k}\in X}{h(\xi,k)}$. 

Consider the set $\Da$ of all \mus s   
$\jfi\in \ms\jpi$ such that there is a number $H>0$
and a \mut\ $\zq\in \ms\jpi$ satisfying conditions
\ben
\nenu
\itlb{uu62}
%$\zp$ occurs in $\jfi$, 
$\abc\zq\cap X=\pu$ and  $\zq\leq\zp\;;$% 
%\imar{uu62}

\itlb{uu61}
if $\ang{\xi ,k}\in X$ then
%\imar{uu61}
$\ang{\xi,k,m_{\xi k}}\in\abc\jfi$;

\itlb{uu61+}
if $\ang{\xi,k,m}\in\abc\jfi$ then
%\imar{uu61+}
$\qh\jfi\xi km=H$;
\een
(but not \ref{uu43} though)
as in the proof of Theorem~\ref{uu4}, 
along with two more requirements
%$\zp\leqq\zq$ and
\ben
\nenu
\atc\atc\atc\atc 
\itlb{ww45}
$\ang{\et,K,{\xm}}\in \abc\jfi$ ---
hence still $\qh\jfi\et K{\xm}=H$ by \ref{uu41+}; 
%\imar{ww45}

\itlb{ww46}
if $s\in2^H$ and $\tau\in(2^H){}^X$ 
%\imar{ww46}
then $\ttt(\jfi,\tau)\cup \zq$
directly forces $\rc\nin[\qT\jfi{\et}K{{\xm}}s]$. 
\een

\ble
\lam{Ka}
$\Da$ is dense in\/ $\ms\jpi$. 
\ele
\bpf
Suppose that $\jsi\in\ms\jpi$; we can assume that  
$\jsi$ already satisfies 
\ben
\aenu 
\itlb{uu6a}
if $\ang{\xi,k}\in X$ then 
%\imar{uu6a}
$\ang{\xi,k,m_{\xi k}}\in\abc\jsi$;

\itlb{uu6b}
there is a number $g<\om$
such that $\qh\jsi\xi km=g$ for all $\ang{\xi,k,m}\in\abc\jsi$;
\een
as in Lemma~\ref{uu4L}, 
and in addition $\ang{\et,K,{\xm}}\in\abc\jsi$.

Let $H=g+1$. 
Define a \mus\ 
$\jhi\in \ms\jpi$ so that 
$\abc\jhi=\abc\jsi$, and 
$\qh\jhi\xi km=H$, $\qT\jhi\xi km{s\we i}=\raw{\qT\jsi\xi km s}i$ 
for all $\ang{\xi,k,m}\in\abc\jsi$ and $s\we i\in2^H$; 
then $\jsi\cle\jhi$.
%
%Pick strings $t_{\xi k}\in 2^H$ with
%$s_{\xi k}\su t_{\xi k}$ for all
%$\ang{\xi ,k}\in X$. 
We claim that there is a \mut\ 
$\zq\in\mt\jpi$ satisfying \ref{uu62}, and a \mus\
$\jfi\in \ms\jpi$ satisfying \ref{ww46}  
and such that still $\abc{\jfi}=\abc\jsi$ 
and $\qh{\jfi}\xi km=H$ for all $\ang{\xi,k,m}\in\abc\jsi$, 
and in addition

\ben
\aenu
\atc\atc
\itlb{gip1}
if $\ang{\xi,k}\in X$ and $s\in2^H$ then 
%\imar{gip1}
$\qT{\jfi}\xi k{m_{\xi k}}s\sq\qT{\jhi}\xi k{m_{\xi k}}s$, 
and we also have $\qT{\jfi}\et K{\xm}s\sq\qT{\jhi}\et K{\xm}s$; 

\itlb{gip2}
$\qT{\jfi}\xi k{m}s=\qT{\jhi}\xi kms$ for all 
%\imar{gip2}
applicable $\xi,k,m,s$ not covered by \ref{hip1}. 
\een
To achieve \ref{ww46} in one step for one particular  
$\tau\in(2^H){}^X,$  consider the \mut\ 
$\zr=\ttt(\jhi,\tau)\cup \zp$.
By Lemma~\ref{avo} and the density assumption 
of the theorem, 
there is a \mut\ $\zr'\in\mt\jfi$, $\zr'\leq\zr$, which 
directly forces $\rc\nin[\zd{\zr'}\et K]$, and there are 
\mut s $U_s\in\mt\jfi$, $s\in 2^H$, such that 
$U_s\sq \qT\jhi\et K{\xm} s$ and $\zr'$ directly forces 
$\rc\nin[U_s]$, $\kaz s$. 
Let $\jhi'$ be obtained from $\jhi$ by the following 
reassignment. 
\ben
\Renu
\itlb{R1}
We set 
$\qT{\jhi'}\xi k{m_{\xi k}}{\tau(\xi,k)}=\zd{\zr'}\xi k$ 
for all $\ang{\xi,k}\in X$.

\itlb{R2}
If $s\in2^H,$ and either $\ang{\et,K}\nin X$, 
or ${\xm}\ne m_{\et K}$, or $s\ne \tau(\et,K)$ then 
we set 
$\qT{\jhi'}\et K{\xm}s=U_s$.
(Note that if $\ang{\et,K}\in X$ and ${\xm}=m_{\et K}$  
then the tree 
$\qT{\jhi'}\et K{\xm}{\tau(\et,K)}=\zd{\zr'}\et K$ 
is already defined by \ref{R1}.)
\een  
Let $\zp'=\zr'\res Y,$ so that $\zr'=\ttt(\jhi',\tau)\cup\zp'$.
By construction the tree $\zp'$ satisfies \ref{ww46}, 
for the system $\tau$ chosen, 
in the case $\ang{\et,K}\in X$, ${\xm}=m_{\et K}$, $s=\tau(\et,K)$ 
by \ref{R1} and in all other cases by \ref{R2}.

Now consider another   
$\tau'\in(2^H){}^X$ and the \mut\ 
$\zr'=\ttt(\jhi',\tau')\cup \zp'$. 
There is a \mut\ $\zr''\in\mt\jpi$, $\zr''\leq\zr'$, 
which which directly forces $\rc\nin[\zd{\zr'}\et K]$ 
and $\rc\nin[U'_s]$ for each $s\in 2^H,$ where 
$U'_s\in\mt\jfi$ and  
$U'_s\sq \qT{\jhi'}\et K{\xm} s$ . 
Let $\jhi''$ be obtained from $\jhi'$ by the the same  
reassignment (for $\tau'$ instead of $\tau$). 

And so on. 
The final \mus\ and \mut\ of this construction will be $\jfi$ 
and $\zq$ satisfying \ref{uu62}, \ref{uu61}, \ref{uu61+}, 
\ref{ww45}, \ref{ww46}.  
\epF{Lemma}

Come back to the theorem.
Note that $\Da\in{\cmp},$ similarly to the proof of 
Theorem~\ref{uu4}.
Therefore, by the lemma,  there is an index $j$ 
such that the system $\jfi(j)$  belongs to $\Da$. 
Let this be witnessed by a number $H>0$ and a 
\mut\ $\zq\in\mt\jpi$, such that conditions 
\ref{uu62}, \ref{uu61}, \ref{uu61+}, 
\ref{ww45}, \ref{ww46} 
are satisfied for $\jfi=\jfi(j)$. 

It remains to prove that 
$\ju\sqf\bigvee\duz{\dqu\rc Q\jpi}\ju\zq$.
Let $V$ consist of all \mut s  
$\jv=\ttt(\jfi(j),\tau)$, where 
$\tau\in(2^H){}^X$; 
$[\ju]\sq\bigcup_{\jv\in V}[\jv]$ 
by construction.

Further, if $s\in2^H$ and $\jv\in V$ then $\jv\cup\zq$ 
directly forces $\rc\nin[\qT{\jfi(j)}{\et}K{\xm}s]$ by 
\ref{ww46}, that is, directly forces 
$\rc\nin[\tfd{\et}K{\xm}s]$ in the notation of 
Definition~\ref{dPhi}. 
Therefore $\jv\cup\zq$ directly forces 
$\rc\nin[\qfd{\et}K{\xm}s]$ since 
$\qfd{\et}K{\xm}s\sq\tfd{\et}K{\xm}s$ 
by Lemma~\ref{dj}\ref{dj5}. 
However $Q=\qfc{\et}K{\xm}=\bigcup_{s\in2^H}\qfd{\et}K{\xm}s$ 
by Lemma~\ref{dj}\ref{dj4}. 
It follows that $\jv\cup\zq$ directly forces 
$\rc\nin[Q]$, that is, $\jv\in\duz{\dqu\rc Q\jpi}\ju\zq$. 

We conclude that $V$ is a (finite) subset of 
$\duz{\dqu\rc Q\jpi}\ju\zq$. 
And this accomplishes the proof of 
$\ju\sqf\bigvee\duz{\dqu\rc Q\jpi}\ju\zq$.
\epf

\parf{Consequences for reals in generic extensions} 
\las{crn}

We first prove a result saying that all reals in 
\dd{\mt\jpi}generic extensions are adequately represented 
by real names. 
Then Theorem~\ref{npn} will show effects of the 
property of being a non-principal name.
 
\bpro
\lam{r2n}
Suppose that\/ $\jpi$ is a \mre\ \muf, $G\sq\mt\jpi$ is 
generic over the ground set universe\/ $\rV$,   
and\/ $x\in\rV[G]\cap\bn$.  Then 
\ben
\renu
\itlb{r2n1}
there is a true\/ \rn\jpi\ $\rc\in\rV$ 
such that\/ $x=\rc[G]\;;$ 

\vyk{
\itlb{r2n2}
if\/ $x\in\rV[G]\cap\dn$ there is a dyadic
true\/ \rn\jpi\  
$\rc=\sis{\kkc ni}{n,i<\om}\in\rV$ such that\/ $x=\rc[G]\;;$ 
%and\/ $\kkc ni=\pu$ whenever\/ $i\ge2\;;$ 
}

\itlb{r2n3}
if\/ $\mt\jpi$ is a CCC forcing in\/ $\rV$  
then there is a\/ {\ubf small} true 
\rn\jpi\ $\rd\in\rV$ with\/ $x=\rd[G]$.
\een
\epro
\bpf
\ref{r2n1}
%, \ref{r2n2} are 
is an instances of a general forcing theorem  
(see Remark~\ref{ai} on the effect of \mre ity).
To prove \ref{r2n3}, pick a \qn\  
$\rc$ by \ref{r2n1}, 
extend each set $\kkc n{}=\bigcup_i\kkc ni$ 
to an open dense set $O_n$ by closing strongwards,
choose maximal antichains $A_n\sq O_n$ in those sets
--- which have to be countable by CCC, and then let
$A_{ni}=A_n\cap \kkc ni$ and 
%$\rd=\sis{A_{ni}}{n,i<\om}$.
$\rd=\ens{\ang{\zp,n,i}}{\zp\in A_{ni}}$.
\epf

\bte
\lam{npn}
Let\/ $\jpi$ be a \mre\ \muf. 
Then
% and\/ $\rc$ be a \/ \dd{\mt\jpi}name. 
\ben
\renu
\itlb{npn1} 
if a set\/ $G\sq\mt\jpi$ is generic over the ground set 
universe\/ $\rV$, $\ang{\xi,k}\in\abc\jpi$, 
and\/ $x\in\rV[G]\cap\bn,$
then\/ $x\ne x_{\xi k}[G]$ if and only if 
there is a 
%dyadic 
true\/ \rn\jpi\ $\rc$, 
non-principal over\/ $\jpi$ at\/ $\xi,k$  
and such that\/ $x=\rc[G]$.

\itlb{npn2} 
if\/ $\rc$ is a \/ \rn\jpi, $\ang{\xi,k}\in\abc\jpi$, 
$\jqo$ is a \muf,   
$\jpi\sse\rc\xi k\jqo$, and a set\/ $G\sq\mt{\jpi\kw\jqo}$ is generic 
over\/ $\rV$ then\/ 
$\rc[G]\nin
%\bigcup_{\ang{\xi,k}\in\abc\jpi}
\bigcup_{Q\in\jqo(\xi,k)}[Q]$.
\een
\ete
\bpf
\ref{npn1}
Suppose that $x\ne x_{\xi k}[G]$.
%Let $\bo\in\dn$ be defined by $\bo(n)=0$, $\kaz n$. 
%Thus $x_{\xi k}\ne\bo$. 
By a known forcing theorem, there is a 
%dyadic
true \rn\jpi\ $\rc$ such that $x=\rc[G]$ and 
$\mt\jpi$ forces that $\rc\ne x_{\xi k}[\uG]$.
It remains to show that 
$\rc$ is a non-principal name over $\jpi$ at $\xi,k$. 
We have to prove that the set
$$
\ddi\xi k\rc\jpi
=\ens{\zp\in\mt\jpi}{\zp\,\text{ directly forces }\,
\rc\nin[\zd\zp\xi k]}\,.
$$
is open dense in\/ $\mt\jpi$. 
The openness is clear, let us prove the density. 
Consider an arbitrary $\zq\in\mt\jpi$. 
Then $\zq$ \dd{\mt\jpi}forces $\rc\ne x_{\xi k}[\uG]$ by the 
choice of $\rc$, hence we can assume that, for some $n$, 
it is \dd{\mt\jpi}forced by $\zq$ that 
$\rc(n)\ne x_{\xi k}[\uG](n)$.
Then by Lemma~\ref{avo}\ref{avo1} there is a \mut\ 
$\zp\in\mt\jpi$, $\zp\leq\zq$, and a string $s\in \om^{n+1},$ 
such that $\zp$ \dd{\mt\jpi}forces $s\sq \rc$. 
Now it suffices to show that $s\nin \zd\zp\xi k$.
Suppose otherwise:  $s\in \zd\zp\xi k$. 
Then the tree $T=\req{\zd\zp\xi k}s$ still belongs to $\mt\jpi$. 
Therefore the \mut\ $\zr$ defined by 
$\zd\zr\xi k=T$ and 
$\zd\zr{\xi'}{k'}=\zd\zp{\xi'}{k'}$ for each pair 
$\ang{\xi',k'}\ne \ang{\xi,k}$, belongs to $\mt\jpi$ and 
satisfies $\zr\leq\zp\leq\zq$.
However $\zr$ directly forces both 
$\rc(n)$ and $x_{\xi k}[\uG](n)$ 
to be equal to one and the same value $\ell=s(n)$, which 
contradicts to the choice of $n$.

To prove the converse let $\rc\in\rV$ 
be a real name non-principal over $\jpi$ at $\xi,k$, 
and $x=\rc[G]$.
Assume to the contrary that $\ang{\xi,k}\in\abc\jpi$ and 
$x=x_{\xi k}[G]$. 
There is a \mut\ $\zq\in G$ which \dd{\mt\jpi}forces 
$\rc=x_{\xi k}[\uG]$.
As $\rc$ is non-principal, 
there is a stronger \mut\ $\zp\in G\cap \ddi\xi k\rc\jpi$, 
$\zp\leq\zq$. 
Thus $\zp$ directly forces $\rc\nin[\zd\zp\xi k]$, and hence 
\dd{\mt\jpi}forces the same statement.
Yet $\zp$ \dd{\mt\jpi}forces $\rpi_{\xi k}\in[\zd\zp\xi k]$, 
of course, and this is a contradiction.\vom

\ref{npn2}
Suppose towards the contrary that 
%$\ang{\xi,k}\in\abc\jpi$, 
$Q\in\jqo(\xi,k)$ and 
$\rc[G]\in[Q]$.
By definition, $\jqo$ locks, over\/ $\jpi$, the set
$$
\dqu\rc Q\jpi
=
\ens{\zr\in\mt\jpi} 
{\zr\,\text{\rm\ directly forces }\, \rc\nin[Q]}\,. 
$$
Therefore in particular $\dqu\rc Q\jpi$ is pre-dense in 
$\mt{\jpi\kw\jqo}$ by Lemma~\ref{pqn}.  
We conclude that $G\cap \dqu\rc Q\jpi\ne \pu$. 
In other words, there is a \mut\ $\zr\in\mt\jpi$  
which directly forces $\rc\nin[Q]$. 
It easily follows that $\rc[G]\nin[Q]$, 
which is a contradiction.
\epf

\parf{Combining refinement types} 
\las{comb}  

Here we summarize the properties 
of generic refinements considered above. 
The next definition combines the 
refinement types $\ssa{D}$, $\ssb {\zD}$, $\sse\rc\xi k$.

\index{zzpissqMqo@$\jpi\ssm\cM\jqo$}%
\index{multiforcing!refinementM@refinement, $\jpi\ssm\cM\jqo$}%
\index{refinementM@refinement, $\jpi\ssm \cM\jqo$}%

\bdf
\lam{extM}
Suppose that $\jpi\bssq\jqo$ are \muf s 
and $\cM\in\hc$ is any set. 
Let $\jpi\ssm\cM\jqo$ 
\kmar{ssm mm}
mean that the four following requirements hold:
\ben
\nenu
\itlb{extM1}
if $\ang{\xi,k}\in\abc\jpi$, $D\in\cM$, $D\sq\jpi(\xi,k)$, 
$D$ is pre-dense in 
$\jpi(\xi,k)$, then 
$\jpi(\xi,k)\ssa D\jqo(\xi,k)$;

\itlb{extM2}
if ${\zD}\in\cM$, ${\zD}\sq\mt\jpi$, 
%\imar{extM2}
${\zD}$ is open dense in $\mt\jpi$,
then $\jpi\ssb {\zD}\jqo$; 

\itlb{extM2c}
if $\rc\in\cM$ is a \rn\jpi\  
%\imar{extM2c}
then $\jpi\ssb {\rc}\jqo$; 

\itlb{extM3}
if $\ang{\xi,k}\in\abc\jpi$ and $\rc\in\cM$ is a 
%\imar{extM3}
\rn\jpi, \rit{non-principal over\/ $\jpi$ at\/ $\xi,k$}, 
then $\jpi\sse\rc\xi k\jqo$.\qed 
\een
\eDf

\bcor
[of lemmas \ref{pqm}, \ref{pqn}, \ref{pqs}, \ref{pqo}]
\lam{mpq}
Let\/ $\jpi,\jqo,\jsg$ be \muf s and\/ 
$\cM$ be a countable set.
Then$:$
\ben
\renu
\itlb{mpq1}
if\/ $\jpi\ssm\cM\jqo\bssq\jsg$ then\/ $\jpi\ssm\cM\jsg\;;$ 
%\imar{pqo1}

\itlb{mpq2}
if\/ $\sis{\jpi_\al}{\al<\la}$  
%\imar{pqo2}
is a\/ \dd\bssq increasing sequence in\/ \mfp, 
$0<\mu<\la$, $\jpi=\bkw_{\al<\mu}\jpi_\al$, 
and\/ $\jpi\ssm{\cM}\jpi_\mu$, then\/
$\jpi\ssm{\cM}\jqo=\bkw_{\mu\le\al<\la}\jpi_\al$.\qed 
\een                                             
\ecor

\bcor
\lam{geec} 
If\/ $\jpi$ is a small \muf, $\cM\in\hc$, and\/ 
$\jqo$ is an\/ \dd\cM generic refinement of\/ $\jpi$ 
{\rm(exists by Proposition~\ref{gee}!)}, 
then\/ 
$\jpi\ssm\cM\jqo$.
\ecor
\bpf
We have $\jpi\ssm\cM\jqo$ by a combination of 
\ref{dj}\ref{dj3}, \ref{uu4}, \ref{uu4c}, and 
\ref{K}.
\epf

%\ble\lam{gele} ????????????\ele

%%%%%%%%%%%%%%

\gla{The forcing notion} 
\las{sek4}  

In this chapter we define the forcing notion to prove 
the main theorem. 
It will have the form $\mt\jPi$, for a certain \muf\ 
$\jPi$ with $\abc\jPi=\omi\ti\om$. 
The \muf\ $\jPi$ will be equal to the componentwise 
union of terms of a certain 
increasing sequence $\vjPi$ of small \muf s.
And quite a complicated construction of this sequence 
in $\rL$ will make use 
of some ideas related to diamond-style constructions, 
as well as to some sort of definable genericity.

\parf{Increasing sequences of small \muf s} 
\las{inc}

Recall that $\mfp$ is the set of all \muf s (Section~\ref{muft}). 
\imar{mfp}
Let $\mf\sq\mfp$ be the set of all 
\imar{mf}
\rit{small \mdi} \muf s; s accounts for both
\rit{small} and \rit{\mdi}. 
\index{multiforcing!MF@$\mf$}%
\index{zzMF@$\mf$}%
Thus a \muf\ $\jpi\in\mfp$ belongs to $\mf$ if $\abc\jpi$ 
is (at most) countable and if $\ang{\xi,k}\in \abc\jpi$ then 
$\jpi(\xi,k)$ is a small special (Definition~\ref{ptf2}) 
forcing in $\ptf$.

\bdf
\lam{sdf}
Let $\vmf$, resp., $\vmi$ be the set of all 
\index{zzMF-@$\vmf$}%
\index{zzMF-w@$\vmi$}% 
\index{length!$\len\vjpi$}% 
\index{zzlenpi@$\len\vjpi$}% 
\dd\bssq increasing sequences 
$\vjpi=\sis{\nor\jpi\al}{\al<\ka}$ of 
\imar{vmf vmi}
\muf s $\nor\jpi\al\in\mf$, 
of length $\ka=\len\vjpi<\omi$, resp., $\ka=\omi$,
which are \rit{domain-continuous}, in the sense that 
\index{domain-continuous}%
%\index{dom-continuous}
if $\la<\ka$ is a limit ordinal then 
$\abc{\nor\jpi\la}=\bigcup_{\al<\la}\abc{\nor\jpi\al}$.
%All terms of sequences in $\vmf\cup \vmi$ 
%are small special \muf s.
\index{multisequence}%
Sequences in $\vmf\cup\vmi$ are called \rit{\muq s}.
We order $\vmf\cup\vmi$ by the usual 
relations $\sq$ and $\su$ of extension of sequences. 
\bit
\item 
Thus 
$\vjpi\su\vjqo$ iff $\ka=\len\vjpi<\la=\len\vjqo$ and 
$\nor\jpi\al=\nor\jqo\al$ for all $\al<\ka$.

\item
In this case, if $\cM$ is any set, and $\nor\jqo\ka$ 
(the first term of $\vjqo$ absent in $\vjpi$) 
satisfies $\jpi\ssm\cM\nor\jqo\ka$, where 
$\jpi=\bkw_{\al<\ka}\nor\jpi\al$, 
then we write $\vjpi\su_\cM\vjqo$. 
\index{multiforcing!extension, $\vjpi\su\vjqo$}%
\index{zzpisuqo-@$\vjpi\su\vjqo$}%
\index{multiforcing!Mextension@\dd\cM extension, 
$\vjpi\su_\cM\vjqo$}%
\index{zzpisuMqo-@$\vjpi\su_\cM\vjqo$}%
\index{zzsuM@$\su_\cM$}%
\eit
If $\vjpi$ is a \muq\ in $\vmf\cup\vmi$ then 
let $\mt\vjpi= \mt\jpi$, where 
\index{multitree!$\mt\vjpi$}% 
\index{zzMTpi-@$\mt\vjpi$}% 
$\jpi=\bkw\vjpi=\bkw_{\al<\ka}\nor\jpi\al$ 
(componentwise union), and $\ka=\dom\vjpi$.
Accordingly a (true) \rit{\rn\vjpi} will mean 
\index{real name!pi-realname@\rn\vjpi}%
\index{real name!pi-realname@\rn\vjpi!true}%
a (true) \rn\jpi. 
\edf

\bcor
%[of Corollaries~\ref{geec} and \ref{mpq}]
\lam{xisc}
Suppose that\/ $\ka<\la<\omi$, $\cM$ is a countable set, 
and\/ $\vjpi=\sis{\nor\jpi\al}{\al<\ka}$ is a \muq\ 
in\/ $\vmf$. 
Then$:$
\ben
\renu
\itlb{xisc0}
the componentwise union\/ 
$\jpi=\bkw\vjpi=\bkw_{\al<\ka}\nor\jpi\al$ is 
a \mre\ \muf$;$ 

\itlb{xisc1}
there is a \muq\/ $\vjqo\in\vmf$ satisfying\/ 
$\len\vjqo=\la$ and\/ $\vjpi\su_\cM\vjqo\;;$

\itlb{xisc1+}
if moreover\/ $\sis{s_\al}{\al<\la}$ is
an\/ \dd\su increasing sequence of countable sets\/
$s_\al\sq\omi\ti\om$,
$s_\al=\abc{\nor\jpi\al}$ for all\/ $\al<\ka$,
and\/ $s_\ga=\bigcup_{\al<\ga}s_\al$ for all
limit\/ $\ga<\la$, then
there is a \muq\/ $\vjqo\in\vmf$ satisfying\/ 
$\len\vjqo=\la$,
$\abc{\nor\jqo\al}=s_\al$ for all\/ $\al<\la$,
and\/ $\vjpi\su_\cM\vjqo\;;$

\itlb{xisc1*}
if\/ $\vjpi,\vjro,\vjqo\in\vmf$
and\/ $\vjpi\su_\cM\vjro\sq\vjqo$ then\/ 
$\vjpi\su_\cM\vjqo\;;$ 

\itlb{xisc2}
if\/ $\vjqo=\sis{\nor\jqo\al}{\al<\la}\in\vmf$ and\/ 
$\vjpi\su_\cM\vjqo$ then\/ 
$\jpi=\bkw_{\al<\ka}\nor\jpi\al\ssm\cM\nor\jqo\ba$ 
whenever\/ $\la\le\ba<\mu$, and also\/ 
$\jpi\ssm\mm\jqo'=\bkw_{\la\le\ba<\mu}\nor\jqo\ba$, 
therefore
\ben
\aenri
\itlc{xisc2x}
$\mt{\jqo'}$ is open dense in\/ $\mt\vjqo\,,$

\itlc{xisc2a}
if\/ $\ang{\xi,k}\in\abc\jpi$, $D\in\mm$, $D\sq\jpi(\xi,k)$, 
$D$ is pre-dense in\/ $\jpi(\xi,k)$, then\/ 
$D$ remains pre-dense in\/ $\jpi(\xi,k)\cup\jqo(\xi,k)\,,$ 

\itlc{xisc2b}
if ${\zD}\in\mm$, ${\zD}\sq\mt\vjpi$, 
${\zD}$ is open dense in\/ $\mt\vjpi$,
then ${\zD}$ is pre-dense in $\mt{\jpi\kw\jqo'}=\mt\vjqo\,.$ 
%
%\itlc{xisc2c}
%any true\/ \rn\vjpi\  ${\rc}\in\mm$  
%is a true\/ \rn\vjqo\ as well. 
\een
\een
\ecor
\bpf
\ref{xisc0} 
Make use of Lemma~\ref{pqr}\ref{pqr3}.

\ref{xisc1} 
We define terms $\nor\jqo\al$ of the \muq\ $\jqo$ required
by induction.

Naturally put $\nor\jqo\al=\nor\jpi\al$ for each $\al<\ka$.

Now suppose that $\ka\le\ga<\la$, \muf s
$\nor\jqo\al$, $\al<\ga$, are defined, and 
$\vjro=\sis{\nor\jqo\al}{\al<\ga}$ is a \muq\ in $\vmf$.
To define $\nor\jqo\ga$, assume first that $\ga$ is limit.
Let $\jro=\bkw\vjro=\bkw_{\al<\ga}\nor\jqo\al$ 
(componentwise union). 
We can assume that $\cM$ contains $\vjro$ and satisfies 
$\ga\sq\cM$ (otherwise take a bigger set). 
By Proposition~\ref{gee}, there is an \dd\cM generic 
refinement $\jqo$ of $\jro$. 
By Theorem~\ref{dj}, $\jqo$ is small \mdi\ \muf,  
$\jro\bssq\jqo$, and $\nor\jro\al\bssq\jqo$ for all $\al<\ga$.
In addition $\jro\ssm\cM\jqo$ by Corollary~\ref{geec}.
We let $\nor\jro\ga=\jqo$. 
The extended \muq\  
$\vjro_+=\sis{\nor\jro\al}{\al<\ga+1}$ 
belongs to $\vmf$ and satisfies $\vjro\su_\cM\vjro_+$.\vom

\vyk{
Iterating this procedure with suitable sets $\cM_\ga$ 
at each step $\ka\le\ga<\la$, we extend $\vjro$ to a 
sequence $\vjqo\in\vmf$ with $\len\vjqo=\la$ and 
$\vjro\su \vjqo$.
Then by definition $\vjpi\su_\cM\vjqo$.
}

\ref{xisc1+}
The proof is similar, with the extra care of
$\abc{\nor\jqo\al}=s_\al$.\vom

To prove the main claim of \ref{xisc2} make use of
Corollary \ref{mpq}.\vom

\ref{xisc1*}
The relation $\vjpi\su_\cM\vjqo$ involves only the first
term of $\vjqo$ absent in $\vjpi$. 
\vom

To prove \ref{xisc2x} apply Corollary~\ref{pqrC}.\vom 

\ref{xisc2a}
As $\jpi\ssm\mm\jqo'$ and $D\in \mm$, 
we have $\jpi(\xi,k)\ssa D\jqo(\xi,k)$. 
Therefore $D$ is pre-dense in\/ $\jqo(\xi,k)$ by 
Lemma~\ref{pqm}\ref{pqm1}.\vom 

\ref{xisc2b}
Similarly $\jpi\ssb \zD\jqo'$, $\zD$ is pre-dense in 
$\mt\vjqo$ by 
Lemma~\ref{pqn}\ref{pqn0}.
%\vom 
%
%\ref{xisc2c}
%Similarly, make use of Lemma~\ref{pqs}\ref{pqs1}.
\epf

Our plan regarding the forcing notion for Theorem~\ref{mt} 
will be to define a certain \muq\ $\vjPi$ in $\vmi$ and 
the ensuing \muf\ $\jPi=\bkw\vjPi$ with 
remarkable properties related to definability and its own 
genericity of some sort. 
But we need first to introduce an important notion 
involved in the construction.

%\cite{ian2}

\parf{Layer restrictions of \muf s and deciding sets} 
\las{bl}

The construction of the mentioned \muf\ $\jPi$ will 
be maintained in such a way that different 
\rit{layers} $\sis{\jPi(k,\xi)}{\xi<\omi}$, $k<\om$, 
appear rather independent of each other, albeit the
principal inductive parameter will be $\xi$ rather than $k$. 
To reflect this feature, we introduce here a suitable 
notation related to layer restrictions.
%
%Recall that $\mf$ is the set of all \rit{small special} \muf s, 
%Section~\ref{dia}.
%Let $\mfp$ be the bigger set of \rit{all} \muf s. 
%\imar{mfp}%
%\index{multiforcing!MF+@$\mfp$}%
%\index{zzMF+@$\mfp$}%
%
If $m<\om$ then, using a special ``layer restriction'' 
symbol $\mes{}$ to provide a transparent 
distinction from the ordinary restriction $\res$, 
we define sets of \mut s:%
$$
\bay{lcl}
\mtl m%
\kmar{mtl m}%
&=& \text{% 
\index{layer restriction!MTI<m@$\mtl m$}%
\index{zzMTI<m@$\mtl m$}%
all \mut s $\zp\in\md$ 
such that $\abc\zp\sq\omi\ti m$, } \\[1ex]

\mtg m&=&\text{%  
\kmar{mtg m}%
\index{layer restriction!MTI>m@$\mtg m$}%
\index{zzMTI>m@$\mtg m$}%
all \mut s $\zp\in\md$ 
with $\abc\zp\sq\omi\ti {(\om\bez m)}$,} \\[1ex]

\mte m&=&\text{% 
\kmar{mte m}%
\index{layer restriction!MTIm@$\mte m$}%
\index{zzMTIm@$\mte m$}%
all \mut s $\zp\in\md$ 
such that $\abc\zp\sq\omi\ti \ans m$,}
\eay
$$
%}%
and, given a \muf\ $\jpi$, define 
\kmar{ntl-ge jpi m}%
$\ntl\jpi m$, $\ntg\jpi m$, $\nte\jpi m$ similarly.
\index{layer restriction!MTpI<m@$\ntl\jpi m$}%
\index{zzMTpI<m@$\ntl\jpi m$}%
\index{layer restriction!MTpI>m@$\ntg\jpi m$}%
\index{zzMTpI>m@$\ntg\jpi m$}%
\index{layer restriction!MTpI=m@$\nte\jpi m$}%
\index{zzMTpI=m@$\nte\jpi m$}%
Accordingly if $\zp\in \md$ then define 
the \rit{layer restriction}
$\prl\zp m\in \mtl m$ so that 
\kmar{prl zp m}%
\index{layer restriction!pI<m@$\prl\zp m$}%
\index{zzpI<m@$\prl\zp m$}%
$\abc{\prl\zp m}=\ens{\ang{\xi,k}\in\abc\zp}{k<m}$  
and $\zdp{\prl\zp m}\xi k=\dpi\xi k$ whenever 
$\ang{\xi,k}\in\abc{\prl\zp m}$. 
Define 
$\prg\zp m\in \mtg m$, 
\kmar{prg-e zp m}%
$\pre\zp m\in \mte m$  similarly.%
\index{layer restriction!pI>m@$\prg\zp m$}%
\index{zzpI>m@$\prg\zp m$}%
\index{layer restriction!pIm@$\pre\zp m$}%
\index{zzpIm@$\pre\zp m$}%

The same definitions are maintained with \muf s:
$$
\bay{lcl}
\nfl m%
\kmar{nfl m}%
&=& \text{% 
\index{layer restriction!MFI<m@$\nfl m$}%
\index{zzMFI<m@$\nfl m$}%
all \muf s $\jpi\in\mf$ 
such that $\abc\jpi\sq\omi\ti m$, } \\[1ex]

\nfg m&=&\text{%  
\kmar{nfg m}%
\index{layer restriction!MFI>m@$\nfg m$}%
\index{zzMFI>m@$\nfg m$}%
all \muf s $\jpi\in\mf$ 
with $\abc\jpi\sq\omi\ti {(\om\bez m)}$,} \\[1ex]

\nfe m&=&\text{% 
\kmar{nfe m}%
\index{layer restriction!MFIm@$\nfe m$}%
\index{zzMFIm@$\nfe m$}%
all \muf s $\jpi\in\mf$ 
such that $\abc\jpi\sq\omi\ti \ans m$,}
\eay
$$
%}%
and $\npl m$, $\npg m$, $\npe m$ are defined similarly.
\kmar{npl-g-e m}%
\index{layer restriction!MFpI<m@$\npl m$}%
\index{zzMFpI<m@$\npl m$}%
\index{layer restriction!MFpI>m@$\npg m$}%
\index{zzMFpIm@$\npg m$}%
\index{layer restriction!MFpI<m@$\npe m$}%
\index{zzMFpIm@$\npe m$}%

Accordingly if $\jpi\in \mfp$ (in particular if $\jpi\in \mf$) 
and $m<\om$ then define 
the \rit{layer restriction}
$\prl\jpi m\in \npl m$ (resp., $\in \nfl m$), so that 
\imar{prl jpi m}%
\index{layer restriction!piI<m@$\prl\jpi m$}%
\index{zzpiI<m@$\prl\jpi m$}%
$\abc{\prl\jpi m}=\ens{\ang{\xi,k}\in\abc\jpi}{k<m}$  
and $\zdp{\prl\jpi m}\xi k=\dpi\xi k$ whenever 
$\ang{\xi,k}\in\abc{\prl\jpi m}$. 
Define 
$\prg\jpi m\in \npg m$, 
\imar{prg-e jpi m}%
$\pre\jpi m\in \npe m$  similarly.%
\index{layer restriction!piI>m@$\prg\jpi m$}%
\index{zzpiI>m@$\prg\jpi m$}%
\index{layer restriction!piIm@$\pre\jpi m$}%
\index{zzpiIm@$\pre\jpi m$}%

%\lam{sdf2}
A similar notation applies to \muq s   
(Definition~\ref{sdf}). 
If $m<\om$ then we let $\nvl m$, $\nvg m$, $\nve m$ 
\imar{nvl-g-e m}%
\index{layer restriction!MF-I<m@$\nvl m$}%
\index{zzMF-I<m@$\nvl m$}%
\index{layer restriction!MF-I>m@$\nvg m$}%
\index{zzMF-I>m@$\nvg m$}%
\index{layer restriction!MF-Im@$\nve m$}%
\index{zzMF-Im@$\nve m$}%
be the set of all \muq s in $\vmf$
whose all terms belong to resp.\ 
$\nfl m$, $\nfg m$, $\nfe m$. 
Define similarly $\ivl m$, $\ivg m$, $\ive m$
(\muq s of length $\omi$).%
\imar{ivl-g-e m}%
\index{layer restriction!MFi-I<m@$\ivl m$}%
\index{zzMFi-I<m@$\ivl m$}%
\index{layer restriction!MFi-I>m@$\ivg m$}%
\index{zzMFi-I>m@$\ivg m$}%
\index{layer restriction!MFi-Im@$\ive m$}%
\index{zzMFi-Im@$\ive m$}%

And further, 
if $\vjpi=\sis{\nor\jpi\al}{\al<\ka}\in\vmf$ and $m<\om$ then 
\imar{prl vjpi m}%
define  
$\prl\vjpi m=\sis{\prl{\nor\jpi\al}m}{\al<\ka}\in \nvl m$, 
\imar{prg-e vjpi m}%
and define 
$\prg{\vjpi}m\in \nvg m$, $\pre{\vjpi}m\in \nve m$ similarly.
\index{layer restriction!pi-I<m@$\prl\vjpi m$}%
\index{zzpi-I<m@$\prl\vjpi m$}%
\index{layer restriction!pi-I>m@$\prg\vjpi m$}%
\index{zzpi-I>m@$\prg\vjpi m$}%
\index{layer restriction!pi-Im@$\pre\vjpi m$}%
\index{zzpi-Im@$\pre\vjpi m$}%
The same for $\vjpi=\sis{\nor\jpi\al}{\al<\omi}\in\vmi$
%\edf

\bdf
\lam{blo}
Assume that $m<\om$. 
A \muq\ $\vjpi%=\sis{\nor\jpi\al}{\al<\vt}
\in\vmf$ 
\dd m\rit{decides} a set $W$ if either  
\index{decision!mdecides@\dd mdecides}%
\index{decision!positive}%
\index{decision!negative}%
$\prg{\vjpi}m$ 
belongs to $W$ (\rit{positive} decision) 
or there is no \muq\  
$\vjqo\in W\cap\nvg m$ 
extending $\prg{\vjpi}m$ (\rit{negative} decision).  
\edf

\ble
\lam{bexi}
If\/ $\vjpi\in\vmf$, 
$\cM$ is countable, $W$ is any set, and\/ $m<\om,$ 
then there is a \muq\/ $\vjqo\in\vmf$ such that\/ 
$\vjpi\su_\cM\vjqo$ and $\vjqo$ \dd mdecides\/ $W.$ 
\ele
\bpf
By Corollary~\ref{xisc}, there is  a \muq\  
$\vjro\in\vmf$ such that $\vjpi\su_\cM\vjro$. 
Then either $\vjro$ outright \dd mdecides $W$ negatively, 
or there is a sequence $\vjsg\in W\cap\nvg m$ 
satisfying $\prg{\vjro}m\sq\vjsg$.

On the other hand, using Corollary~\ref{xisc}\ref{xisc1+},
we get a \muq\ $\vjsg'\in\nvl m$ of the same length as $\vjsg$, 
such that $\prl{\vjro}m\sq\vjsg'$.
Therefore there exists a \muq\ $\vjqo\in\vmf$ of that
same length, satisfying $\prg{\vjqo}m=\vjsg$ and
$\prl{\vjqo}m=\vjsg'$ --- then obviously $\vjro\sq\vjqo$
and by definition $\vjqo$ decides $W$ positively.
Finally we have $\vjpi\su_\cM\vjro\sq\vjqo$, 
and hence $\vjpi\su_\cM\vjqo$ by Corollary~\ref{xisc}\ref{xisc1*}.
\epf

%\np

\parf{Auxiliary diamond sequences} 
\las{diam}

Recall that $\hc$ is the set of all 
\rit{hereditarily countable} sets (those with finite or 
countable transitive closures).

The next theorem employs the technique of diamond sequences
in $\rL$.

\bte
[in $\rL$]
\lam{gret}
There exist\/ $\id\hc1$ sequences\/ 
\index{diamond sequences!$\nos\vjpi\mu\yd\nos D\mu\yd\nos z\mu$}%
\index{zzpimu@$\nos\vjpi\mu$}%
\index{zzDmu@$\nos D\mu$}%
\index{zzzmu@$\nos z\mu$}%
$\sis{\nos\vjpi\mu}{\mu<\omi}$,
$\sis{\nos D\mu}{\mu<\omi}$,
$\sis{\nos z\mu}{\mu<\omi}$, 
such that, for every\/ $\mu$, 
$\nos D\mu$ and\/ $\nos z\mu$ are sets in $\hc$, 
$\nos\vjpi\mu\in\vmf$, 
$\len{\nos\vjpi\mu}=\mu$, 
and in addition if\/ 
$\vjPi=\sis{\nor\jPi\nu}{\nu<\omi}\in\vmi$, 
%$\jPi=\bkw_{\nu<\omi}{\nor\jPi\nu}$,
$z\in\hc$, and $D\sq \mt\vjPi$,
then the set\/ $M$ of all ordinals\/ $\mu<\omi$ such that\/ \ 
\ben
\aenr
%\atc
\itlb{gret1}
$\nos z\mu=z\;;$

\itlb{gret2}
$\nos\vjpi\mu$ is equal to the restricted sub-\muq\/ 
$\vjPi\res\mu=\sis{\nor\jPi\nu}{\nu<\mu}\;;$

\itlb{gret3}
%$\nos D\mu=D\cap\mt{\norl\jPi\mu}$, 
$\nos D\mu=D\cap\mt{\vjPi\res\mu}\;;$ 
%where $\norl\jPi\mu=\bkw_{\nu<\mu}{\nor\jPi\nu}\;;$
\een
is stationary in\/ $\omi$.
\ete
\bpf
\rit{Arguing in $\rL$, the constructible universe}, 
we  
%\lam{gred}
let $\lel$ be the canonical wellordering of $\rL$.
\index{zz<L@$\lel$}%
\index{zzHC@$\hc$}%
It is known that $\lel$ orders $\hc$ similarly to 
$\omi$, and that $\lel$ is $\id\hc1$ and has the  
\rit{goodness} property: the set of 
all \dd\lel initial segments 
$I_x(\lel)=\ens{y}{y\lel x}$, $x\in\hc$, is still $\id\hc1$.

We begin with a  $\id\hc1$ sequence of sets\/ 
$S_\al\sq\al$, $\al<\omi$, such that 
\ben
\Aenu
\itlb{1gret}
if\/ 
$X\sq\hc$ then the set\/ $\ens{\al<\omi}{S_\al=X\cap\al}$ 
is stationary in\/ $\omi$.
\een
This is a well-known instance of the diamond 
principle true in $\rL$. 
The additional definability property can be achieved 
by taking the \dd\lel least possible $S_\al$ at each 
step $\al$.
We get the following two results as easy corollaries. 

First, let $A_\mu=\ens{c_\al}{\al\in S_\mu}$, where $c_\al$ is 
the \dd\al th element of $\hc$ in the sense of the 
ordering $\lel$. 
Then $\sis{A_\mu}{\mu<\omi}$ is still a $\id\hc1$ sequence, 
and
\ben
\Aenu
\atc
\itlb{2gret}
if $X_\al\in\hc$ for all $\al<\omi$ then 
the set\/ $\ens{\mu}{A_\mu=\ens{X_\al}{\al<\mu}}$ 
is stationary in\/ $\omi$.
\een

Second, for any $\al$, if   
$A_\al=\sis{a_\ga}{\ga<\al}$, where each $a_\ga$ itself is 
equal to an \dd\om sequence $\sis{a^n_\ga}{n<\om}$, 
then let $B^n_\al=\sis{a^n_\ga}{\ga<\al}$ for all $n$. 
Otherwise let $B^n_\al=\pu$, $\kaz n$.
Then $\sid{B_\al}{n<\om}{\al<\omi}$ is still a\/ $\id\hc1$ 
system of sets in $\hc$, such that 
\ben
\Aenu
\atc
\atc
\itlb{3gret}
if $X^n_\al\in\hc$ for all $\al<\omi$, $n<\om$, then, 
for every $\mu<\omi$, the set \break
$\ens{\mu}{\kaz n\,(B^n_\mu=\ens{X^n_\al}{\al<\mu}}$ 
is stationary in $\omi$.
\een

Now things become more routinely complex.  

Let $\mu<\omi$. 
We define $\nos z\mu=\bigcup B^0_\mu$. 
If $B^1_\mu\in\vmf$ and $\len{B^1_\mu}=\mu$ 
then let $\nos\vjpi\mu=B^1_\mu$; 
otherwise let $\nos\vjpi\mu$ be equal to the \dd\lel least 
\muq\ in $\vmf$ of length $\mu$. 
(Those exist 
by Corollary~\ref{xisc}\ref{xisc1}.) 
Finally we let $\nos D\mu=\bigcup B^2_{\mu+1}$. 
%\cap\mt{\bkw\nos\vjpi\mu}$.

Let's show that the sequences of sets $\nos\vjpi\mu$,  
$\nos D\mu$, $\nos z\mu$ prove the theorem. 
Suppose that $\vjPi=\sis{\nor\jPi\nu}{\nu<\omi}\in\vmi$, 
%$\jPi=\bkw_{\nu<\omi}{\nor\jPi\nu}$,
$z\in\hc$, and $D\sq \mt\vjPi$.
Let $X^0_\al=z$, $X^1_\al=\ang{\al,\nor\jPi\al}$, 
$X^2_\al=D\cap\mt{\vjPi\res\al}$ for all $\al$.
The set 
$$
M=\ens{\mu<\omi}{B^n_\mu=\ens{X^n_\al}{\al<\mu}\:
\text{ for }\:n=0,1,2}
$$
is stationary by \ref{3gret}.
Assume that $\mu\in M$. 
Then $B^0_\mu=\ens{X^0_\al}{\al<\mu}=\ans z$, 
therefore $\nos z\mu=z$. 
Further 
$B^1_\mu=\ens{X^1_\al}{\al<\mu}
=\ens{\ang{\al,\nor\jPi\al}}{\al<\mu}=\vjPi\res\mu\in\vmf$, 
therefore $\nos\vjpi\mu=\vjPi\res\mu$.
Finally we have 
$\nos D\mu=\bigcup B^2_{\mu+1}=\bigcup_{\al\le\mu}{X^2_\al}
=D\cap\mt{\vjPi\res\mu}$, as required.
\epf

\parf{Key sequence theorem} 
\las{keyst}

Now we prove a theorem which 
introduces the key \muq\ $\vjPi$.

\bte
[$\rV=\rL$]
\lam{kyt}
There exists a \muq\/ 
$\vjPi=\sis{\jpn\al}{\al<\omi}\in\vmi$ satisfying 
the following requirements$:$ \ 
\ben
\renu
\itlb{kyt1}
if\/ $m<\om$ then the \muq\/ 
$\pre{\vjPi}m$ belongs to the class\/ $\id\hc{m+2}\;;$ 

\itlb{kyt2}
if\/ $m'<\om$ and\/ $W\sq\vmf$ is a\/ $\fs\hc{m'+1}$ set  
then there is an ordinal\/ $\ga<\omi$ such that the \muq\/ 
$\vjPi\res\ga$ \dd{m'}decides\/ $W\;;$ 

\itlb{kyt3}
if a set\/ $D\sq \mt\vjPi$ is dense in\/ $\mt\vjPi$, 
%where\/ $\jPi=\bkw\vjPi$, 
then the set\/ $Z$ of all ordinals\/ $\ga<\omi$ such that\/ 
$\vjPi\res\ga\su_{\ans{D\cap\mt{\vjPi\res\ga}}}\vjPi$, 
%$\vjPi\res\ga\su_{\ans{D\cap\mt{\pilg \ga}}}\vjPi$, 
%where\/ $\pilg \ga=\bkw{(\vjPi\res\ga)}$, 
is stationary in\/ $\omi$. 
\een 
\ete
\bpf
If $m<\om$ then let $\ufo_m(p,x)$ be a canonical universal 
\index{universal formula, $\ufo_m(p,x)$}%
$\is{}{m+1}$ formula, so that 
the family of all $\fs\hc{m+1}$ sets $X\sq\hc$ 
(those definable in $\hc$ by $\is{}{m+1}$ formulas with  
parameters in $\hc$)
is equal to  the family of all sets 
of the form $\ups_m(p)=\ens{x\in\hc}{\hc\models \ufo_m(p,x)}$, 
\index{sets $\ups_m(p)$}%
\index{zzYpsmz@$\ups_m(p)$}%
$p\in\hc$. 

\ben
\Renu
\itlb{uxo1} 
Fix $\id\hc1$ sequences\/ 
%\imar{uxo1}
$\sis{\nos\vjpi\mu}{\mu<\omi}$,
$\sis{\nos D\mu}{\mu<\omi}$, and 
$\sis{\nos z\mu}{\mu<\omi}$ 
satisfying Theorem~\ref{gret}; the terms  
$\nos D\mu\yi\nos z\mu\yi\nos\vjpi\mu$ 
of the sequences belong to $\hc$, and 
in addition   
$\nos\vjpi\mu\in\vmf$, $\len{\nos\vjpi\mu}=\mu$.\vom

\itlb{uxo2} 
Let $\mu<\omi$. 
If $\nos z\mu$ is a pair of the form 
$\nos z\mu=\ang{m,p}$ then let $\nos m\mu=m$ and  
$\nos p\mu=p$, otherwise let 
%triple of the form 
%$\nos z\mu=\ang{m,p,\da}$ then let $\nos m\mu=m$, 
%$\nos p\mu=p$, $\nos\da\mu=\da$, otherwise let 
%
\index{diamond sequences!$\nos m\mu$, $\nos p\mu$, $\nos\da\mu$}%
\index{zzmmu@$\nos m\mu$}%
\index{zzpmu@$\nos p\mu$}%
\index{zzqmu@$\nos \da\mu$}%
$\nos m\mu=\nos p\mu=0$.
%\snos
%{The dummy component $\nos\da\mu$ is involved only to 
%make sure that any pair of active components 
%$\nos m\mu$, $\nos p\mu$ occurs uncountably many times.}
\vom

\itlb{uxo3}
If $m<\om$ then let, by Lemma~\ref{bexi}, 
$\nom\vjpi m\mu\in\vmf$ be the \dd\lel least \muq\ in 
\index{diamond sequences!$\nom\vjpi m\mu$, $\mup\mu m$}%
\index{zzpimmu@$\nom\vjpi m\mu$}%
$\vmf$ which satisfies 
$\nos\vjpi\mu\su_{\ans{\nos D\mu}}\nom\vjpi m\mu$ and 
\dd mdecides the set $\ups_{m}(\nos p\mu)$. 
Let $\mup\mu m=\len{\nom\vjpi m\mu}$; 
then $\mu<\mup\mu m<\omi$. 
\index{zzmu+m@$\mup\mu m$}%
\een

\bpro
[in $\rL$]
\lam{d11}
The sequences\/ 
$\sis{\nos m\mu}{\mu<\omi}$ and\/ 
$\sis{\nos p\mu}{\mu<\omi}$ 
% and\/ $\sis{\nos \da\mu}{\mu<\omi}$ 
belong to the definability class\/ $\id\hc1$.
If\/ $m<\om$ then the sequences\/ 
$\sis{\nom\vjpi m\mu}{\mu<\omi}$ 
and\/ $\sis{\mup\mu m}{\mu<\omi}$ 
belong to the class\/ $\id\hc{m+2}$.
\epro
\bpf
Routine. 
Note that $\nom\vjpi m\mu$ and $\mup\mu m$ depend on $m$ 
through the formulas $\ufo_m(\cdot,\cdot)$, whose complexity 
strictly increases with $m\to\iy$.
\epf

Now define a \muq\ 
$\vjPi=\sis{\jpn\al}{\al<\omi}\in\vmi$ and a 
family of strictly increasing, continuous  
maps\/ $\mu_m:\omi\to\omi$, $m<\om$,
%based on \ref{uxo1}, \ref{uxo2}, \ref{uxo3}, 
as follows$:$ \ \ 
\ben
\cenu
\itlb{kez1}
Let $\mu_m(0)=0$ and $\mu_m(\la)=\tsup_{\ga<\la}\mu_m(\ga)$ 
for all $m$ and all limit $\la<\omi$.

\itlb{kez2}
Suppose that $m<\om$, $\ga<\omi$, $\mu=\mu_m(\ga)$, and
the twofold-restricted sequence   
$\pre{(\vjPi\res\mu)}m=(\pre{\vjPi}m)\res\mu$ is already
defined.
If the following holds:
\ben
\fenv
\itlc{kez2*}
$m\ge m'=\nos m\mu$ and 
$\pre{(\vjPi\res\mu)}m$ 
coincides with\/ $\pre{\nos\vjpi\mu}m$,  
%of Definition~\ref{uxo}, 
\een
then let  
$\mu_m(\ga+1)=\mup\mu{m'}$ and\/  
$\pre{(\vjPi\res{\mup\mu{m'}})}m=\pre{\nom\vjpi{m'}\mu}m$.

\itlb{kez3}
In the assumptions of \ref{kez2}, if \ref{kez2*} fails, 
%\imar{kez3}
then let $\vjro$ is the\/ \dd\lel least \muq\  
%in\/ $\nve m$ with\/ $\pre{(\vjPi\res\mu)}m\su\vjro$\/ 
in\/ $\vmf$ with\/ $\pre{(\vjPi\res\mu)}m\su\vjro$\/ 
{\rm(we refer to Corollary~\ref{xisc})}, 
and define $\mu_m(\ga+1)=\len{\vjro}$ and 
$\pre{(\vjPi\res{\mu_m(\ga+1)})}m=\pre\vjro m\;.$ 
\een 
To conclude, given $\ga<\omi$ and $m$, if an ordinal 
$\mu=\mu_m(\ga)$, and a \muq\    
$\pre{(\vjPi\res\mu)}m=(\pre{\vjPi}m)\res\mu$
are defined, then items \ref{kez2}, \ref{kez3} define 
a bigger ordinal $\mu_m(\ga+1)>\mu=\mu_m(\ga)$ and 
a longer \muq\ $\pre{(\vjPi\res{\mu_m(\ga+1)})}m$ 
satisfying 
$\pre{(\vjPi\res\mu)}m\su \pre{(\vjPi\res{\mu_m(\ga+1)})}m$.
Thus overall items \ref{kez1}, \ref{kez2}, \ref{kez3} 
of the definition 
contain straightforward instructions as how to uniquely define  
the layers $\pre{\vjPi}m$ and maps $\mu_m$ for 
different $m<\om$, independently from each other.   

From now on, fix a \muq\ 
$\vjPi=\sis{\jpn\al}{\al<\omi}\in\vmi$ 
\imar{jpn al}%
\index{key elements!Pi-al@$\jpn\al$}%
\index{multiforcing!Pi-al@$\jpn\al$}% 
\index{zzPi-al@$\jpn\al$}% 
\index{key elements!Pi-@$\vjPi$}%
\index{multiforcing!Pi-@$\vjPi$}% 
\index{zzPi-@$\vjPi$}% 
of \muf s $\jpn\al\in\mf$ and 
increasing continuous maps $\mu_m:\omi\to\omi$ 
\index{key elements!key maps $\mu_m$}%
\index{key elements!mumga@$\mu_m(\ga)$}%
\index{zzmumga@$\mu_m(\ga)$}%
\imar{mu \ m(ga)}%
defined by \ref{kez1}, \ref{kez2}, \ref{kez3}.
As the maps $\mu_m$ are continuous, the following holds:

\bpro
[in $\rL$]
\lam{cclub}
$\dC=\ens{\ga<\omi}{\kaz m\,(\ga=\mu_m(\ga))}$ is a club
\index{key elements!C@$\dC$}%
\index{set!C@$\dC$}% 
\index{zzCd@$\dC$}% 
in\/ $\omi$.\qed
\epro

To show that $\vjPi$ proves Theorem~\ref{kyt}, 
we check items \ref{kyt1}, \ref{kyt2}, \ref{kyt3}.\vom

\ref{kyt1} 
Let $m<\om$. 
Then the \muq\/ 
$\pre{\vjPi}m$ 
and the map\/ $\mu_m$ 
belong to the class\/ $\id\hc{m+2}$ by Proposition~\ref{d11};
a routine proof is omitted.\vom

\ref{kyt2} 
Suppose that\/ $m'<\om$ and\/ $W\sq\vmf$ is a\/ $\fs\hc{m'+2}$ set. 
Pick $p\in\hc$ such that $W=\ups_{m'}(p)$.
Let $z=\ang{m',p}$. 
%Let $z=\ang{m',p,0}$. 
As $\dC$ is a club, it follows 
from the choice of terms 
$\nos\vjpi\mu$, $\nos D\mu$, and $\nos z\mu$, 
by \ref{uxo1} and Theorem~\ref{gret}, 
that there is an ordinal $\ga\in\dC$ such that 
$\nos\vjpi\ga=\vjPi\res\ga$ and 
$\nos z\ga=z$ --- hence, $\nos m\ga=m'$ and $\nos p\ga=p$.

Let $\mu=\ga$; then also $\mu=\mu_m(\ga)$, $\kaz m$ --- since 
$\ga\in\dC$, and $\vjPi\res\mu=\nos\vjpi\mu$.

Then it follows from the choice of $\vjPi$ that 
item \ref{kez2} 
of the construction applies for the ordinal $\ga$ chosen 
and all $m\ge m'$. 
It follows that the \muq\ $\vjro=\nom\vjpi{m'}\mu$ and 
the ordinal $\nu=\mu_m(\ga+1)=\mup\mu{m'}$ satisfy 
$\nu=\len{\vjro}$ and 
$\pre{(\vjPi\res{\nu})}m=\pre{\vjro}m$ for all $m\ge m'$. 
In other words, $\prg{(\vjPi\res{\nu})}{m'}=\prg{\vjro}{m'}$.

However by definition $\vjro$ \dd{m'}decides the set 
$W=\ups_{m'}(p)$, and the definition of this property 
depends only on $\prg{\vjro}{m'}$.\vom

\ref{kyt3} 
%Let $\jPi=\bkw\vjPi$. 
Assume that a set\/ $D\sq \mt\vjPi$ is dense in $\mt\vjPi$, 
and $C\sq\dC$ is a club in $\omi$. 
Following the proof of \ref{kyt2}, we find 
an ordinal $\ga\in C$ such that 
$\nos\vjpi\ga=\vjPi\res\ga$, $\nos m\ga=0$, and 
$\nos D\ga=D\cap\mt{\vjPi\res\ga}$, 
%where $\norl\jPi\ga=\bkw{(\vjPi\res\ga)}$. 
%=\bkw_{\nu<\ga}{\nor\jPi\nu}$. 
%
Note that $\ga=\mu_m(\ga)$, $\kaz m$.
We have $\nos\vjpi\ga\su_{\ans{\nos D\ga}}\nom\vjpi 0\ga$  
by \ref{uxo3} (with $\mu=\ga$), 
that is, 
$$
\nos\vjpi\ga\su_{\ans{D\cap\mt{\vjPi\res\ga}}}\nom\vjpi 0\ga.
\eqno(\dag)
$$
Yet it follows from the choice of $\ga$ that 
condition \ref{kez2*} holds (for $\mu=\ga$) for all $m\ge0$. 
Then, by definition \ref{kez2}, the ordinal 
$\mu^+=\mup\ga m$ satisfies $\mu^+=\mu_m(\ga+1)$  
and $(\vjPi\res\mu^+)\reb m=(\nom\vjpi{0}\ga)\reb m$
for all $m$, that is, just $\vjPi\res\mu^+=\nom\vjpi{0}\ga$.
We conclude that 
$\vjPi\res\ga\su_{\ans{D\cap\mt{\vjPi\res\ga}}}\vjPi\res\mu^+$ 
by $(\dag)$, therefore we have
$\vjPi\res\ga\su_{\ans{D\cap\mt{\vjPi\res\ga}}}\vjPi$, 
as required.
\epF{Theorem~\ref{kyt}}

\bdf
[in $\rL$]
\lam{keys}
From now on we fix a \muq\/ 
$\vjPi=\sis{\jpn\al}{\al<\omi}\in\vmi$ satisfying 
\index{multisequence!key \muq}%
requirements 
%\ref{kyt1},  \ref{kyt2},  \ref{kyt3}  
of Theorem~\ref{kyt}, that is, 
\ben
\renu
\itlb{keys1}
if\/ $m<\om$ then the \muq\/ 
$\pre{\vjPi}m$ belongs to the class\/ $\id\hc{m+2}\;;$ 

\itlb{keys2}
if\/ $m'<\om$ and\/ $W\sq\vmf$ is a\/ $\fs\hc{m'+1}$ set  
then there is an ordinal\/ $\ga<\omi$ such that the \muq\/ 
$\vjPi\res\ga$ \dd{m'}decides\/ $W\;;$ 

\itlb{keys3}
if a set\/ $D\sq \mt\vjPi$ is dense in\/ $\mt\vjPi$, 
%where\/ $\jPi=\bkw\vjPi$, 
then the set\/ $Z$ of all ordinals\/ $\ga<\omi$ such that\/ 
$\vjPi\res\ga\su_{\ans{D\cap\mt{\vjPi\res\ga}}}\vjPi$, 
%$\vjPi\res\ga\su_{\ans{D\cap\mt{\pilg \ga}}}\vjPi$, 
%where\/ $\pilg \ga=\bkw{(\vjPi\res\ga)}$, 
is stationary in\/ $\omi$. 
\een 
We call $\vjPi$ 
\rit{the key \muq}. 
\edf

A set $U\sq \nfg m$ is \rit{dense in\/ $\nfg m$} if 
for each $\vjpi\in\nfg m$ there is a \muq\ $\vjqo\in U$ 
satisfying $\vjpi\sq\vjqo$.

\ble
\lam{18a}
If\/ $m<\om$ and\/ $W\sq\nfg m$ is a\/ $\fs\hc{m+1}$ set  
dense in\/ $\nfg m$ then there is an ordinal\/ 
$\ga<\omi$ such that\/ $\prg{(\vjPi\res\ga)}m\in W$. 
In particular, if\/ $W\sq\vmf$ is a\/ $\fs\hc{1}$ set  
dense in\/ $\vmf$ then there is\/ 
$\ga<\omi$ such that\/ $\vjPi\res\ga\in W$. 
\ele
\bpf
Apply \ref{keys}\ref{keys2}. 
The negative decision is impossible by the density. 
\epf

\vyk{
In the next lemma, if $x\in\hc$ then let $\mm(x)$ be 
the least CTM of $\zfm$ containing $x$. 

\ble
\lam{18b}
The set\/ $M$ of all ordinals\/ $\ga<\om$ such that\/ 
$\vjPi\res\ga\su_{\mm(\vjPi\res\ga)}\vjPi$ is unbounded\/ 
in\/ $\omi$. 
\ele
\bpf
Let $\ga_0<\omi$. 
Consider the set $W$ of all \muq s $\vjpi\in\vmf$ of a 
successor length $\ga+1>\ga_0$ 
(and $\vjpi(\ga)$ as the last term), such that 
$\vpi\res\ga\su_{\mm(\vjpi\res\ga)}\vjpi$.
Simple evaluation shows that $W$ is $\fs\hc1$ 
(with $\ga_0$ as the only parameter). 
The set $W$ is dense in $\vmf$ by 
Corollary~\ref{xisc}\ref{xisc1}.
Therefore $\vjPi\res{(\ga+1)}\in W$ for some $\ga<\omi$ 
by Lemma~\ref{18a}. 
Then $\ga\in M$ and $\ga\ge\ga_0$.
\epf
}

\parf{Key product forcing} 
\las{bpf}

We continue to argue in $\rL$, and we'll  make use of 
the key \muq\ $\vjPi=\sis{\jPi_\al}{\al<\omi}$  
introduced by Definition~\ref{keys}. 
% and of the club $\dC\sq\omi$  of Proposition~\ref{cclub}. 

\bdf
[in $\rL$]
\lam{bang}
Define the \muf s 
$$
\bay{lclcllllll}
\kmar{jPi}% 
\jPi&=&\bkw\vjPi&=&\bkw_{\al<\omi}\jPi_\al
&\in&\mfp,& \\[1ex] 
\index{key elements!Pi@$\jPi$}%
\index{multiforcing!Pi@$\jPi$}% 
\index{zzPi@$\jPi$}% 
\index{key elements!Pi<ga@$\pilg\ga$}%
\index{multiforcing!Pi<ga@$\pilg\ga$}% 
\index{zzPi<ga@$\pilg\ga$}% 
\kmar{pilg ga}% 
\pilg\ga&=&\bkw(\vjPi\res\ga)
&=&\bkw_{\al<\ga}\jPi_\al&\in& \mf, & 
\text{for each $\ga<\omi$}\\[1ex]
\index{key elements!Pi>ga@$\pigg\ga$}%
\index{multiforcing!Pi>ga@$\pigg\ga$}% 
\index{zzPi>ga@$\pigg\ga$}% 
\kmar{pigg ga}% 
\pigg\ga&=&\bkw(\vjPi\res(\omi\bez\ga))
&=&\bkw_{\ga\le\al<\omi}\jPi_\al&\in& \mfp, & 
\text{for each $\ga<\omi$}.
\eay%
$$%
We further define $\fP=\mt{\jPi}=\mt{\vjPi}$, 
\kmar{fP}% 
\index{key elements!Pd@$\fP$}%
\index{multiforcing!Pd@$\fP$}% 
\index{zzPd@$\fP$}% 
and, for all $\ga<\omi$,
$$
\fpl\ga=\mt{\pilg\ga}=\mt{\vjPi\res\ga}\,,\quad 
\kmar{fpl ga fpg ga}% 
\index{key elements!Pd<ga@$\fpl\ga$}%
\index{multiforcing!Pd<ga@$\fpl\ga$}% 
\index{zzPd<ga@$\fpl\ga$}% 
\fpg\ga=\mt{\pigg\ga}=\mt{\vjPi\res(\omi\bez\ga)}\,.
\eqno\qed
$$
%
%If $\xi<\omi$ and $k<\om$ then let $\fP_{\xi k}=\jPi(\xi,k)$; 
%\index{key elements!Pdxik@$\fP_{\xi k}$}%
%\index{multiforcing!Pdxik@$\fP_{\xi k}$}% 
%\index{zzPdxik@$\fP_{\xi k}$}% 
%this is a perfect tree forcing in $\ptf$ and obviously 
\eDf
\index{key elements!Pd>ga@$\fpg\ga$}%
\index{multiforcing!Pd>ga@$\fpg\ga$}% 
\index{zzPd>ga@$\fpg\ga$}% 

The \muf\ $\fP$ will be our principal forcing notion, 
\rit{the key forcing}. 
%It follows from the next lemma that
%$\fP=\prod_{\xi<\omi,\,k<\om}\jPi(\xi,k)$ (with finite support). 

\ble
[in $\rL$]
\lam{doml}
$\jPi$ is a \mre\ \muf. 
In addition,
$\abc\jPi=\omi\ti\om$, thus if\/ $\xi<\omi$ 
and $k<\om$ then there is an ordinal\/ 
$\al<\omi$ such that\/ 
$\ang{\xi,k}\in\abc{\jpn\al}$.
Therefore\/ 
$\fP=\prod_{\xi<\omi,\,k<\om}\jPi(\xi,k)$ 
(with finite support).  
\ele
\bpf
To prove the additional claim, note that
the set $W$ of all \muq s $\vjpi\in\vmf$ satisfying 
$\ang{\xi,k}\in\abc{\bkw\vjpi}$ is $\fs\hc1$ 
(with $\xi$ as a parameter of definition). 
In addition $W$ is dense in $\vmf$. 
(First extend $\vjpi$ by Corollary~\ref{xisc} 
so that is has a non-limit length and the last term, 
then make use of Corollary~\ref{xistt}.) 
Therefore by Lemma~\ref{18a}  
there is an ordinal $\ga<\omi$ such that 
$\vjPi\res\ga\in W$, as required.
\epf

If $\xi<\omi$ and $k<\om$ then, following the lemma,
let $\al(\xi,k)<\omi$
be the least ordinal $\al$ satisfying  
$\ang{\xi,k}\in\abc{\jpn{\al}}$. 
Thus a forcing $\jpn{\al}(\xi,k)\in\ptf$
is defined whenever $\al$ satisfies 
$\al(\xi,k)\le\al<\omi$, and    
$\sis{\jpn{\al}(\xi,k)}{\al(\xi,k)\le\al<\omi}$ is a 
\dd\bssq increasing sequence of countable special forcings 
in\/ $\ptf$.

Note that
$\jPi(\xi,k)= \bigcup_{\al(\xi,k)\le\al<\omi} \jpn{\al}(\xi,k)$
by construction.

\bcor
[in $\rL$]
\lam{doml2}
If\/ $k<\om$ then the sequence of ordinals\/ 
$\sis{\al(\xi,k)}{\xi<\omi}$ and the sequence of
\muf s\/ 
$\sis{\jpn{\al}(\xi,k)}{\xi<\omi,\,\al(\xi,k)\le\al<\omi}$
are\/ $\id\hc{k+2}$. 
\ecor
\bpf
By construction the following double equivalence holds:
$$
\bay{cclc}
\al<\al(\xi,k)
&\eqv&
\sus\jpi(\jpi=\jpn\al\reb k\land \ang{\xi,k}\in\dom\jpi)
&\eqv 
\\[1ex]
&\eqv&
\kaz\jpi(\jpi=\jpn\al\reb k\imp \ang{\xi,k}\in\dom\jpi)
&. \eay
$$
However $\jpi=\jpn\al\reb k$ is a $\id\hc{k+2}$ relation by 
Theorem~\ref{kyt}\ref{kyt1}.
It follows that so is the sequence $\sis{\al(\xi,k)}{\xi<\omi}$.
The second claim easily follows by the same
Definition~\ref{keys}\ref{keys1}.
\epf

\bcor
[in $\rL$, of Lemma~\ref{pqr}\ref{pqr4}]
\lam{xiden}
If\/ $\xi<\omi$, $k<\om$, and\/ $\al(\xi,k)\le\al<\omi$ 
then the set\/ $\jpn\al(\xi,k)$ 
is pre-dense in\/ $\jPi(\xi,k)$ and in\/ $\jPi$.\qed  
\ecor

In spite of Lemma~\ref{doml}, the sets $\abc{\pilg\ga}$ can be 
quite arbitrary (countable) subsets of $\omi\ti\om$. 
However we easily get the next corollary:

\bcor
[in $\rL$, of Lemma~\ref{doml}]
\lam{domc}
The set\/ $\dcp=\ens{\ga<\omi}{\abc{\pilg\ga}=\ga\ti\om}$ 
\index{key elements!C'@$\dcp$}%
\index{set!C'@$\dcp$}% 
\index{zzCd'@$\dcp$}% 
\imar{dcp}% 
is a club in\/ $\omi$.\qed
\ecor

\ble
[in $\rL$]
\lam{ccc}
$\fP$ is CCC.
\ele
\bpf
Let $A\sq\fP$ be a maximal antichain in $\fP$. 
The set 
$$
C=\ens{\ga<\omi} 
{A\cap\fpl\ga\text{ is a maximal antichain in }\fpl\ga}
$$ 
is a club in $\omi$. 
Let $D=\ens{\zp\in\fP}{\sus\zq\in A\,(p\leq q)}$; 
this is an open dense set. 
By Definition~\ref{keys}\ref{keys3}, there is an ordinal  
$\ga\in C$ such that 
$\vjPi\res\ga\su_{\ans{D\cap\fpl\ga}}\vjPi$. 
Recall that $\ga\in C$, hence $A\cap\fpl\ga$ is 
a maximal antichain in $\fpl\ga$, thus 
$D\cap\fpl\ga$ is open dense in $\fpl\ga$. 
Therefore the set $D\cap\fpl\ga$  
is pre-dense in the forcing $\mt{\vjPi}=\fP$ by 
Corollary~\ref{xisc}\ref{xisc2b}.
%%%%%
\vyk{
It follows that 
$\pilg\ga\ssm{\ans{D\cap\fpl\ga}}\jPi(\ga)$, 
by Definition \ref{sdf}, and finally,  
by Corollary~\ref{xisc}\ref{xisc2},
$$
%\label{ccc2}
\pilg\ga\ssm{\ans{D\cap\fpl\ga}}\pigg\ga \,.
\eqno(\ddag)
$$
Recall that $\ga\in C$, hence $A\cap\fpl\ga$ is 
a maximal antichain in $\fpl\ga$, thus 
$D\cap\fpl\ga$ is open dense in $\fpl\ga$. 
Therefore we have $\pilg\ga\ssb{D\cap\fpl\ga}\jPi(\ga)$
(Definition~\ref{extM}) by $(\ddag)$. 
It follows by Lemma~\ref{pqn} that the set $D\cap\fpl\ga$  
is pre-dense in the forcing $\mt{\pigg\ga\lor\pilg\ga}=\fP$.
}%
%%%%%%
We claim that $A=A\cap\fpl\ga$, so $A$ is countable.

Indeed suppose that $\zr\in A\bez\fpl\ga$.
Then $\zr$ is compatible with some $\zq\in D\cap\fpl\ga$; 
let $\zp\in D\cap\fpl\ga$, $\zp\leq\zq$, $\zp\leq\zr$.  
As $\zq\in D$, there is some $\zr'\in A$ with $\zq\leq \zr'$. 
Then $\zr=\zr'$ as $A$ is an antichain; 
thus $\zq\leq \zr\in A\bez\fpl\ga$.
However $\zq\in\fpl\ga$ and $A\cap\fpl\ga$ is 
a maximal antichain in $\fpl\ga$, thus $\zq$, and hence $\zr$, 
is compatible with some $\zr''\in A\cap\fpl\ga$.
Which is a contradiction.
\epf

\bcor
[in $\rL$]
\lam{Cjden}
If a set\/ $D\sq\fP$ is pre-dense in\/ $\fP$ 
then there is an ordinal\/ 
$\ga<\omi$ such that\/ $D\cap\fpl\ga$ is already pre-dense 
in\/ $\fP$.
\ecor
\bpf
We can assume that in fact $D$ is dense.
Let $A\sq D$ be a maximal antichain in $D$; then $A$ is 
a maximal antichain in $\fP$ because of the density of $D$.
Then $A\sq \fpl\ga$ for some $\ga<\omi$ by Lemma~\ref{ccc}.
But $A$ is pre-dense in\/ $\fP$.
\epf

\gla{Auxiliary forcing relation}
\las{afn}

Recall that $\fP=\mt\jPi$, the key forcing, 
is a product forcing notion defined (in $\rL$) in 
Section~\ref{bpf}.
Its components $\jPi(\xi,k)$ have different complexity
in $\hc$, depending on $k$ by Corollary~\ref{doml2}, 
hence there
is no way the forcing notion $\fP$ (or $\jPi$)
as a whole is definable in $\hc$. 
Somewhat surprisingly, the \dd\fP forcing relation turns 
out to be definable in $\hc$ when restricted to 
analytic formulas 
of a certain level of complexity within the usual hierarchy.
This will be established on the base of an auxiliary
forcing relation.

\parf{Auxiliary forcing: preliminaries}
\las{auxA}

{\ubf We argue in $\rL$.}
Consider the 2nd order arithmetic language, with
variables $k,l,m,n,\dots$ of type $0$ over $\omega$ and 
variables $a,b,x,y,\dots$ of type $1$ over $\bn$,  
whose atomic formulas are those of the form $x(k)=n$.
Let $\cL$ be the extension of this language,
which allows to substitute
free variables of type $0$ with natural numbers
(as usual) and free variables of type $1$ with 
%dyadic
{\ubf small real names} $\rc\in\rL$.
By \rit{\dd\cL formulas} we understand formulas of this
\index{formula!Lformula@\dd\cL formula}%
\index{zzLformula@\dd\cL formula}%
extended language.

We define natural classes $\ls1n$, $\lp1n$ ($n\ge1$)
\index{formula!LSP@$\ls1n$, $\lp1n$, $\lsp1n$}%
\index{zzLSP@$\ls1n$, $\lp1n$, $\lsp1n$}%
of \dd\cL formulas.
Let $\lsp11$ be the closure of $\ls11\cup\lp11$ under 
${\neg}\yi{\land}\yi{\lor}$ and quantifiers over $\om$.
If $\vpi$ is a formula in $\ls1n$ (resp., $\lp1n$), then 
let $\vpi^-$ be  
\index{formula!phi-@$\vpi^-$}%
\index{zzphi-@$\vpi^-$}%
the result of canonical transformation of $\neg\:\vpi$ 
to the $\lp1n$ (resp., $\ls1n$) form.

If $\vpi$ is a \dd\cL formula and $G\sq\md$ is a
pairwise compatible set of \mut s then let $\vpi[G]$
be the result of substitution of $\rc[G]$ for any
name $\rc$ in $\vpi$.
(Recall Definition~\ref{preta}.)
Thus $\vpi[G]$ is an ordinary 2nd order arithmetic
formula, which may include natural numbers and
elements of $\bn$ as parameters.

We are going to define 
a relation $\zp\foe\vjpi\vpi$ between \mut s $\zp$, 
\muq s $\vjpi$, and \dd\cL formulas $\vpi$, 
which suitably approximates the true \dd\fP forcing 
relation. 
But it depends on a two more definitions. 
%The first of them introduces classes of \muq s which 
%agree with the key \muq\ $\vjPi$ up to a certain 
%layer $m<\om$.

\bdf
\lam{wmf}
If $m<\om$ then 
$\wmf m$ consists of all \muq s $\vjpi\in\vmf$ such that
$\prl\vjpi m\su \prl\vjPi m$, that is, 
$\prl\vjpi m=(\prl\vjPi m)\res\da$, where $\da=\len\vjpi$
--- \muq s which agree with the key \muq\ $\vjPi$ 
on layers below $m$.
Obviously $\wmf{m+1}\sq\wmf m\sq\wmf 0=\vmf$.
\edf

If $\ga<\omi$ then the subsequence $\vjPi\res\ga$ of the 
key \muq\ $\vjPi$ belongs to $\bigcap_m\wmf m$, of course.
To prove the next lemma use \ref{keys}\ref{keys1}.

\ble
\lam{wmfd}
$\wmf m$ is a subset of\/ $\hc$ of definability class\/ 
$\id\hc{m+1}$.\qed
\ele

The other definition deals with models of a 
subtheory of $\zfc$.

\bdf
\lam{zfm} 
Let $\zfm$ be the theory containing all axioms of 
$\zfc^-$ (minus the Power Set axiom) plus the axiom 
of constructibility $\rV=\rL$. 
% saying that all sets are constructible. 
Any transitive model (TM) of $\zfm$ has the form 
$\rL_\al$, where $\al\in\Ord$. 
Therefore it is true in $\rL$ that for any set $x$ 
there is a least TM $\mm=\mm(x)$ of $\zfm$ 
\index{model!$\mm(x)$}%
\index{zzMx@$\mm(x)$}%
\index{model!CTM, countable transitive model}%
\index{CTM, countable transitive model}%
containing $x$.
If $x\in\hc$ 
($\hc$= all hereditarily countable sets) then 
$\mm(x)$ is a \rit{countable} transitive model (CTM), 
equal to the least CTM of $\zfm$ containing $x$. 
\edf

\parf{Auxiliary forcing}
\las{auxB}

The definition of $\zp\foe\vjpi\vpi$ goes on by 
induction on the complexity of $\vpi$.   

\ben
\cenu
\itlb{fo2}
Let $\vpi$ is a closed $\lsp11$ formula, $\vjpi\in\vmf$, 
$\zp\in\md$, but $\zp\in\mt\vjpi$ is not necessarily assumed.
We define:
\ben
\itlb{fo21}
$\zp\foe\vjpi\vpi$ iff   
\index{forcing!forc@$\fof$}%
\index{zzforc@$\fof$}%
%$\zp\in\mt{\vjpi}$, and 
\vyk{
there is an ordinal $\vt<\dom\vjpi$ and 
a \mut\  
$\zpo\in\mt{\vjpi\res\vt}$
such that $\zp\leq \zpo$ (meaning: $\zp$ is stronger),
the model $\mm=\mm(\vjpi\res\vt)$ contains $\vpi$ 
(that is, all \qn s in $\vpi$), 
$\vjpi\res\vt\su_\mm\vjpi$, and 
$\zpo$ \dd{\mt{\vjpi\res\vt}}forces $\vpi[\uG]$ 
over $\mm$ (in the usual sense)
}%
there is a CTM $\mm\mo\zfm$, 
an ordinal $\vt<\dom\vjpi$, 
and a \mut\ $\zpo\in\mt{\vjpi\res\vt}$,   
such that $\zp\leq \zpo$ (meaning: $\zp$ is stronger),
the model $\mm$ contains $\vjpi\res\vt$ 
(then contains $\mt{\vjpi\res\vt}$ as well) 
and contains $\vpi$ 
(that is, all names in $\vpi$), 
$\vjpi\res\vt\su_\mm\vjpi$, and 
$\zpo$ \dd{\mt{\vjpi\res\vt}}forces $\vpi[\uG]$ 
over $\mm$ (in the usual sense)
;

\itlb{fo22}
$\zp\wfoe\vjpi\vpi$ (weak forcing) iff   
\index{forcing!wforc@$\wfof$}%
\index{zzwforc@$\wfof$}%
there is no \muq\
$\vjpi'\in\vmf$ and $\zp'\in\mt{\vjpi'}$ 
such that $\vjpi\sq\vjpi'$, $\zp'\leq\zp$, 
and $\zp'\foe{\vjpi'}\neg\:\vpi$.
%(Compare to \ref{fo4}.) 
\een

\itlb{fo3}
If $\vpi(x)$ is a $\lp1n$ formula, $n\ge1$, then 
we define $\zp\foe\vjpi\sus x\,\vpi(x)$ iff there is a 
%stable true 
small \qn\ $\rc$ 
%\in\mm(\vjpi)$ 
such that $\zp\foe\vjpi\vpi(\rc)$.
%$\zp\in\mt\vjpi$ directly forces $\rc(n)=k$ 

\itlb{fo4}
If $\vpi$ is a closed $\lp1n$ formula, $n\ge2$, 
then we define $\zp\foe\vjpi\vpi$ iff  
%$\vpi\in\mm(\vjpi)$, 
%$\zp\in\mt\vjpi$, 
$\vjpi\in\wmf{n-2}$, and there is no \muq\   
%$\vjqo\in\vmf$  
$\vjpi'\in\wmf{n-2}$ and \mut\ $\zp'\in\mt{\vjpi'}$ 
such that $\vjpi\sq\vjpi'$,
$\zp'\leq\zp$, 
%($\zp'$ is stronger),
and $\zp'\foe{\vjpi'}\vpi^-$.\snos
{If $\vjpi$ does not belong to $\wmf{n-2}$
in \ref{fo4}, then $\zp\foe\vjpi\vpi$ holds
for any $\lp1n$ formula $\vpi$ 
by default as a \muq\ not in $\wmf{n-2}$ is 
definitely not 
extendable to a \muq\ in $\wmf{n-2}$.
This motivates the condition $\vjpi\in\wmf{n-2}$ 
in \ref{fo4}.
} 
%(but not necessarily $\zp\in\mt\vjpi$).
\een 

\bre
\lam{22c}
With $\zpo$ and $\vt$ given, the premise 
``$\zpo$ \dd{\mt{\vjpi\res\vt}}forces $\vpi[\uG]$ 
over $\mm$'' 
of \ref{fo21} does not depend on the 
choice of a CTM $\mm$ containing $\vjpi\res\vt$ 
and $\vpi$, 
% and satisfying $\vjpi\res\vt\su_\mm\vjpi$.
since if $\vpi$ is $\lsp11$ then the formula $\vpi[G]$ 
(in which all names are evaluated by some 
\dd{\mt{\vjpi\res\vt}}generic set $G$ as in \ref{preta}) 
in simultaneously true or false in all transitive 
models by the Shoenfield absoluteness theorem.
\ere

\bre
\lam{22b}
It easily holds by induction that if $\zp\foe\vjpi\vpi$ 
then $\vjpi\in\vmf$, 
%$\zp\in\mt\vjpi$, 
$\vpi$ is a closed formula in one of 
the classes $\lsp11,\,\ls1n,\,\lp1n$, $n\ge2$, 
%$\vpi$ (that is, each \qn\ in $\vpi$) belongs to 
%$\mm(\vjpi)$, 
and if $n\ge2$ and $\vpi\in\lp1n\cup\ls1{n+1}$ 
then $\vjpi\in\wmf{n-2}$.
\ere

\vyk{\bre
\lam{22a}
If $\vjpi$ does not belong to $\wmf{n-2}$
in \ref{fo4}, then $\zp\foe\vjpi\vpi$ holds
for any $\lp1n$ formula $\vpi$ 
by default as a \muq\ not in $\wmf{n-2}$ is 
definitely not 
extendable to a \muq\ in $\wmf{n-2}$.
This motivates the condition $\vjpi\in\wmf{n-2}$ 
in \ref{fo4}.
\ere
}

\ble
%[monotonicity, in $\rL$]
\lam{mont}
Assume that \muq s\/ $\vjpi\sq\vjqo$ belong to\/ $\vmf$, 
%$\zp\in\mt\vjpi$, $\zq\in\mt\vjqo$, 
$\zq,\zp\in\md$, 
$\zq\leq\zp$, 
$\vpi$ is an\/ \dd\cL formula as in\/ \ref{22b}, 
and if\/ $n\ge2$ and\/ $\vpi\in\lp1n\cup\ls1{n+1}$ 
then $\vjpi,\vjqo\in\wmf{n-2}$.
%over\/ $\vjpi$, 
Then\/ $\zp\foe\vjpi\vpi$ implies\/ $\zq\foe\vjqo\vpi$, 
and if\/ $\vpi$ belongs to\/ $\lsp11$ then\/ 
$\zp\wfoe\vjpi\vpi$ implies\/ $\zq\wfoe\vjqo\vpi$ 
as well.
\ele
\bpf
If $\vpi$ is a $\lsp11$ formula,  
$\zp\foe\vjpi\vpi$, and this is witnessed 
by $\mm$, $\vt$, $\zpo$ as in \ref{fo21},
then the exactly same $\mm$, $\vt$, $\zpo$ 
witness $\zq\foe\vjqo\vpi$. 

The induction step $\sus$, as in \ref{fo3}, is elementary.
%since any \dd{\mt\vjpi}real name  
%$\rc\in\mm$ is a \dd{\mt\vjqo}real name, see above.

Now the induction step $\kaz$, as in \ref{fo4}.
Let $\vpi$ be a closed \dd{\lp1n}formula, $n\ge2$, and
%$\zp\in\mt{\vjpi}$, 
$\zp\foe\vjpi\vpi$. 
Assume to the contrary that $\zq\foe\vjqo\vpi$
fails.
Then by \ref{fo4} there exist:
a \muq\ $\vjqo'\in\wmf{n-2}$ and \mut\ $\zq'\in\mt{\vjqo'}$ 
such that $\vjqo\sq\vjqo'$, $\zq'\leq\zq$,
and $\zq'\foe{\vjqo'}\vpi^-$.
But then $\vjpi\sq\vjqo'$ and $\zq'\leq\zp$,
hence $\zp\foe\vjpi\vpi$ fails by \ref{fo4}, which is
a contradiction.

The additional result for $\wfof$ and $\lsp11$ formulas 
is entirely similar to the induction step $\kaz$ as 
just above.
\epf

If $K$ is one of the classes
$\lsp11$, $\ls1n$, $\lp1n$ ($n\ge2$), then let $\fos K$ consist
of all triples $\ang{\vjpi,\zp,\vpi}$ such that
$\vjpi\in\vmf$, $\zp\in\mt\vjpi$, $\vpi$ 
is a closed \dd\cL formula of class $K$, 
and if $n\ge2$ and $\vpi\in\ls1n\cup\lp1n$ 
then $\vjpi\in\wmf{n-2}$,  
and finally $\zp\foe\vjpi\vpi$.
Then $\fos K$ is a subset of $\hc$.

\ble
[definability, in $\rL$]
\lam{deff}
$\fos{\lsp11}$ 
%and\/ $\fos{\ls12}$ 
belongs to\/ $\id\hc1$.
If\/ $n\ge2$ then\/ 
$\fos{\ls1n}$ belongs to\/ $\is\hc{n-1}$ 
and\/ $\fos{\lp1{n}}$ belongs to\/ $\ip\hc{n-1}$. 
\ele
\bpf
Relations like ``being an MSP'',  
``being a formula in $\lsp11$, $\ls1n$, $\lp1n$'', 
$\zp\in \mt{\vjro}$, forcing over a CTM, etc.\  
are definable in $\hc$ by bounded formulas, hence 
$\id\hc1$. 
On the top of this, the model $\mm$ can be 
tied by both $\sus$ and $\kaz$ in \ref{fo21}, 
see Remark~\ref{22c}.
%`being $\mm(x)$' is a $\id\hc1$ relation. 
This wraps up the $\id\hc1$ estimation for $\lsp11$. 

%Now proceed by induction. 
The inductive step by \ref{fo3} is quite simple.  
% since the bounded quantifier $\sus\rc$ is 
%essentially 
%relativized to $\mm(\vjpi)$, 
%(Remark~\ref{22b}), 
%so it keeps the definability class over $\hc$. 

Now the step by \ref{fo4}.
Assume that 
%$\vpi$ is a closed \dd{\lp1n}formula, 
$n\ge2$, and it is already established that 
$\fos{\ls1n}\in\is\hc{n-1}$. 
Then $\ang{\vjpi,\zp,\vpi}\in\fos{\lp1n}$
iff $\vjpi\in\wmf{n-2}$, $\zp\in\md$, $\vpi$ 
is a closed $\lp1n$ formula in $\mm(\vjpi)$, 
% over $\msp\mm\vjpi$, 
and, by \ref{fo4}, there exist no triple 
$\ang{\vjpi',\zp',\psi}\in\fos{\ls1n}$
such that $\vjpi'\in\wmf{n-2}$, 
$\vjpi\sq\vjpi'$, $\zp'\in\mt{\vjpi'}$, $\zp'\leq\zp$, 
and $\psi$ is $\vpi^-$.
Evaluating the key conjunct $\vjpi'\in\wmf{n-2}$ by 
Lemma~\ref{wmfd} as $\id\hc{n-1}$, we get a required 
estimation $\ip\hc{n-1}$ of $\fos{\lp1n}$. 
\epf

\parf{Forcing simple formulas}
\las{mainf}

{\ubf We still argue in $\rL$.} 
The following results are mainly related to 
the relation $\fof$  with respect to 
formulas in the class $\lsp11$.

\ble
[in $\rL$]
\lam{lsp}
Assume that\/ $\vjpi\in\vmf$, $\vjqo\in\vmf\cup\vmi$, 
$\vjpi\sq\vjqo$, $\zp\in\mt\vjpi$, 
$\vpi$ is a formula in\/ $\lsp11$, 
$\nn\mo\zfm$ is a TM containing\/ $\vjqo$ and\/ $\vpi$, 
and\/ $\zp\foe\vjpi\vpi$. 
%by\/ \ref{fo21} via some\/ $\vt<\dom\vjpi$.
% and\/ $\zpo\in\mt{\vjpi\res\vt}$. 
Then\/ $\zp$ \dd{\mt\vjqo}forces\/ $\vpi[\uG]$ 
over\/ $\nn$.
%
%In particular if\/ $\nn$ contains\/ $\vjpi,\vpi$ then\/    
%$\zp$ \dd{\mt\vjpi}forces\/ $\vpi[\uG]$ over\/ $\nn$.
\ele
\bpf
By definition there is an ordinal $\vt<\dom\vjpi$, 
a \mut\ $\zpo\in\mt{\vjpi\res\vt}$, 
and a CTM $\mm\mo\zfm$ containing $\vpi$ and 
$\vjpi\res\vt$, such that  
$\vjpi\res\vt\su_\mm\vjpi$, $\zp\leq\zpo$, 
%the model $\mm=\mm(\vjpi\res\vt)$ contains $\vpi$, 
%$\vjpi\res\vt\su_\mm\vjpi$, 
and 
$\zpo$ \dd{\mt{\vjpi\res\vt}}forces $\vpi[\uG]$ 
over $\mm$. 
We can \noo\ assume that $\mm\sq\nn$ 
(by the same reference to Shoenfield as in Remark~\ref{22c}).

Now suppose that $G\sq\mt\vjqo$ is a set 
\dd{\mt\vjqo}generic 
over $\nn$ and $\zp\in G$ --- then $\zpo\in G$, too. 
We have to prove that $\vpi[G]$ is true in $\nn[G]$. 

We claim that the set $G'=G\cap\mt{\vjpi\res\vt}$ is 
\dd{\mt{\vjpi\res\vt}}generic over $\mm$.
Indeed, let a set $\zD\in\mm$, $\zD\sq\mt{\vjpi\res\vt}$, 
be open dense in\/ $\mt{\vjpi\res\vt}$.
Then, as ${\vjpi\res\vt}\su_\mm\vjqo$ by 
Corollary~\ref{xisc}\ref{xisc1*}, 
$\zD$ is pre-dense in $\mt\vjqo$ 
by \ref{xisc}\ref{xisc2b}, 
and hence $G\cap\zD\ne\pu$ by the choice of $G$. 
It follows that  $G'\cap\zD\ne\pu$.

We claim that $\rc[G]=\rc[G']$ for any name 
$\rc\in\mm$, in particular, for any name in $\vpi$. 
Indeed, as $G'\sq G$, the otherwise occurs  
by Definition~\ref{preta} only if for some $n,i$ 
and $\zq'\in\kkc ni$ there is $\zq\in G$ satisfying 
$\zq\leq \zq'$, but there is no such $\zq$ in $G'$. 
Let $\zD$ consist of 
all \mut s $\zr\in\mt{\vjpi\res\vt}$ either satisfying 
$\zr\leq \zq'$ or somewhere \ad\ with $\zq'$. 
Then $\zD\in\mm$ and $\zD$ is open dense in 
$\mt{\vjpi\res\vt}$. 
Therefore $\zD\cap G'\ne\pu$ by the above, so let 
$\zr\in \zD\cap G'$. 
If $\zr\leq \zq'$ then we get a contradiction with the 
choice of $\zq'$. 
If $\zr$ is somewhere \ad\  with $\zq'$ then 
we get a contradiction with the choice of $\zq$ as 
both $\zq,\zr$ belong to the generic filter $G$.

It follows that $\vpi[G]$ coincides with $\vpi[G']$.

Note also that $\zpo\in G'$.
We conclude that $\vpi[G']$ is true in $\mm[G']$ as 
$\zpo$ forces $\vpi[\uG]$ over $\mm$. 
The same formula $\vpi[G]$ is true in 
$\nn[G]$ by Shoenfield.
\epf  

\ble
\lam{wf1}
Let\/ $\vjpi\in\vmf$, $\zp\in\mt\vjpi$, 
$\vpi$ be a formula in\/ $\lsp11$. 
Then
%and\/ 
\ben
\renu
\itlb{wf11}\msur
$\zp\foe\vjpi\vpi$ and\/ $\zp\foe\vjpi\neg\:\vpi$ 
cannot hold together$;$

\itlb{wf12}
if\/ $\zp\foe\vjpi\vpi$ then\/ $\zp\wfoe\vjpi\vpi\;;$ 

\itlb{wf13}
if\/ $\zp\wfoe\vjpi\vpi$ then there is a \muq\ 
$\vjqo\in\vmf$ such that\/ $\vjpi\su_{\mm(\vjpi)}\vjqo$ 
and\/ $\zp\foe\vjqo\vpi$;

\itlb{wf14}\msur
$\zp\wfoe\vjpi\vpi$ and\/ $\zp\wfoe\vjpi\neg\:\vpi$ 
cannot hold together.
\een
\ele
\bpf 
\ref{wf11}
Otherwise $\zp$ \dd{\mt\vjpi}forces both $\vpi[\uG]$ 
and $\neg\:\vpi[\uG]$ over a large enough CTM $\mm$, by 
Lemma~\ref{lsp}, which cannot happen.

\ref{wf12}
Assume that $\zp\wfoe\vjpi\vpi$ fails, hence there is 
a \muq\ $\vjqo\in\vmf$ and a \mut\ $\zq\in\mt\vjqo$ 
such that $\zq\leq\zp$ and $\zq\foe\vjpi\neg\:\vpi$.
But Lemma~\ref{mont} implies $\zq\foe\vjpi\vpi$, 
which contradicts to \ref{wf11}. 

\ref{wf13}
Let $\mm\mo\zfm$ be a CTM containing $\vpi$ and $\vjpi$.
By Corollary~\ref{xisc}\ref{xisc1},   
there is a \muq\ \muq\ 
$\vjqo\in\vmf$ with\/ $\vjpi\su_{\mm}\vjqo$. 
%Let\/ $\vt=\dom\vjpi$. 
We claim that $\zp\,$ \dd{\mt\vjpi}forces 
$\vpi[\uG]$ 
over $\mm$ in the usual sense --- then by 
definition $\zp\foe\vjqo\vpi$ (via $\vt=\dom\vjpi$), 
and we are done. 
To prove the claim suppose otherwise. 
Then there is a \mut\ $\zq\in\mt\vjpi$ such that 
$\zq\leq\zp$ and $\zq$
\dd{\mt\vjpi}forces $\neg\:\vpi[\uG]$ 
over $\mm$, thus 
$\zq\foe\vjqo\neg\,\vpi$.
But this contradicts to $\zp\wfoe\vjpi\vpi$.

\ref{wf14}
There is a \muq\ $\vjqo\in\vmf$ by \ref{wf13}, such that 
$\vjpi\su\vjqo$ and $\zp\foe\vjqo\vpi$. 
Note that still $\zp\wfoe\vjqo\neg\:\vpi$ by 
Lemma~\ref{mont}. 
Extend $\vjqo$ once again, getting a contradiction 
with \ref{wf11}.
\epf 

\bcor
\lam{wf2}
Let\/ $n\ge2$, $\vjpi\in\vmf$, $\zp\in\mt\vjpi$, 
$\vpi$ be a formula in\/ $\ls1n$. 
Then\/ $\zp\foe\vjpi\vpi$ and\/ $\zp\foe\vjpi\vpi^-$ 
cannot hold together.
\ecor
\bpf
If $n=1$ then apply Lemma~\ref{wf1}\ref{wf11}. 
If $n\ge2$ then the result immediately follows 
by definition (\ref{fo4} in Section~\ref{auxB}).
\epf    

\bcor
[in $\rL$]
\lam{lspC}
Assume that\/ $\vjpi\in\vmf$, $\zp\in\mt\vjpi$, 
$\vpi$ is a formula in\/ $\lsp11$, 
$\nn\mo\zfm$ is a TM containing\/ $\vjpi$ and\/ 
$\vpi$, and\/ $\zp\wfoe\vjpi\vpi$. 
%by\/ \ref{fo21} via some\/ $\vt<\dom\vjpi$.
% and\/ $\zpo\in\mt{\vjpi\res\vt}$. 
Then\/ $\zp$ \dd{\mt\vjpi}forces\/ $\vpi[\uG]$ 
over\/ $\nn$ in the usual sense.
\ecor

This looks like the case $\vjro=\vjpi$ in 
Lemma~\ref{lsp}, but $\fof$ is weakened to 
$\wfof$, and $\vpi\in\mm$ 
(automatic in Lemma~\ref{lsp}) is added, in the premise.
 
\bpf
Otherwise there is a \mut\ $\zq\in\mt\vjpi$, 
$\zq\leq\zp$, that \dd{\mt\vjpi}forces\/ $\neg\:\vpi[\uG]$ 
over\/ $\nn$.
On the other hand, by Lemma~\ref{wf1}\ref{wf13}, 
there is a \muq\ $\vjqo\in\vmf$ such that\/ 
$\vjpi\su_{\mm(\vjpi)}\vjqo$ 
and\/ $\zp\foe\vjqo\vpi$, hence, $\zq\foe\vjqo\vpi$ by 
Lemma~\ref{mont}.
However we have $\zq\foe\vjqo\neg\,\vpi$ by 
definition 
(\ref{fo21} in Section~\ref{auxB} with $\vt=\dom\vjpi$), 
which contradicts to Lemma~\ref{wf1}\ref{wf11}.
\epf

\parf{Tail invariance}
\las{tail}

If $\vjpi=\sis{\jpi_\al}{\al<\la}\in\vmf$ and   
$\ga<\la=\dom\vjpi$ then let the 
\dd\ga\rit{tail} 
\index{tail@\dd\ga tail}%
$\vjpi\reg\ga$ be the restriction 
$\vjpi\res{[\ga,\la)}$ to the ordinal semiinterval 
$[\ga,\la)=\ens{\al}{\ga\le\al<\la}$.
Then the \muf\ 
$\mt{\vjpi\reg\ga}=\bkw_{\ga\le\al<\la}\vjpi(\al)$ 
is open dense in $\mt\vjpi$ 
by Corollary~\ref{xisc}\ref{xisc2x}.
Therefore it can be expected that if 
$\vjqo$ is another \muq\ of the same length 
$\la=\dom\vjqo$, and $\vjqo\reg\ga=\vjpi\reg\ga$, 
then the relation 
$\foe\vjpi$ coincides with $\foe\vjqo$.
And indeed this turns out to be the case (almost).

%Recall \ref{fo22} in Section~\ref{auxB} regarding $\wfof$ 
%in \ref{aut1} of the next theorem.

\bte
\lam{tat}
Assume that\/ $\vjpi,\vjqo$ are \muq s in\/ $\vmf$, 
$\ga<\la=\dom\vjpi=\dom\vjqo$, 
$\vjqo\reg\ga=\vjpi\reg\ga$, 
$\zp\in\md$, and\/ $\vpi$ is an\/ \dd\cL formula. 
Then
\ben
\renu
\itlb{tat1} 
if\/ $\vpi\in\lsp11$ then\/ 
$\zp\wfoe\vjpi\vpi$ iff\/ $\zp\wfoe{\vjqo}\vpi\;;$  

\itlb{tat2} 
if\/ $n\ge2$, $\vjpi,\vjqo\in\wmf{n-2}$, and\/ 
$\vpi\in\lp1n\cup\ls1{n+1}$, 
then\/ $\zp\foe\vjpi\vpi$ iff\/ 
$\zp\foe{\vjqo}\vpi$. 
\een
\ete
\bpf
\ref{tat1}
Suppose to the contrary that $\zp\wfoe\vjpi\vpi$, 
but $\zp\wfoe{\vjqo}\vpi$ fails, 
so there is a \muq\ 
$\vjqo'\in\vmf$ and $\zp'\in\mt{\vjqo'}$ 
such that $\vjqo\su\vjqo'$,
$\zp'\leq\zp$, and $\zp'\foe{\vjqo'}\neg\:\vpi$.
Let $\la'=\dom\vjqo'$. 
%Let $\jqo=\bkw_{\la\le\al<\la'}\vjqo'(\al)$.
By Corollary~\ref{xisc}\ref{xisc2x},
there is a \mut\ $\zr\in\mt{\vjqo'\reg\ga}$, $\zr\leq\zp'$.
Then still $\zr\leq\zp$ and $\zr\foe{\vjqo'}\neg\:\vpi$, 
by Lemma~\ref{mont}.

Define a \muq\ $\vjpi'$ so that 
$\dom{\vjpi'}=\la'=\dom{\vjqo'}$, $\vjpi\sq\vjpi'$, 
and $\vjpi'\reg\la=\vjqo'\reg\la$. 
Then $\zr\in\mt{\vjpi'}$, and $\zr\wfoe{\vjpi'}\vpi$  
by Lemma~\ref{mont}.  

Consider any CTM $\nn\mo\zfm$ containing $\vpi$,  
$\vjpi'$, $\vjqo'$. 
%hence containing  $\zr$, $\jqo$, $\mt\jqo$ as well. 
Then, by Corollary~\ref{lsp}, one and the same 
\mut\ $\zr$ \dd{\mt{\vjpi'}}forces $\vpi[\uG]$ but 
\dd{\mt{\vjqo'}}forces $\neg\:\vpi[\uG]$ over $\nn$. 
But this contradicts to the fact that the forcing 
notions $\mt{\vjpi'}$, $\mt{\vjqo'}$ contain one and 
the same dense set 
$\mt{\vjpi'\reg\la}=\mt{\vjqo'\reg\la}$.  

\ref{tat2}
Consider first the $\lp12$ case. 
Assume that $\vpi(x)$ is a $\ls11$ formula, 
$\zp\foe\vjpi\kaz x\,\vpi(x)$, but 
to the contrary $\zp\foe{\vjqo}\kaz x\,\vpi(x)$ fails.
Thus there is a \muq\ 
$\vjqo'\in\vmf$ and a \mut\ $\zp'\in\mt{\vjqo'}$ 
such that $\vjqo\sq\vjqo'$, $\zp'\leq\zp$, and 
$\zp'\foe{\vjqo'}\sus x\,\vpi^-(x)$.
By definition there is a small \qn\ $\rc$ 
such that $\zp'\foe{\vjqo'}\vpi^-(\rc)$.
There is a \mut\ $\zr\in\mt{\vjqo'\reg\ga}$, $\zr\leq\zp'$.
Then still $\zr\leq\zp$ and $\zr\foe{\vjqo'}\vpi^-(\rc)$. 
As above there is a \muq\ $\vjpi'$ such that 
$\dom{\vjpi'}=\la'=\dom{\vjqo'}$, $\vjpi\sq\vjpi'$, 
and $\vjpi'\reg\la=\vjqo'\reg\la$. 
Then $\zr\in\mt{\vjpi'}$ and $\zr\wfoe{\vjpi'}\vpi^-(\rc)$  
by the inductive hypothesis. 
By Lemma~\ref{wf1}, there is a \muq\ $\vjsg$ 
such that $\vjpi'\sq\vjsg$ and 
$\zr\wfoe{\vjsg}\vpi^-(\rc)$, 
hence, $\zr\wfoe{\vjsg}\sus x\,\vpi^-(x)$. 
But this contradicts to 
%the assumption 
$\zp\foe\vjpi\kaz x\,\vpi(x)$, since $\zr\leq\zp$.

To carry out the step $\lp1n\to\ls1{n+1}$, $n\ge2$, 
let $\vpi(x)$ be a formula in $\lp1n$. 
Assume that $\zp\foe\vjpi\sus x\,\vpi(x)$. 
By definition (see \ref{fo3} in Section~\ref{auxB}), 
there is a small \qn\ $\rc$ such that 
$\zp\foe\vjpi\vpi(\rc)$. 
Then we have $\zp\foe{\vjqo}\vpi(\rc)$ 
by the inductive assumption, thus  
$\zp\foe{\vjqo}\sus x\,\psi(x)$.  

To carry out the step $\ls1n\to\lp1{n}$, $n\ge3$, 
assume that $\vpi$ is a $\lp1{n}$ formula, 
$\zp\foe\vjpi\vpi$, but 
to the contrary $\zp\foe{\vjqo}\vpi$ fails.
Then by \ref{fo4} of Section~\ref{auxB}, 
as $\vjqo\in\wmf{n-2}$, there is a \muq\ 
$\vjqo'\in\wmf{n-2}$ and a \mut\ $\zp'\in\mt{\vjqo'}$ 
such that $\vjqo\sq\vjqo'$, $\zp'\leq\zp$, and 
$\zp'\foe{\vjqo'}\vpi^-$.
By Corollary~\ref{xisc}\ref{xisc2x},
there is a \mut\ $\zr\in\mt{\vjqo'\reg\ga}$, $\zr\leq\zp'$.
Then still $\zr\leq\zp$ and $\zr\foe{\vjqo'}\vpi^-$.
As above in the proof of \ref{tat1}, 
there is a \muq\ $\vjpi'$ such that 
$\dom{\vjpi'}=\la'=\dom{\vjqo'}$, $\vjpi\sq\vjpi'$, 
and $\vjpi'\reg\la=\vjqo'\reg\la$. 

We claim that $\vjpi'\in\wmf{n-2}$. 
Indeed if $\al<\dom\vjpi$ then 
$\vjpi'(\al)=\vjpi(\al)$ (as $\vjpi\sq\vjpi'$), 
hence 
$\prl{\vjpi'(\al)}{m-2}=\prl{\vjpi(\al)}{m-2}
=\prl{\jPi_\al}{m-2}$
(as $\vjpi$ belongs to $\wmf{n-2}$).
If $\dom\vjpi\le\al<\dom{\vjpi'}$ then 
$\vjpi'(\al)=\vjqo'(\al)$ 
(as $\vjpi'\reg\la=\vjqo'\reg\la$), 
hence 
$\prl{\vjpi'(\al)}{m-2}=\prl{\vjqo'(\al)}{m-2}
=\prl{\jPi_\al}{m-2}$
(as $\vjqo'$ belongs to $\wmf{n-2}$).
Thus $\prl{\vjpi'(\al)}{m-2}=\prl{\jPi_\al}{m-2}$
for all $\al$, meaning that 
$\prl{\vjpi'}{m-2}\su\prl{\vjPi}{m-2}$ 
and $\vjpi'\in\wmf{n-2}$. 
To conclude, 
%we have 
$\vjpi'\in\wmf{n-2}$, 
$\vjpi\sq\vjpi'$, $\zr\in\mt{\vjpi'\reg\ga}$, 
$\zr\leq\zp$, and also $\zr\foe{\vjpi'}\vpi^-$   
by the inductive hypothesis.  
But this contradicts to the assumption 
$\zp\foe\vjpi\vpi$. 
\epf

\parf{Permutations}
\las{aut}

Still {\ubf arguing in $\rL$}, 
we let $\aut$ be the set of all 
bijections $\hh:\omi\ti\om\onto\omi\ti\om$, such that 
the \rit{kernel} 
$\abs\hh=\ens{\ang{\xi,k}}{\hh(\xi,k)\ne\ang{\xi,k}}$ 
is at most countable. 
Elements of $\aut$ will be called \rit{permutations}.
\index{permutation!$\aut$}%
\index{zzperm@$\aut$}%
If $m<\om$ then let $\aut_m$ consist of those 
permutations $\hh\in\aut$ satisfying 
$\abs\hh\sq\omi\ti(\omi\bez m)$. 
In other words, if $\hh\in\aut_m$ and $\xi<\omi$, $k<m$, 
then $\hh(\xi,k)=\ang{\xi,k}$.

Let $\hh\in\aut$. 
\index{permutation!action}%
We extend the action of $\hh$  as follows.
%is naturally extended to some   objects considered above. 
\bit
\item 
if $\zp$ is a \mut\ then $\hh\zp$ is 
a \mut, 
$\abs{\hh\zp}=\ima\hh\zp=
\ens{\hh(\xi,k)}{\ang{\xi,k}\in\abs\zp}$, 
and $(\hh\zp)(\hh(\xi,k))=\zp(\xi,k)$
whenever $\ang{\xi,k}\in\abs\zp$, 
in other words, $\hh\zp$ coincides with the 
superposition $\zp\circ{(\hh\obr)}$;

\item 
if $\zpi\in\md$ is a \muf\ then 
$\hh\ap\zpi=\zpi\circ{(\hh\obr)}$ is 
a \muf, 
%satisfying 
$\abs{\hh\ap\zpi}=\ima\hh\zpi$ 
and $(\hh\ap\zpi)(\hh(\xi,k))=\zpi(\xi,k)$
whenever $\ang{\xi,k}\in\abs\zpi$;  

%\item
%if $P\sq\md$ then $\hh P=\ens{\hh\zp}{\zp\in P}$, still 
%a set of \mut s; 

\item
if $\rc\sq\md\ti(\om\ti\om)$ is a \qn, 
%so that $\kkc ni$ are subsets of $\md$, 
then put 
$\hh\rc=\ens{\ang{\hh\zp, n,i}}{\ang{\zp, n,i}\in\rc}$, 
thus easily $\hh\rc$ is a \qn\ as well;

\item
if $\vjpi=\sis{\jpi_\al}{\al<\ka}$ is a \muq, then   
$\hh\vjpi=\sis{\hh\ap\jpi_\al}{\al<\ka}$, 
still a \muq.

\item
if $\vpi:=\vpi(\rc_1,\dots,\rc_n)$ is a \dd\cL formula 
(with all names explicitly indicated), then   
$\hh\vpi$ is $\vpi(\hh\rc_1,\dots,\hh\rc_n)$.
\eit

Many notions and relations defined 
above are clearly \dd\aut invariant, \eg, 
$\zp\in\mt\jpi$ iff $\hh\zp\in\mt{\hh\ap\jpi}$, 
$\jpi\ssb{}\jqo$ iff ${\hh\ap\jpi}\ssb{}{\hh\ap\jqo}$, 
{\sl et cetera}. 
%The non-invariance can be expected in cases of reference 
%to a CTM $\mm$, like \eg\ the relation $\vjpi\su_\mm\vjqo$. 
%Fortunately t
The invariance persists even with respect to 
the relation $\fof$, at least to some extent.

\vyk{ get invariance here, one naturally has to apply 
$\hh\in\aut$ to $\mm$, too, but $\hh\mm$ 
(whatever way you define it) 
will not necessarily be a CTM, which is not good. 
The next two results pursue a more appropriate 
way to deal 
with the problem: assume that $\hh$ essentially belongs 
to $\mm$!
Note that the action of $\hh\in\aut$ is fully determined 
by the restriction $\hh\res\abc\hh$.
}

\vyk{
\ble
\lam{autl}
Assume that\/ $\mm$ is a CTM of\/ $\zfm$, $\hh\in\aut$, 
and\/ $\hh\res\abc\hh\in\mm$. 
Let\/ $\vjpi\su_\mm\vjqo$ be \muq s in\/ $\vmf$, and\/ 
$\vjpi$\/ 
{\rm(but not necessarily $\vjqo$!)} 
belongs to\/ $\mm$. 
Then\/ $\hh\vjpi,\,\hh\vjqo$ are \muq s in\/ $\vmf$ 
and\/ $(\hh\vjpi)\su_\mm(\hh\vjqo)$. 
\ele
\bpf
In fact $\vjpi\su_\mm\vjqo$ means $\vjpi\su_M\vjqo$, 
where $M$ consists of all dense sets and \qn s in $\mm$ 
involved in the definition of $\vjpi\su_\mm\vjqo$ by 
Definition~\ref{extM}. 
Therefore we have $(\hh\vjpi)\su_{\hh M}(\hh\vjqo)$
by routine observation, where $\hh M$ consists of 
all \dd\hh transforms of elements of $M$. 
Yet $\hh M\sq\mm$ since $\hh\res\abc\hh\in\mm$. 
\epf
}

%Recall \ref{fo22} in Section~\ref{auxB} regarding $\wfof$ 
%in \ref{aut1} of the next theorem.

\bte
\lam{auT}
Assume that\/ $\vjpi\in\vmf$, $\zp\in\mt\vjpi$, 
$\vpi$ 
%\in\mm(\vjpi)$ 
is an\/ \dd\cL formula, and\/ $\hh\in\aut$. 
%, $\gao<\dom\vjpi$, and\/ $\hh\res\abc\hh\in\mm(\gao)$. 
Then
\ben
\renu
\itlb{aut1} 
if\/ $\vpi$ belongs to\/ $\lsp11$ and\/ 
$\zp\foe\vjpi\vpi$, then\/ 
$(\hh\zp)\wfoe{\hh\vjpi}(\hh\vpi)\;;$  

\itlb{aut2} 
if\/ $n\ge2$, $\hh\in\aut_{n-2}$, and\/ 
$\vpi$ belongs to\/ $\lp1n\cup\ls1{n+1}$, 
then\/ $\zp\foe\vjpi\vpi$ iff\/ 
$(\hh\zp)\foe{\hh\vjpi}(\hh\vpi)$. 
\een
\ete
\bpf
\vyk{
Note first of all that 
%the \muq\ $\vjqo=\hh\vjpi$  
%satisfies $\mm(\vjqo)=\mm(\vjpi)$ 
%(since $\gao<\dom\vjpi=\dom\vjqo$ and      
%$\hh\res\abc\hh\in\mm(\gao)$), 
the \mut\ $\zq=\hh\zp$ belongs to $\mt\vjqo$, 
where $\vjqo=\hh\vjpi$, and 
$\vpi$ belongs to $\mm(\vjpi)$ iff 
the formula $\psi:=\hh\vpi$ belongs to $\mm(\vjqo)$.\vom
}
Let $\vjqo=\hh\vjpi$, $\zq=\hh\zp$, $\psi:=\hh\vpi$.

\ref{aut1}
Suppose to the contrary that $\zp\wfoe\vjpi\vpi$, 
but $\zq\wfoe{\vjqo}\psi$ fails, 
so that there is a \muq\ 
$\vjqo'\in\vmf$ and $\zq'\in\mt{\vjqo'}$ 
such that $\vjqo\su\vjqo'$,
$\zq'\leq\zq$, and $\zq'\foe{\vjqo'}\neg\:\psi$.
The \muq\ $\vjpi'=\hh\obr\vjqo'$ then satisfies 
$\vjpi\su\vjpi'$, and the \mut\ $\zp'=\hh\obr\zq'$ 
belongs to $\mt{\vjpi'}$ and $\zp'\leq\zp$, hence 
we have $\zp'\wfoe{\vjpi'}\vpi$ by Lemma~\ref{mont}.

Now let $\mm\mo\zfm$ be an arbitrary CTM containing 
$\vjpi',\vjqo',\vpi,\psi,\hh\res{\abc\hh}$. 
Then, by Corollary~\ref{lspC}, 
$\zp'$ \dd{\mt{\vjpi'}}forces $\vpi[\uG]$, but 
$\zq'$ \dd{\mt{\vjqo'}}forces $\psi[\uG]$, 
over $\mm$. 
However the sets $\mt{\vjpi'}$, $\mt{\vjqo'}$ 
belong to the same model $\mm$, 
where they are order-isomorphic 
via the isomorphism induced by $\hh\res\abc\hh$. 
Therefore, and since $\zq=\hh\zp$ and $\psi=\hh\vpi$, 
it cannot happen that both 
$\zp\,$ \dd{\mt{\vjpi'}}forces $\vpi[\uG]$ and 
$\zq\,$ \dd{\mt{\vjqo'}}forces $\neg\:\psi[\uG]$.  
But this contradicts to the above.
\vom 

\ref{aut2}
Consider first the $\lp12$ case. 
Assume that $\vpi(x)$ is a $\ls11$ formula, 
$\psi(x):=\hh\vpi(x)$, 
$\zp\foe\vjpi\kaz x\,\vpi(x)$, but 
to the contrary $\zq\foe{\vjqo}\kaz x\,\psi(x)$ fails.
Thus there is a \muq\ 
$\vjqo'\in\vmf$ and a \mut\ $\zq'\in\mt{\vjqo'}$ 
such that $\vjqo\su\vjqo'$, $\zq'\leq\zq$, and 
$\zq'\foe{\vjqo'}\sus x\,\psi^-(x)$.
By definition there is a small \qn\ $\rd$ 
such that $\zq'\foe{\vjqo'}\psi^-(\rd)$.
The \muq\ $\vjpi'=\hh\obr\vjqo'$ then satisfies 
$\vjpi\su\vjqo$, the \mut\ $\zp'=\hh\obr\zq'$ 
belongs to $\mt{\vjpi'}$ and $\zp'\leq\zp$, 
$\rc=\hh\obr{\rd}$ 
%\in\mm(\vjpi')$ 
is a small \qn, 
and we have $\zp'\wfoe{\vjpi'}\vpi^-(\rc)$ 
by \ref{aut1}.
Then by Lemma~\ref{wf1} there is a longer \muq\ 
$\vjsg\in\vmf$ satisfying $\vjpi'\su\vjsg$ and 
$\zp'\foe{\vjsg}\vpi^-(\rc)$, that is, we have 
$\zp'\foe{\vjsg}\sus x\,\vpi^-(x)$. 
But by definition (\ref{fo4} in Section~\ref{auxB}) 
this contradicts to the assumption 
$\zp\foe\vjpi\kaz x\,\vpi(x)$.\vom 

\vyk{
To prove the converse, run the same arguments for 
$\vjqo,\zq,\psi,\hh\obr$ 
instead of $\vjpi,\zp,\vpi,\hh$; 
the setup is pretty symmetric.\vom
}

To carry out the step $\lp1n\to\ls1{n+1}$, $n\ge2$, 
let $\vpi(x)$ be a formula in $\lp1n$, 
$\psi(x):=\hh\vpi(x)$, and $\hh\in\aut_{n-2}$. 
Assume that $\zp\foe\vjpi\sus x\,\vpi(x)$. 
By definition (see \ref{fo3} in Section~\ref{auxB}), 
there is a small \qn\ $\rc$ such that 
$\zp\foe\vjpi\vpi(\rc)$. 
Then we have $\zq\foe{\vjqo}\psi(\rd)$ 
by inductive assumption, where 
$\rd=\hh\rc$ is a small \qn\ itself. 
Thus  
$\zq\foe{\vjqo}\sus x\,\psi(x)$.\vom 

To carry out the step $\ls1n\to\lp1{n}$, $n\ge3$, let 
$\vpi$ be a formula in $\lp1n$, and $\hh\in\aut_{n-2}$. 
Let $\zp\foe\vjpi\vpi$, 
in particular $\vjpi\in\wmf{n-2}$, 
% and $\vpi$ belongs to $\mm(\vjpi)$ by \ref{22b},  
but, to the contrary, 
$\zq\foe{\vjqo}\psi$ fails, 
where $\zq=\hh\zp$, $\vjqo=\hh\vjpi$, 
and $\psi$ is $\hh\vpi$, as above. 
Then in our assumptions, 
$\prl{\vjqo}{m-2}=\prl{\vjpi}{m-2}$, hence  
$\vjqo\in\wmf{n-2}$ as well. 
%Note also that $\psi\in\mm(\vjqo)=\mm(\vjpi)$. 
Therefore by definition 
(\ref{fo4} in Section~\ref{auxB})  
there is a \muq\ 
$\vjqo'\in\wmf{n-2}$ and $\zq'\in\mt{\vjqo'}$ 
such that $\vjqo\sq\vjqo'$,
$\zq'\leq\zq$, and $\zq'\foe{\vjqo'}\psi^-$.
%Moreover by Lemma~\ref{stL} there is a \rit{strong} 
%\muq\ $\vjqo'$ with the properties mentioned. 

\vyk{
Now consider the action of $\hh\obr$; note that 
$\abc{\hh\obr}=\abs\hh$, 
$\hh\obr$ belongs to $\aut_{n-2}$ since 
so does $\hh$, 
and also $\hh\obr\res\abc{\hh\obr}$ belongs 
to $\mm(\gao)$ since so does $\hh\res\abc\hh$. 
We also note that

Note that  the action of $\hh\obr$ is the inverse 
of the action of $\hh$ in all cases, 
so that in particular 
$\hh\obr\zq=\zp$, $\hh\obr\vjqo=\vjpi$, and  
$\hh\obr\psi$ is $\vpi$. 
}

Now let $\zp'=\hh\obr\zq'$ and $\vjpi'=\hh\obr\vjqo'$, 
so that $\zp'\leq\zp$, $\vjpi\sq\vjpi'$, and, 
that is most important, $\vjpi'$ belongs to $\wmf{n-2}$
since so does $\vjqo'$ and $\hh\obr\in\aut_{n-2}$. 
Moreover we have $\zp'\foe{\vjpi'}\vpi^-$ 
by inductive assumption. 
We conclude that $\zp\foe\vjpi\vpi$ fails, which is a 
contradiction.  
\epf

\parf{Forcing with subsequences of the key \muq}
\las{fkm}

{\ubf We argue in $\rL$.} 
The key \muq\ 
$\vjPi=\sis{\jpn\al}{\al<\omi}\in\vmi$, satisfying 
%requirements 
\ref{keys1}, \ref{keys2}, \ref{keys3} of 
Theorem~\ref{kyt}, was fixed by 
\ref{keys}, and 
$\fP=\mt\vjPi$ is our forcing notion. 
If $\ga<\omi$ then the subsequence $\vjPi\res\ga$ 
belongs to $\wmf m$, $\kaz m$. 

\bdf
\lam{foe}
We write 
$\zp\foe\al\vpi$ instead of $\zp\foe{\vjPi\res\al}\vpi$, 
for the sake of brevity. 
Let $\zp\fof\vpi$ mean: $\zp\foe\al\vpi$ for some $\al<\omi$.
\edf

\ble
[in $\rL$]
\lam{pi}
Assume that\/ $\zp\in\fP$, $\al<\omi$, and\/ 
$\zp\foe\al\vpi$. 
Then$:$
\ben
\renu
\itlb{pi1}
if\/ $\al\le\ba<\omi$, $\zq\in\fpl\ba=\mt{\vjPi\res\ba}$, 
and\/ $\zq\leq\zp$, 
then\/ $\zq\foe\ba\vpi\;;$ 

\itlb{pi2}
if\/ $\zq\in\fP$, $\zq\leq\zp$, then\/ $\zq\foe\ba\vpi$ 
for some\/ $\ba\,;$ $\al\le\ba<\omi\;;$ 

\itlb{pi3}
if\/ $\zq\in\fP$  
and\/ $\zq\fof\vpi^-$ then\/ $\zp,\zq$ 
are somewhere almost disjoint$;$

\itlb{pi4}
therefore, 1st, if\/ $\zp,\zq\in\fP$, $\zq\leq\zp$, 
and\/ $\zp\fof\vpi$ then\/ $\zq\fof\vpi$, and 2nd,\/  
$\zp\fof\vpi$, $\zp\fof\vpi^-$ cannot hold together.
\een
\ele
\bpf
To prove \ref{pi1} apply Lemma~\ref{mont}.
To prove \ref{pi2} let $\ba$ satisfy $\al<\ba<\omi$ 
and $\zq\in\mt{\vjPi\res\ba}$, and apply \ref{pi1}.
Finally to prove \ref{pi3} note that $\zp,\zq$ 
have to be incompatible in $\fP$, as otherwise 
\ref{pi1} leads to contradiction, but the 
incompatibility in $\fP$ implies being somewhere  \ad\ 
by Corollary~\ref{sad}.
\epf

\vyk{
\ble
[density, in $\rL$]
\lam{densi}
Let\/ $\vpi$ be a closed formula in\/ 
%one of the classes\/ 
$\lsp11$ or\/ $\ls1n,\,\lp1n$, $n\ge2$, 
with small\/ \qn s as parameters. 
Then the set\/ $\zD_\vpi$ of all\/ $\zp\in\fP$,  
such that\/ $\zp\fof\vpi$ or\/ $\zp\fof\vpi^-$, 
%for some\/ $\al<\omi$,
is open dense in\/ $\fP$.
\ele
\bpf
Let $\zpo\in\fP$. 
Fix an ordinal $\al_0<\omi$ such that 
$\zp_0\in\fpl{\al_0}=\mt{\vjPi\res{\al_0}}$.\vom 

{\sl Case 1\/}: $\vpi$ is a $\lsp11$ formula. 
Consider the set $W$ of all \muq s $\vjsg\in\vmf$ 
such that $\vjPi\res{\al_0}\sq\vjsg$ and there is 
a \mut\ $\zp\in\mt\vjsg$ satisfying $\zp\leq\zp_0$ 
and $\zp\foe\vjsg\vpi$ or $\zp\foe\vjsg\vpi^-$.
It follows from Lemma~\ref{deff} that $W$ belongs to 
$\fd\hc{1}$
(with $\al_0,\,\vpi,\,\zp_0$ as parameters). 
By \ref{keys}\ref{keys2} 
there is an ordinal $\ga<\omi$ such that 
the subsequence $\vjpi=\vjPi\res\ga$ \ 
\dd{0}decides $W$.\vom

{\sl Subcase 1.1\/}: $\vjpi\in W$. 
Let this be witnessed by a \mut\ 
$\zp\in\mt{\vjpi}$, so that $\zp\leq\zpo$, 
and $\zp\foe\ga\vpi$ or $\zp\foe\ga\vpi^-$.
Thus $\zp\in\zD_\vpi$, as required.  
\vom

{\sl Subcase 1.2\/}: negative decision, that is,  
there is no \muq\ $\vjsg$ in $W$ 
extending $\vjpi=\vjPi\res\ga$.
Show that this cannot happen. 
By Corollary~\ref{xisc}\ref{xisc1}, there is a \muq\ 
$\vjqo\in\vmf$ satisfying $\vjpi\sq\vjqo$, of 
arbitrarily big $\dom\vjqo$,
in particular, big enough for the model 
$\mm=\mm(\vjqo)$ to contain 
%(all names in) 
$\vpi$. 
And again by Corollary~\ref{xisc}\ref{xisc1}, 
there is a \muq\ $\vjsg\in\vmf$ satisfying 
$\vjqo\su_\mm\vjsg$. 

Now viewing $\mt{\vjqo}$ as a forcing in $\mm$, we 
find a \mut\ $\zp\in\mt{\vjsg}$, $\zp\leq\zpo$, 
which \dd{\mt{\vjqo}}forces $\vpi[\uG]$ or 
\dd{\mt{\vjqo}}forces $\vpi^-[\uG]$
over $\mm$ in the usual sense. 
Then accordingly 
$\zp\foe\vjsg\vpi$ or $\zp\foe\vjsg\vpi^-$ 
holds by \ref{fo21} in Section \ref{auxB}.
Thus we have $\vjsg\in W$, 
contrary to the subcase assumption.
\vom

{\sl Case 2\/}: $\vpi$ is a $\ls1n$ formula, 
$n\ge2$. 
Consider the set $U$ of all 
%layer-restricted 
\muq s 
of the form $\prg\vjpi{n-2}$, where $\vjpi\in\wmf {n-2}$, 
$\dom\vjpi>\al_0$, and 
there is a \mut\ $\zp\in\mt\vjpi$ satisfying 
$\zp\leq\zpo$ ($\zp$ is stronger) and 
$\zp\foe\vjpi\vpi$.
It follows from Lemma~\ref{wmfd} and Lemma~\ref{deff}
that $U$ belongs to $\fs\hc{n-1}$ 
(with $\vpi$ and $\zpo$ as parameters). 
Therefore by \ref{keys}\ref{keys2} 
there is an ordinal $\ga<\omi$ such that 
the subsequence $\vjPi\res\ga$ \ \dd{(n-2)}decides $U$. 
%We can wlog assume that $\ga\ge\al_0$, that is, 
%in particular, $\zp_0\in \mt{\vjPi\res\ga}$.
\vom

{\sl Subcase 2.1\/}: $\prg{\vjPi\res\ga}{n-2}\in U$. 
Let this be witnessed by a \muq\ $\vjpi\in\wmf {n-2}$  
and a \mut\ $\zp\in\mt\vjpi$, so that in particular 
$\prg{\vjPi\res\ga}{n-2}=\prg{\vjpi}{n-2}$ and 
$\dom\vjpi=\ga>\al_0$. 
Then by definition (Definition~\ref{wmf}) 
we also have $\prl{\vjpi}{n-2}=\prl{\vjPi\res\ga}{n-2}$, 
so that overall $\vjpi=\vjPi\res\ga$.
Thus $\zp\in\mt{\vjPi\res\ga}$, $\zp\leq\zp_0$, and 
$\zp\foe{\vjPi\res\ga}\vpi$, that is, 
$\zp\foe{\ga}\vpi$ and $\zp\in\zD_\vpi$, as required.  
\vom

{\sl Subcase 2.2\/}: negative decision,  
no \muq\ in $U$ extends $\prg{\vjPi\res\ga}{n-2}$.
Take an ordinal $\da$, $\ga<\da<\omi$, big enough for 
the model $\mm(\vjPi\res\da)$ to contain $\vpi$. 
We claim that $\zpo\foe{\da}\vpi^-$, so that 
$\zpo$ itself belongs to $\zD$. 
Indeed otherwise by \ref{fo4} 
there is a \muq\ $\vjqo\in\wmf{n-2}$   
and a \mut\ $\zp\in\mt{\vjqo}$, such that
$\vjPi\res\da\sq\vjqo$,
$\zp\leq\zpo$, and $\zp\foe{\vjqo}\vpi$.
But then $\vjqo$ and $\zp$ witness that 
$\vjsg=\prg{\vjqo}{n-2}$ belongs to $U$.  
On the other hand, $\vjsg$ obviously 
extends $\prg{\vjPi\res\da}{n-2}$, since 
$\vjPi\res\da\sq\vjqo$, 
contrary to the subcase assumption.
\epf

\parf{Connections to the ordinary forcing} 
\las{orf}
}

Now we are going to prove that the auxiliary relation 
$\fof$ essentially 
coincides with the usual \dd\fP forcing relation 
over $\rL$.

\ble
\lam{fofo1}
If\/ $n<\om$, $\vpi$ is a closed formula as in\/ 
\ref{22b}, and\/ $\zp\in\fP$, then\/ 
$\zp$ \dd\fP forces\/ $\vpi[\uG]$ over\/ $\rL$ 
in the usual sense if and only if\/ $\zp\fof\vpi$.
\ele
\bpf
Let $\ofo$ denote the usual \dd\fP forcing relation 
over $\rL$.\vom 

{\it Part 1\/}: $\vpi$ is a formula in $\lsp11$. 
If $\zp\fof\vpi$ then $\zp\foe{\vjPi\res\ga}\vpi$
for some $\ga<\omi$, and then $\zp\ofo\vpi[\uG]$ 
by Lemma~\ref{lsp} with $\vjqo=\vjPi$ and $\nn=\rL$. 

Suppose now that $\zp\ofo\vpi[\uG]$. 
There is an ordinal $\gao<\omi$ such that 
$\zp\in\fP_{\gao}=\mt{\vjPi\res\gao}$ and $\vpi$ 
belongs to $\mm(\vjPi\res\gao)$. 
The set $U$ of all \mut s $\vjpi\in\vmf$ such that 
$\gao<\dom\vjpi$ and there is an ordinal $\vt$, 
$\gao<\vt<\dom\vjpi$, such that 
$\vjpi\res\vt\su_{\mm(\vjpi\res\vt)}\vjpi$, is 
dense in $\vjpi$ by Corollary~\ref{xisc}\ref{xisc1}, 
and is $\fd\hc1$. 
Therefore by Corollary~\ref{18a} there is an 
ordinal $\ga<\omi$ such that  
$\vjpi=\vjPi\res\ga\in U$. 

Let this be witnessed by an ordinal $\vt$, 
so that $\gao<\vt<\ga=\dom\vjpi$ and  
$\vjpi\res\vt\su_{\mm(\vjpi\res\vt)}\vjpi$. 
We claim that $\zp$ \dd{\mt{\vjpi\res\vt}}forces 
$\vpi[\uG]$ over $\mm(\vjpi\res\vt)$ in the 
usual sense --- then 
by definition $\zp\foe\vjpi\vpi$, and we are done. 

To prove the claim, suppose otherwise. 
Then there is a \mut\ $\zq\in\mt{\vjPi\res\vt}$, 
$\zq\leq\zp$, which \dd{\mt{\vjpi\res\vt}}forces 
$\neg\:\vpi[\uG]$ over $\mm(\vjpi\res\vt)$. 
Then by definition we have $\zq\foe\vjpi\neg\,\vpi$, 
hence $\zq\fof\neg\,\vpi$, which implies 
$\zq\ofo\neg\:\vpi[\uG]$ (see above), 
with a contradiction to $\zp\ofo\vpi[\uG]$.\vom 

{\it Part 2\/}:  the step $\lp1n\to\ls1{n+1}$ ($n\ge1$). 
Consider a $\lp1n$ formula $\vpi(x)$. 
Assume $\zp\fof\sus x\,\vpi(x)$. 
By definition there is a small \qn\ $\rc$ such that 
$\zp\fof\vpi(\rc)$. 
By inductive hypothesis, $\zp\ofo\vpi(c)[\uG]$, that is, 
$\zp\ofo\sus x\,\vpi(x)[\uG]$. 
Conversely, assume that $\zp\ofo\sus x\,\vpi(x)[\uG]$. 
As $\fP$ is CCC, there is a small \qn\ $\rc$ (in $\rL$) 
such that $\zp\ofo\vpi(\rc)[\uG]$. 
We have $\zp\fof\vpi(\rc)$ by the inductive hypothesis, 
hence $\zp\fof\sus x\,\vpi(x)$.\vom

{\it Part 3\/}:  the step $\ls1n\to\lp1{n}$ ($n\ge2$). 
Consider a closed $\ls1n$ formula $\vpi$. 
Assume that $\zp\fof\vpi^-$.
By Lemma~\ref{pi}\ref{pi4}, there is no \mut\ 
$\zq\in\fP$, $\zq\leq\zp$, with $\zq\fof\vpi$.
This implies $\zp\ofo\vpi^-$ by 
the inductive hypothesis. 

Conversely, suppose that $\zp\ofo\vpi^-$.
There is an ordinal $\gao<\omi$ such that 
$\zp\in\fP_{\gao}=\mt{\vjPi\res\gao}$ and $\vpi$ 
belongs to $\mm(\vjPi\res\gao)$. 
Consider the set $U$ of all 
%layer-restricted 
\muq s 
of the form $\prg\vjpi{n-2}$, where $\vjpi\in\wmf {n-2}$, 
$\dom\vjpi>\gao$, and 
there is a \mut\ $\zq\in\mt\vjpi$ satisfying 
$\zq\leq\zp$ ($\zq$ is stronger) and 
$\zq\foe\vjpi\vpi$.
It follows from Lemma~\ref{wmfd} and Lemma~\ref{deff}
that $U$ belongs to $\fs\hc{n-1}$ 
(with $\vpi$ and $\zpo$ as parameters). 
Therefore by \ref{keys}\ref{keys2} 
there is an ordinal $\ga<\omi$ such that 
the subsequence $\vjPi\res\ga$ \ \dd{(n-2)}decides $U$. 
%We can wlog assume that $\ga\ge\al_0$, that is, 
%in particular, $\zp_0\in \mt{\vjPi\res\ga}$.
%\vom

{\sl Case 1\/}: $\prg{(\vjPi\res\ga)}{n-2}\in U$. 
Let this be witnessed by a \muq\ $\vjpi\in\wmf {n-2}$  
and a \mut\ $\zq\in\mt\vjpi$, so that in particular 
$\prg{(\vjPi\res\ga)}{n-2}=\prg{\vjpi}{n-2}$ and 
$\dom\vjpi=\ga>\gao$. 
Then by definition (Definition~\ref{wmf}) 
we also have $\prl{\vjpi}{n-2}=\prl{(\vjPi\res\ga)}{n-2}$, 
so that overall $\vjpi=\vjPi\res\ga$.
Thus $\zq\in\mt{\vjPi\res\ga}$, $\zq\leq\zp$, and 
$\zq\foe{\vjPi\res\ga}\vpi$, that is, 
$\zq\ofo\vpi[\uG]$ by the inductive hypothesis, 
contrary to the choice of $\zp$.  
Therefore Case 1 cannot happen, and we have:

{\sl Case 2\/}: negative decision,  
no \muq\ in $U$ extends $\prg{(\vjPi\res\ga)}{n-2}$.
We can assume that $\ga>\gao$. 
(Otherwise replace $\ga$ by $\gao+1$.) 
We claim that $\zp\foe{\ga}\vpi^-$. 
Indeed otherwise by \ref{fo4} 
there is a \muq\ $\vjpi\in\wmf{n-2}$   
and a \mut\ $\zq\in\mt{\vjpi}$, such that
$\vjPi\res\ga\sq\vjpi$,
$\zq\leq\zp$, and $\zq\foe{\vjqo}\vpi$.
But then $\vjpi$ and $\zq$ witness that 
$\vjsg=\prg{\vjpi}{n-2}$ belongs to $U$.  
On the other hand, $\vjsg$ obviously 
extends $\prg{\vjPi\res\ga}{n-2}$, since 
$\vjPi\res\ga\sq\vjpi$, 
contrary to the Case 2 assumption.
Thus indeed $\zp\fof\vpi^-$, as required.  
\epf

The next lemma provides a useful strengthening.

\ble
\lam{fofo2}
If\/ $\Phi$ is a\/ $\id\hc1$ collection of 
closed\/ $\ip1{n+2}$ formulas, $\zpo\in\fP$, 
and\/ $\zpo$ \dd\fP forces\/ $\vpi[\uG]$ 
over\/ $\rL$ for each\/ $\vpi\in\Phi$, 
then there is an ordinal\/ $\ga<\omi$  
such that if\/ $\vpi\in\Phi$ then\/ 
$\zpo\foe{\vjPi\res\ga}\vpi$.
{\rm(Same $\ga$ for all $\vpi$.)}
\ele
\bpf 
Let $U$ consist of all \muq s of the form $\prg\vjpi{m}$, 
where $\vjpi\in\wmf {m}$, 
%$\vjPi\res\gao\su\vjpi$, 
%(and hence $\zpo\in\mt\vjpi$ by \ref{eet*}), 
and there is a formula $\vpi\in\Phi$    
and $\zp\in\mt\vjpi$ such that $\zp\leq\zpo$ and 
$\zp\foe\vjpi \vpi^-$.
It follows from Lemma~\ref{wmfd} and \ref{deff} 
that $U$ is a $\fs\hc{m+1}$ set, so by 
\ref{keys}\ref{keys2} 
there is an ordinal $\ga<\omi$ such that 
$\vjPi\res\ga$ \ \dd{m}decides $U$. 

{\sl Case 1\/}: $\prg{(\vjPi\res\ga)}m\in U$,  
that is, the \muq\ $\vjpi=\vjPi\res\ga$ satisfies the 
condition that there exist $\vpi\in\Phi$    
and a \mut\ $\zp\in\vjpi$ such that $\zp\leq\zpo$ and 
$\zp\foe\vjpi \vpi^-$, 
and hence  $\zp\,$ \dd\fP forces 
$\vpi^-[\uG]$ over $\rL$ by Lemma~\ref{fofo1}, 
contrary to the choice of $\zpo$.
Therefore Case 1 cannot happen, and we have 

{\sl Case 2\/}: no \muq\ in 
$U$ extends $\prg{(\vjPi\res\ga)}m$. 
We can assume that $\ga>\gao$. 
(Otherwise replace $\ga$ by $\gao+1$.) 
We claim that $\ga$ is as required. 
Indeed otherwise 
$\zpo\foe{\vjPi\res\ga}\vpi$ 
\rit{fails} for a formula  $\vpi\in\Phi$, 
thus (\ref{fo4} in Section~\ref{auxB}), 
there is a \muq\   
$\vjpi\in\wmf{m}$ and $\zp\in\mt{\vjpi}$ 
such that $\vjPi\res\ga\sq\vjpi$,
$\zp\leq\zpo$, 
and $\zp\foe{\vjpi}\vpi^-$.
It follows that $\prg\vjpi{m}\in U$. 
In addition, $\prg\vjpi{m}$ extends 
$\prg{(\vjPi\res\ga)}m$  
by construction. 
But this contradicts to the Case 2 assumption.
\epf

\gla{The model}
\las{them}

In this conclusive section we gather the results 
obtained above towards the proof of Theorem~\ref{mt}. 
We begin with the analysis of definability of key 
generic reals in \dd\fP generic extensions 
of $\rL$, which will lead to \ref{mt1} and \ref{mt2} 
of Theorem~\ref{mt}. 
Then we proceed to \ref{mt3} (elementary equivalence) 
and \ref{mt4} (the non-wellorderability).

\parf{Key generic extension and subextensions} 
\las{kge}

Recall that the key \muq\ $\vjPi=\sis{\jPi_\al}{\al<\omi}$ 
of small \muf s $\jPi_\al$ is defined in $\rL$ by \ref{keys},   
$\jPi=\bkw_{\al<\omi}$ is a \muf, $\abc\jPi=\omi\ti\om$ 
in $\rL$, and $\fP=\mt\vjPi=\mt\jPi\in\rL$ is our key 
forcing notion, equal to the finite-support product 
$\prod_{\xi<\omi,k<\om}\jPi(\xi,k)$ of perfect-tree 
forcings $\jPi(\xi,k)$ in $\rL$. 
See Section~\ref{bpf}, where some properties of $\fP$ are 
established, including CCC and definability of the 
factors $\jPi(\xi,k)$.

From now on, we'll typically 
argue in $\rL$ and in \dd\fP generic extensions 
of $\rL$, so by Lemma~\ref{ccc} it will always be true that 
$\omil=\omi$.
This allows us to still think that $\abc\jPi=\omi\ti\om$ 
(rather than $\omil\ti\om$).

Recall that $\jPi\in\rL$ and $\fP=\mt\jPi$ is a forcing 
notion in $\rL$.

\bdf
\lam{bstr}
Let a set $G\sq\fP$ be generic over the 
constructible set universe $\rL$.
If $\ang{\xi,k}\in\omi\ti\om$ 
then following Remark~\ref{adds}, we  
\bit
\item[$-$]
define 
$G(\xi,k)=\ens{\zd \zp\xi k}{\zp\in G\land\ang{\xi,k}\in\abc\zp}
\sq  \jPi(\xi, k)$;

\item[$-$]
let $x_{\xi k}=x_{\xi k}[G]\in\dn$ be the only real in 
$\bigcap_{T\in G(\xi,k)}[T]$.
\eit
Thus $\fP$ adds an array   
$X=\sis{x_{\xi k}}{\ang{\xi,k}\in\omi\ti\om}$ of reals, 
where each real 
$x_{\xi k}=x_{\xi k}[G]\in\in\dn\cap\rL[G]$ 
is a \dd{\jPi(\xi, k)}generic real over $\rL$, 
and $\rL[G]=\rL[X]$.  
\edf

Let $G\sq\fP$ be a set (filter) \dd\fP generic over $\rL$. 
If $m<\om$ then following the notation in 
Section~\ref{bl} we define
$$
\prl G m= G\cap \mtl m=\ens{\prl \zp m}{\zp\in G}\,,
\kmar{prl G m}
$$
so that the set $\prl G m$ is \dd{\prl\fP m}generic 
\kmar{prl fP m}
over $\rL$, where accordingly 
$$
\prl\fP m= \fP\cap \mtl m=\ens{\prl \zp m}{\zp\in\fP}\,.
$$
Each subextension $\rL[\prl G m]\sq\rL[G]$ coincides 
with $\rL[\sis{x_{\xi k}[G]}{\xi<\omi\land k<m}]$.
Our goal will be to demonstrate that the model 
$\rL[X]=\rL[G]$, along with the system of submodels 
$\rL[\sis{x_{\xi k}[G]}{\xi<\omi\land k<m}]$, 
proves Theorem~\ref{mt}.

\parf{Definability of generic reals} 
\las{dgr}
%{bex1}

Recall that the factors $\jPi(\xi,k)$ of the forcing 
notion $\jPi$ are defined by 
$\jPi(\xi,k)=\bigcup_{\al(\xi,k)\le\al<\omi}\jpn{\al}(\xi,k)$, 
where $\al(\xi,k)<\omi$, the sets $\jpn{\al}(\xi,k)$ are 
countable sets of perfect trees, whose  definability in 
$\rL$ is determined by Corollary~\ref{doml2}.

\bte
\lam{xik}
Assume that a set\/ $G\sq\fP$ is\/ 
\dd\fP generic over\/ $\rL$, 
$\xi<\omi$, $k<\om$, and\/ $x\in\rL[G]\cap\bn.$ 
The following are equivalent$:$
\ben
\nenu
\itlb{xik1}
$x=x_{\xi k}[G]\;;$

\itlb{xik2}
$x$ is\/ \dd{\jPi(\xi,k)}generic over\/ $\rL\;;$

\itlb{xik3}
$x\in\bigcap_{\al(\xi,k)\le\al<\omi}
\bigcup_{T\in\jPi_\al(\xi,k)}[T]$.
\een
\ete
\bpf
$\ref{xik1}\imp\ref{xik2}$ is a routine (see Remark~\ref{adds}). 
To check $\ref{xik2}\imp\ref{xik3}$ recall that each set 
$\jPi_\al(\xi,k)$ is pre-dense in $\jPi(\xi,k)$ by 
Lemma~\ref{pqr}\ref{pqr4}.
It remains to establish $\ref{xik3}\imp\ref{xik1}$.
Suppose towards the contrary that a real $x\in\rL[G]\cap\dn$ 
satisfies \ref{xik3} but $x\ne x_{\xi k}[G]$. 
By Theorem~\ref{npn}\ref{npn1} there is a 
%dyadic 
true \rn\jPi\ $\rc=\sis{\kc ni}{n,i<\om}$, 
non-principal over $\jPi$ at $\xi,k$  
and such that\/ $x=\rc[G]$. 
Being non-principal means that the next set
is open dense in\/ $\fP=\mt\jPi$:
$$
\ddi\xi k\rc\jPi
=\ens{\zp\in\fP=\mt\jPi}{\zp\,\text{ directly forces }\,
\rc\nin[\zd\zp\xi k]}\,.
$$
And as $\fP=\mt\jPi$ is a CCC forcing by Lemma~\ref{ccc}, 
we can assume that the name $\rc$ is small, 
that is, each set $\kc ni\sq\fP$ is countable.
Then there is an ordinal $\ga_0<\omi$ such that 
$\kc ni\sq\fpl{\ga_0}$ for all $n,i$.
Then $\rc$ is a true 
%\dd{\fpl{\ga_0}}real name. 
\rn{\pilg{\ga_0}}.
Moreover we can assume by Corollary~\ref{Cjden} 
that $\ddi\xi k\rc\jPi\cap\fpl{\ga_0}$ is 
pre-dense in\/ $\fP$. 

Now consider the set $W$ of all \muq s 
$\vjpi=\sis{\jpi_\al}{\al<\len\vjpi}\in\vmf$ 
such that $\len{\vjpi}>\ga_0$ and 
\bit
\item[$-$] 
either (I) $\vjPi\res\ga_0\not\su\vjpi$;

\item[$-$] 
or (II) $\vjPi\res\ga_0\su\vjpi$ and $\rc$  
is {\ubf not} non-principal over 
$\jpi=\bkw\vjpi$ at $\xi,k$;

\item[$-$] 
or (III) $\vjPi\res\ga_0\su\vjpi$, 
$\len{\vjpi}=\da+1$ is 
a successor, and 
$\bkw_{\al<\da}\jpi_\al\sse\rc\xi k\jpi_\da$.
\eit
We assert that $W$ is \rit{dense} in $\vmf$: 
any \muq\ $\vjpi\in\vmf$ can be extended to some 
$\vjqo\in W$. 
Indeed first extend $\vjpi$ by Corollary~\ref{xisc} 
so that is has a length 
$\len{\vjpi}=\da>\ga_0$. 
If now $\vjPi\res\ga_0\not\su\vjpi$ then immediately 
$\vjpi\in W$ via (I), so we assume that 
$\vjPi\res\ga_0\su\vjpi$. 
We can also assume that $\rc$  
is non-principal over $\jpi=\bkw\vjpi$ at $\xi,k$ 
by similar reasons related to (II).
%Under these assumptions, t
The \muq\ $\vjpi$ can 
be extended, by Corollary~\ref{xisc}, by an extra term 
$\jpi_\da$, so that the extended \muq\ $\vjpi_+$ 
satisfies $\vjpi\su_{\ans{\rc}}\vjpi_+$, that is, 
$\jpi\ssm{\ans{\rc}}\jpi_\da$. 
By definition (Definition~\ref{extM}) and the 
nonprincipality of $\rc$, we get 
$\jpi\sse\rc\xi k\jpi_\da$. 
It follows that  $\vjpi_+\in W$ via (III).  

Since $W$ is $\fs\hc1$, by Definition~\ref{keys}\ref{keys2} 
there is an ordinal $\ga<\omi$ such that the \muq\   
$\vjPi\res\ga$ \dd0decides $W$.
However the negative decision is impossible by the 
density (see the proof of Lemma~\ref{doml}). 
We conclude that $\vjPi\res\ga\in W$; hence, $\ga>\ga_0$.
Option (I) for $\vjpi=\vjPi\res\ga$ clearly fails, and 
(II) fails either because the set
$\ddi\xi k\rc\jPi\cap\fpl{\ga_0}$ is 
pre-dense in $\fP$ and $\ga>\ga_0$. 
Therefore $\vjPi\res\ga$ belongs to $W$ via (III), 
that is, $\ga=\da+1$ and 
$\pilg\da=\bkw_{\al<\da}\jPi_\al\sse\rc\xi k\jPi_\da$.
Then 
$\pilg\da\sse\rc\xi k\pigg\da=\bkw_{\da\le\al<\omi}\jPi_\da$ 
by Lemma~\ref{pqo}\ref{pqo2}. 

Now we make use of Theorem~\ref{npn}\ref{npn2} with 
$\jpi=\pilg\da$ and $\jqo=\pigg\da$; note that 
$\jpi\kw\jqo=\jPi$. 
It follows that 
$x=\rc[G]\nin\bigcup_{Q\in\pigg\da(\xi,k)}[Q]$, 
which clearly contradicts to the assumption \ref{xik3}.
\epf  

\bcor
\lam{mod2}
Assume that\/ $k<\om$ and\/ $G\sq\fP$ is\/ \dd\fP generic 
over\/ $\rL$. 
Then 
$$
\gfu_k=\ens{\ang{\xi,x_{\xi k}[G]}}{\xi<\omi}\sq\omi\ti\dn
$$
is a set of definability class\/ $\ip\hc{k+2}$ in\/ $\rL[G]$ 
and in any transitive model\/ $M\mo\zfc$ satisfying\/ 
$\rL\sq M\sq\rL[G]$ and\/ $\ens{x_{\xi k}[G]}{\xi<\omi}\sq M$. 
\ecor
\bpf
By the theorem, it is true in $\rL[G]$ 
that $\ang{\xi,x}\in \gfu_k$ iff 
$$
\kaz \al<\omi\:\sus T\in\jPi_\al(\xi,k)
\big(\al(\xi,k)\le\al \imp x\in[T]\big),
$$
which can be re-written as 
$$
\kaz \al<\omi\:\kaz \mu<\omi\:\kaz X\:
\sus T\in X
\big(\mu=\al(\xi,k)\land X=\jPi_\al(\xi,k)\land\mu\le\al
\imp x\in[T]\big).
$$
Here the equality $\mu=\al(\xi,k)$ (with a fixed $k$) 
is $\id\hc{k+2}$ by Corollary~\ref{doml2}, and so is 
the equality $X=\jPi_\al(\xi,k)$ by 
Corollary~\ref{doml2}. 
It follows that the whole relation is $\ip\hc{k+2}$, 
since the quantifier $\sus T\in X$ is bounded.
\epf

%\parf{Subextensions and wellorderings} 
%\las{swo}

The next corollary is the first cornerstone in the proof 
of Theorem~\ref{mt}.

\bcor
[= \ref{mt1}, \ref{mt2} of Theorem~\ref{mt}]
\lam{mod3}
Assume that\/ $m<\om$ and a set\/ $G\sq\fP$ is\/ 
\dd\fP generic over\/ $\rL$. 
Then\/ $\bn\cap\rL[\prl G m]$ is a\/ $\is1{m+3}$ 
set in $\rL[G]$,  
and it holds in\/ $\rL[\prl G m]$  that there is 
a\/ $\id1{m+3}$ wellordering of\/ $\bn$ of length\/ $\omi$. 
\ecor
\bpf
If $\ga<\omi$ then let 
$X_{\ga n}=\sis{x_{\xi k}[G]}{\xi<\ga\land k<n}$; 
thus $X=X_{\ga n}$ is a $\ip\hc{n+1}$ relation 
in $\rL[G]$ (with $\ga,n,X$ as arguments)
by Corollary~\ref{mod2}. 
It follows that 
$$
\bn\cap\rL[\prl G m] = 
\ens{x\in\bn}{\sus\ga<\omi(x\in\rL[X_{\ga n}])}
$$
is a set in $\is\hc{n+2}$, hence, a $\is1{m+3}$ 
set in $\rL[G]$.
To define a required wellordering, if 
$x\in\bn\cap\rL[\prl G m]$ then let $\ga(x)$ be 
the least $\ga<\omi$ such that $x\in\rL[X_{\ga n}]$, 
and let $\nu(x)<\omi$ be the index of $x$ in the 
canonical wellordering of $\bn$ in $\rL[X_{\ga n}]$.
Now we wellorder $\bn\cap\rL[\prl G m]$ according to 
the lexicographical ordering of triples 
$\ang{\tmax\ans{\ga(x),\nu(x)},\ga(x),\nu(x)}$.
\epf

\vyk{
Further conclusive steps will be to prove that
\ben
\nenu
\item 
each $\rL[\prl G m]$ is an elementary submodel of 
$\rL[G]$ with respect to $\is1{m+2}$ formulas, and 
\item 
$\rL[G]$ does not admit an analytically definable 
wellordering of the reals.
\een
}

\parf{Elementary equivalence}
\las{ee}

Here we prove the following elementary equivalence 
theorem for key generic extensions.
Compare to \ref{mt3} of Theorem~\ref{mt}.

\bte
\lam{eet}
Assume that\/ $m<\om$ and a set\/ $G\sq\fP$ is\/ 
\dd\fP generic over\/ $\rL$. 
Then\/ $\rL[\prl G m]$ is an elementary submodel of\/ 
$\rL[G]$ \poo\ all\/ $\is1{m+2}$ formulas. 
\ete
\bpf
Suppose that this is not the case. 
Then there is a $\ip1{m+1}$ formula $\vpi(r,x)$ with 
$r\in\bn\cap\rL[\prl G m]$ as the only parameter, 
and a real $x_0\in\bn\cap\rL[G]$
such that $\vpi(r,x_0)$ is true in $\rL[G]$ 
but there is no $x\in\bn\cap\rL[\prl G m]$ 
such that $\vpi(r,x)$ is true in $\rL[G]$. 
By a version of Proposition~\ref{r2n}\ref{r2n3}, 
we have $r=\rco[G]$, where $\rco$ is a small 
true \rn{(\prl\fP m)}. 
(See Section~\ref{kge} on notation.)
And there is a small true \rn{\fP} $\rc$ such that 
$x_0=\rc[\uG]$. 
%Thus $\rc\sq\fP\ti(\om\ti\om)$ and 
%$\rco\sq(\prl\fP m)\ti(\om\ti\om)$ are small names,   
%and $\vpi(r,\rc[\uG])$ is true in $\rL[G]$. 

By Lemma~\ref{fofo1}, there is a \mut\ $\zpo\in G$ 
such that 
\ben
\nenu
\itlb{po1}
$\zpo$ \ \dd\fP forces  
$
\vpi(\rc_0[\uG],\rc[\uG])
\,\land\,
\neg\:\sus x\in \rL[\prl \uG m]\,\vpi(\rc_0[\uG],x) 
$
over $\rL$; 

\itlb{eet*}\msur
$\zpo\fof \vpi(\rc_0,\rc)$, that is, 
$\zpo\foe{\vjPi\res\gao}\vpi(\rc_0,\rc)$, 
where $\gao<\omi$ --- 
and we can assume that $\zpo\in\mt{\vjPi\res\gao}$ 
%and $\rc,\rco\in\mm(\vjPi\res\gao)$ 
as well.   
\een
\vyk{
We can assume in addition that 
\ben
\nenu
\atc
\atc
\itlb{eet*1}\msur
$\rco$ is a small true \rn{(\prl{\fpl\gao}m)},
$\rc$ is a small true \rn{({\fpl\gao})}.
%and $\vjPi\res\gao$ is strong. 
\een
By definition this means that 
$\zpo\in\mt{\pilg{\gao}}$, 
$\rco$ is a small true \rn{(\prl{\pilg{\gao}}m)}
while $\rc$ is a small true \rn{({\pilg{\gao}})},  
where $\pilg{\gao}=\bkw_{\xi<\gao}\jPi_\xi$ is 
a small \muf. 
}%
As $\rc,\rco$ are small names, there is an 
ordinal $\da<\omi$ satisfying 
%$\abc{\pilg{\gao}}\sq\da\ti\om$. 
%It follows then from 
%\ref{eet*1} that 
\ben
\atc
\atc
\nenu
\itlb{eet*2}\msur
$\abc\rco\sq \da\ti m$, $\abc\rc\sq\da\ti\om$, 
 and $\abc\zpo\sq\da\ti\om$, 
\een
and we can enlarge $\gao$, if necessary, 
using the equality $\abc\vjPi=\omi\ti\om$ 
of Lemma~\ref{doml}, 
to make sure that 
\ben
\atc
\atc
\atc
\nenu
\itlb{eet*3}\msur
$\da\ti\om\sq\abc{\vjPi\res\gao}$, 
that is, if $\et<\da$ and $k<\om$ then 
$\ang{\et,k}\in\abc{\jPi_{\al'}}$ for some 
$\al'=\al'(\et,k)<\gao$.
\een
We are starting from here towards a contradiction.  

Let $U$ consist of all \muq s of the form $\prg\vjpi{m}$, 
where 
\ben
\Aenu
%\atc\atc\atc\atc
\itlb{bab0}
$\vjpi\in\wmf {m}$, 
$\vjPi\res\gao\su\vjpi$, 
and hence $\zpo\in\mt\vjpi$ by \ref{eet*};
\een
and 
there is an ordinal $\za<\omi$ and 
a transformation $\hh\in \aut_{m-1}$ such that 
\ben
\Aenu
%\atc\atc\atc\atc
\atc
\itlb{bab2} 
$\hh=\hh\obr$, $\abc\hh=D\cup R$, 
and $\hh$ maps $D$ onto $R$ and $R$ onto $D$, 
where $D=\da\ti[m,\om)$, 
$R=\ens{\ang{\xi,m-1}}{\nuo\le\xi<\nui}$, and  
$\da<\nuo<\nui<\omi$;  
%and $\abc\zp\sq\gad\ti\om$;

%\itlb{bab3}
%there is an ordinal $\gai$, $\gao\le\gai<\dom\vjpi$, 
%such that $\hh\res\abc\hh\in\mm(\gai)$;

\itlb{bab1}\msur
%$\hh\vjpi=\vjpi$;
%
%there is an ordinal $\za$, 
$\gao\le\za<\dom\vjpi$ and  
$(\hh\vjpi)\reg\za=\vjpi\reg\za$, 
or equivalently $\hh(\vjpi(\al))=\vjpi(\al)$ whenever 
$\za\le\al<\dom\vjpi$.
\een
It follows from Lemma~\ref{wmfd} that $U$ is 
a $\fs\hc{m+1}$ set 
(with $\vjPi\res\gao$, $\da$ as parameters). 
Therefore by \ref{keys}\ref{keys2} 
there is an ordinal $\ga<\omi$ such that 
$\vjPi\res\ga$ \ \dd{m}decides $U$.\vom

{\sl Case 1\/}: $\prg{(\vjPi\res\ga)}m\in U$. 
Basically this means that there is a transformation 
$\hh\in\aut_{m-1}$ such that 
\ref{bab0}, \ref{bab2}, 
%\ref{bab3}, 
\ref{bab1} 
hold for $\hh$ and $\vjpi=\vjPi\res\ga$, via  
ordinals $\da<\nuo<\nui$ and $\gao<\za<\ga$ as in 
\ref{bab2}, 
%\ref{bab3}, 
\ref{bab1}.

Now, by Lemma~\ref{mont} and \ref{eet*}, we have 
$\zpo\foe{\vjPi\res\ga}\vpi(\rco,\rc)$. 
We further get 
$\hh\zpo\foe{\hh\vjPi\res\ga}\vpi(\hh\rco,\hh\rc)$
by Theorem~\ref{auT}
because $\vpi$ is a $\lp1{n+1}$ 
formula and $\hh$ belongs to $\aut_{m-1}$. 
% and satisfies \ref{bab3}. 
%
However $\hh\rco=\rco$ since $\abc\rco\cap\abc\hh=\pu$ 
by \ref{bab2}.
Thus  
$\zpo'\foe{\vjPi\res\ga}\vpi(\rco,\rc')$ 
holds by Theorem~\ref{tat} and \ref{bab1}, 
where $\rc'=\hh\rc$, $\zpo'=\hh\zpo$. 

Note that the common part $\abc\zpo\cap\abc{\zpo'}$ 
of the domains of $\zpo,\zpo'$ does not intersect  
$\abc\hh$ by \ref{bab2} since $\abc{\zpo}\sq\da\ti\om$ 
by \ref{eet*2}. 
It follows that $\zpo,\zpo'$ are compatible, basically 
$\zp=\zpo\cup\zpo'$ is a \mut\ in $\mt{\vjPi\res\ga}$. 
Thus $\zp\leq\zpo'$ and still 
$\zp\foe{\vjPi\res\ga}\vpi(\rco,\rc')$.
It follows by Lemma~\ref{fofo1} that 
\ben
\nenu
\atc\atc\atc
\atc
\itlb{999}
$\zp$ \ \dd\fP forces  
$
\vpi(\rco[\uG],\rc'[\uG])
$
over $\rL$.
\een 
However $\abc{\rc'}\sq\omi\ti m$ by construction 
because $\abc\rc\sq\da\ti\om$ by \ref{eet*2}, 
and hence $\rc'[\uG]\in \rL[\prl \uG m]$ is forced. 
Thus $\zp$ \ \dd\fP forces  
$\sus x\in \rL[\prl \uG m]\,\vpi(\rco[\uG],x)$
over $\rL$ by \ref{999}, contrary to \ref{po1}. 
The contradiction ends Case 1.\vom

{\sl Case 2\/}: negative decision, no \muq\ in 
$U$ extends $\prg{(\vjPi\res\ga)}m$. 
We can assume that $\ga>\gao$. 
(Otherwise replace $\ga$ by $\gao+1$.) 
Let $\nuo$ be the lest ordinal, bigger than $\da$ 
and satisfying $\abc{\vjPi\res\ga}\sq\nuo\ti\om$. 
Let $\nui=\nuo+\om$. 
Then countable sets $D=\da\ti[m,\iy)$ and $R$ 
as in \ref{bab2} are 
defined and $D\cap R=\pu$, so we can fix a 
transformation $\hh\in\aut_{m-1}$ satisfying \ref{bab2}. 
Note that $D\sq\da\ti\om\sq\abc{\vjPi\res\ga}$ by 
\ref{eet*3} but $R\cap\abc{\vjPi\res\ga}=\pu$ by 
the choice of $\nuo$.

%Pick an ordinal $\gai$, $\ga<\gai<\omi$, to satisfy 
%$\hh\res\abc\hh\in\mm(\gai)$ in \ref{bab3}.
Pick $\la<\omi$ such that $\la>\ga>\gao$. 
Then the \muq\ $\vjqo=\vjPi\res\la$ clearly satisfies 
\ref{bab0}, \ref{bab2} 
%, \ref{bab3} 
and extends $\vjPi\res\ga$. 
Our plan is now to slightly modify $\vjqo$ 
in order to fulfill \ref{bab1} as well, with $\za=\ga$.
Such a minor modification consists in the replacement 
of the \dd Rpart of $\vjqo$ above $\ga$ 
by the \dd\hh copy of its \dd Dpart.

To present this in detail, recall that 
$\vjqo=\vjPi\res\la=\sis{\jPi_\al}{\al<\la}$, where 
each $\jPi_\al$ is a small \muf, whose domain 
$d_\al=\abc{\jPi_\al}\sq\omi\ti\om$ is countable. 
If $\al<\ga$ then  put $\jpi_\al=\jPi_\al$. 
Suppose that $\ga\le\al<\la$. 
Then $D\sq\abc{\jPi_\al}$ by \ref{eet*3}.
%Let $D_\al=d_\al\cap D$, $R_\al=d_\al\cap R$.
Define a modified \muf\ $\jpi_\al$ such that
\ben
\aenu
%\item[$-$]\msur
\itlb{jpia1}\msur
$\abc{\jpi_\al}=d_\al\cup R$ --- note that 
$D\sq d_\al\sq\abc{\jpi_\al}$ in this case because 
$D\sq\abc{\vjPi\res\ga}$ by \ref{eet*3} 
(as $\gao\le\ga$), and hence 
$D\sq d_\al=\abc{\jPi_\al}$ (as $\al\ge\ga$),  

\itlb{jpia2}
if $\ang{\xi,k}\in d_\al\bez R$ then  
$\jpi_\al(\xi,k)=\jPi_\al(\xi,k)$, 

\itlb{jpia3}
if $\ang{\xi,k}\in D$, so  
%(that is, $\xi<\da$ and $k\ge m$) and 
$\hh(\xi,k)=\ang{\et,m-1}\in R$, then  
$\jpi_\al(\et,m-1)=\jPi_\al(\xi,k)$. 
\een
We claim that $\vjpi=\sis{\jpi_\al}{\al<\la}$ is a 
\muq, that is, if $\al<\ba<\la$ then 
$\jpi_\al\bssq\jpi_\ba$. 
This amounts to the folowing: if 
$\ang{\et,k}\in \abc{\jpi_\al}$ then 
$\jpi_\al(\et,k)\ssq\jpi_\ba(\et,k)$. 
Note that $\jpi_\al(\et,k)=\jPi_\al(\et,k)$ in case 
$\ang{\et,k}\nin R$. 

Thus it remains to check that 
$\jpi_\al(\et,m-1)\ssq\jpi_\ba(\et,m-1)$ whenever 
$\al<\ba<\la$, 
$\ang{\et,m-1}=\hh(\xi,k)\in R\cap \abc{\jpi_\al}$, and 
$\ang{\xi,k}\in D$.
If now $\al<\ga$ then $R\cap\abc{\jpi_\al}=\pu$ 
by the choice of $\nuo$, so it remains to consider  
the case when $\ga\le\al$. 
Then the pairs $\ang{\xi,k}$, $\ang{\et,m-1}$ belong to 
$\abc{\jpi_\al}$ by construction, and we have 
$\jpi_\al(\et,m-1)=\jPi_\al(\xi,k)$ and 
$\jpi_\ba(\et,m-1)=\jPi_\ba(\xi,k)$.  
Therefore $\jpi_\al(\xi,m)\ssq\jpi_\ba(\xi,m)$ 
since $\vjPi$ is a \muq, and we are done.

Now we claim that the \muq\ $\vjpi=\sis{\jpi_\al}{\al<\la}$ 
satisfies \ref{bab0}, \ref{bab2}, 
%\ref{bab3}, 
\ref{bab1}. 
Indeed as the difference between each 
$\jpi_\al$ and the corresponding $\jPi_\al$  
is fully located in the domain 
$R=\ens{\ang{\xi,m-1}}{\nuo\le\xi<\nui}$, 
we have $\prl {\vjpi}{m-1}=\prl {\vjqo}{m-1}$, 
therefore $\vjpi\in\wmf {m}$. 
We also note that $\vjpi\res\ga=\vjqo\res\ga$ by 
construction, hence $\vjPi\res\ga=\vjqo\res\ga\su\vjpi$. 
This implies \ref{bab0}. 

We also have \ref{bab2} 
%, \ref{bab3} 
by construction. 
We finally claim that \ref{bab1} is satisfied with  
$\za=\ga$, that is, if $\ga\le\al<\la$ then 
$\hh\ap\jpi_\al=\jpi_\al$. 
Indeed we have $D\cup R\sq\abc{\jpi_\al}$, see 
\ref{jpia1}. 
Now the invariance of $\jpi_\al$ under $\hh$ holds 
by \ref{jpia2}, \ref{jpia3}.

It follows that $\prg\vjpi{m}\in U$. 
In addition, $\prg\vjpi{m}$ extends 
$\prg{(\vjPi\res\ga)}m$, since 
$\vjPi\res\ga\su\vjpi$. 
But this contradicts the Case 2 assumption.\vom

To conclude, either case leads to a contradiction, 
proving the theorem. 
\epf

\parf{Non-wellorderability}
\las{nwo}

We finally prove that the reals are not wellorderable by 
a (lightface) analytically definable relation in 
\dd\fP generic extensions, that is, \ref{mt4} of 
Theorem~\ref{mt}.

\bte
\lam{nwot}
Assume that\/ $m<\om$ and a set\/ $G\sq\fP$ is\/ 
\dd\fP generic over\/ $\rL$. 
Then it is true in\/ $\rL[G]$ that the reals\/ 
are not wellorderable by an analytically definable 
relation. 
\ete
\bpf
Suppose to the contrary that, in $\rL[G]$, a $\is1{m+2}$ 
relation $\ll$ strictly wellorders $\bn$, $m\ge1$.
Let $\psi(x,y)$ be a parameter-free $\is1{m+2}$ formula, 
which defines $\ll$, so that $x \ll y$ iff $\psi(x,y)$ 
in $\rL[G]$. 
Note that $\ll$ is essentially a $\id1{m+2}$ relation, 
since $x\ll y\leqv y\not \ll x\land x\ne y$.

Of all nonconstructible reals $x_{\xi m}[G]$, $\xi<\omi$, 
there is a \dd\ll least one. 
We suppose that $x_{0 m}[G]$ is such. 
(If it is some $x_{\xi_0 m}[G]$, $\xi_0\ne0$, 
then the arguments suitably change in obvious way.) 
That is, $x_{0 m}[G]\ll x_{\xi m}[G]$
whenever $\xi>0$. 
Accordingly there is a \mut\ $\zpo\in G$ that  
\dd\fP forces, over $\rL$, that\vom 

(1) 
$\ll$ (that is, the relation defined by $\psi$) 
is a wellordering of $\bn$, and\vom
 
(2) 
$\kaz\xi>0\,(x_{0 m}[\uG]\ll x_{\xi m}[\uG])$.\vom

\noi 
Therefore, if $\xi>0$ then $\zpo$ \dd\fP forces 
$\rpi_{0m}[\uG]\ll \rpi_{\xi m}[\uG])$ 
over $\rL$.  
(We make use of the \qn s $\rpi_{\xi k}$ introduced 
by \ref{proj1}, \ref{proj2}.) 

By Lemma~\ref{fofo1}, we can assume that 
$\zpo\fof(1)\land (2)$, so that in fact we have 
$\zpo\foe{\vjPi\res\gao}(1)\land (2)$, 
for some $\gao<\omi$.
Then $\zpo\in\mt{\vjPi\res\gao}=\mt{\pilg{\gao}}$, 
%$\rc$ is a small true \rn{({\pilg{\gao}})},  
where $\pilg{\gao}=\bkw_{\xi<\gao}\jPi_\xi$ is 
a small \muf. 
Let $\da<\omi$ be the least ordinal satisfying 
%$\gai>\gao$ and 
$\abc{\pilg{\gao}}\sq\da\ti\om$. 
It follows then that $\abc\zpo\sq\da\ti\om$.
%Note that (1) is a $\ip1{m+3}$ statement, but (2) is more 
%complex, of course. 

\vyk{

\ble
\lam{nwol1}
There is an ordinal\/ $\ga$, $\gao<\ga<\omi$, 
such that if\/ $\xi<\omi$ then\/ 
$\zpo\foe{\vjPi\res\ga}(\rpi_{\xi m}\ll \rpi_{0 m})^-$.
\ele
\bpf[Lemma]
Let $U$ consist of all \muq s of the form $\prg\vjpi{m}$, 
where $\vjpi\in\wmf {m}$, 
$\vjPi\res\gao\su\vjpi$, 
%(and hence $\zpo\in\mt\vjpi$ by \ref{eet*}), 
and there is an ordinal $\xi<\dom\vjpi$,   
and a \mut\ $\zp\in\vjpi$ such that $\zp\leq\zpo$ and 
$\zp\foe\vjpi (\rpi_{\xi m}\ll \rpi_{0 m})$.

It follows from Lemma~\ref{wmfd} and \ref{deff} 
that $U$ is a $\fs\hc{m+1}$ set.  
Therefore by \ref{keys}\ref{keys2} 
there is an ordinal $\ga<\omi$ such that 
$\vjPi\res\ga$ \ \dd{m}decides $U$.\vom

{\sl Case 1\/}: $\prg{(\vjPi\res\ga)}m\in U$. 
Basically this means that $\ga>\gao$ 
and the \muq\ $\vjpi=\vjPi\res\ga$ satisfies the 
condition that there exist $\xi<\dom\vjpi$,   
and a \mut\ $\zp\in\vjpi$ such that $\zp\leq\zpo$ and 
$\zp\foe\vjpi (\rpi_{\xi m}\ll \rpi_{0 m})$, 
and hence $\zp\fof (\rpi_{\xi m}\ll \rpi_{0 m})$ 
by definition, therefore $\zp$ \dd\fP forces 
$\rpi_{\xi m}[\uG]\ll \rpi_{0 m}[\uG]$, contrary to 
the choice of $\zpo$.
Therefore Case 1 cannot happen, and we have\vom

{\sl Case 2\/}: negative decision, no \muq\ in 
$U$ extends $\prg{(\vjPi\res\ga)}m$. 
We can assume that $\ga>\gao$. 
(Otherwise replace $\ga$ by $\gao+1$.) 
We claim that $\ga$ is as required. 
Indeed otherwise 
$\zpo\foe{\vjPi\res\ga}(\rpi_{\xi m}\ll \rpi_{0 m})^-$ 
\rit{fails} for an ordinal $\xi<\omi$, 
thus (\ref{fo4} in Section~\ref{auxB}), 
there is a \muq\   
$\vjpi\in\wmf{m}$ and \mut\ $\zp\in\mt{\vjpi}$ 
such that $\vjPi\res\ga\sq\vjpi$,
$\zp\leq\zpo$, 
and $\zp\foe{\vjpi}(\rpi_{\xi m}\ll \rpi_{0 m})$.
It follows that $\prg\vjpi{m}\in U$. 
In addition, $\prg\vjpi{m}$ extends 
$\prg{(\vjPi\res\ga)}m$  
by construction. 
But this contradicts the Case 2 assumption.
\epF{Lemma}
}

%To continue the proof of the theorem, fix 

By Lemma~\ref{fofo2}, there is an ordinal 
$\gai$, $\gao<\gai<\omi$, such that  
if $\xi<\omi$ then 
$\zpo\foe{\vjPi\res\gai}(\rpi_{\xi m}\ll \rpi_{0 m})^-$.
We can enlarge $\gai$, if necessary, using Lemma~\ref{doml}, 
to make sure that $\ang{0,m}\in\abc{\vjPi\res\gai}$, 
that is, $\ang{0,m}\in\abc{\jPi_{\al'}}$ for some 
$\al'<\gai$.
%
%The following argument resembles the proof of 
%Theorem~\ref{eet} to some extent.

If $\xi<\omi$ then let $\hh_\xi\in\aut_m$ be the 
permutation of $\ang{0,m}$ and $\ang{\xi,m}$, such that 
$\abc{\hh_\xi}=\ans{\ang{0,m},\ang{\xi,m}}$, 
$\hh_\xi(0,m)=\ang{\xi,m}$, 
$\hh_\xi(\xi,m)=\ang{0,m}$, 
$\hh_\xi(\et,n)=\ang{\et,n}$ for any pair $\ang{\et,n}$ 
different from both $\ang{0,m}$ and $\ang{\xi,m}$.

The remainder of the proof is very similar to the proof 
of Theorem~\ref{eet}. 
Let $U$ consist of all \muq s of the form $\prg\vjpi{m}$, 
where 
\ben
\Aenu
%\atc\atc\atc\atc
\itlb{nwo1}
$\vjpi\in\wmf {m}$, 
$\vjPi\res\gai\su\vjpi$, 
and hence $\zpo\in\mt\vjpi$ by \ref{eet*};
\een
and there exist ordinals $\xi,\za<\omi$ such that 
\ben
\Aenu
\atc
\itlb{nwo2}\msur
$\da<\xi<\dom\vjpi$; 
% --- 
%then by the way $\hh_\xi\res\abc{\hh_\xi}\in\mm(\vjpi)$;

\itlb{nwo3}\msur
%$\hh_\xi\vjpi=\vjpi$ --- oops this is not going to happen, 
%so we require that 
%there exists an ordinal $\za$ with 
$\gai\le\za<\dom\vjpi$ and 
$(\hh_\xi\vjpi)\reg\za=\vjpi\reg\za$, 
or equivalently $\hh_\xi\ap(\vjpi(\al))=\vjpi(\al)$ whenever 
$\za\le\al<\dom\vjpi$.
\een
It follows from Lemma~\ref{wmfd} that $U$ is 
a $\fs\hc{m+1}$ set 
(with $\vjPi\res\gai$ as a parameter). 
Therefore by \ref{keys}\ref{keys2} 
there is an ordinal $\ga<\omi$ such that 
$\vjPi\res\ga$ \ \dd{m}decides $U$.\vom

{\sl Case 1\/}: $\prg{(\vjPi\res\ga)}m\in U$. 
Basically this means that $\ga>\gai$ 
and there are ordinals $\xi<\omi$ and $\za<\ga$ 
such that \ref{nwo1}, \ref{nwo2}, \ref{nwo3} 
hold for $\xi$ and the \muq\ $\vjpi=\vjPi\res\ga$.
By Lemma~\ref{mont} and the choice of $\gai$,  
$\zpo\foe{\vjPi\res\ga}(\rpi_{\xi m}\ll \rpi_{0 m})^-$ 
holds. 
This implies  
$\hh_\xi\zpo\foe{\hh_\xi\vjPi\res\ga}
(\rpi_{0 m}\ll \rpi_{\xi m})^-$
by Theorem~\ref{auT} because 
$(\rpi_{0 m}\ll \rpi_{\xi m})^-$ is a $\lp1{m+2}$ 
formula and $\hh_\xi$ belongs to $\aut_{m}$. 
%
%However $\hh_\xi(\vjPi\res\ga)=\vjPi\res\ga$ by \ref{nwo2}. 
%while $\ba\rco=\rco$ since $\abc\rco\cap\abc\ba=\pu$ by \ref{bab2}.
We conclude that  
$\zpo'\foe{\vjPi\res\ga}(\rpi_{0 m}\ll \rpi_{\xi m})^-$ 
by \ref{nwo3} and Theorem~\ref{tat}, 
where $\zpo'=\hh_\xi\zpo$. 

Note that $\ang{\xi,m}\nin \abc{\zpo}$ by \ref{nwo2}. 
It follows by the definition of $\hh_\xi$ 
that $\zpo,\zpo'$ are compatible, basically 
$\zp=\zpo\cup\zpo'$ is a \mut\ in $\mt{\vjPi\res\ga}$ 
and $\zp\leq\zpo$, $\zp\leq\zpo'$. 
Thus both  
$\zp\foe{\vjPi\res\ga}(\rpi_{0 m}\ll\rpi_{\xi m})^-$ and 
$\zp\foe{\vjPi\res\ga}(\rpi_{\xi m}\ll\rpi_{0 m})^-$ 
hold, 
contrary to Lemma~\ref{pi}\ref{pi4}. 
The contradiction ends Case 1.\vom

{\sl Case 2\/}: negative decision, no \muq\ in 
$U$ extends $\prg{(\vjPi\res\ga)}m$. 
We can assume that $\ga>\gai$  
(otherwise replace $\ga$ by $\gai+1$). 
Note that 
$\vjPi\res\ga=\sis{\jPi_\al}{\al<\ga}$, where 
each $\jPi_\al$ is a small \muf, whose domain 
$d_\al=\abc{\jPi_\al}\sq\omi\ti\om$ is countable. 
Therefore the set 
$d=\abc{\vjPi\res\ga}=\bigcup_{\al<\ga}d_\al$ is 
still countable. 
Pick an ordinal $\xi$, $\da<\xi<\omi$, 
such that $\ang{\xi,m}\nin d$. 
Finally pick an ordinal $\la$, $\ga<\la<\omi$. 
%such that $\xi<\la$ and $\ang{\xi,m}\in\abc{\vjqo}$, 
%where $\vjqo=\vjPi\res\la$.
Then $\vjqo=\vjPi\res\la$ (as $\vjpi$) and $\xi$ 
clearly satisfy \ref{nwo1} and \ref{nwo2}, 
and $\vjqo$ extends $\vjPi\res\ga$. 
Let's somewhat modify $\vjqo$ 
in order to fulfill \ref{nwo3} as well.

As above,   
$\vjqo=\vjPi\res\la=\sis{\jPi_\al}{\al<\la}$,  
each $\jPi_\al$ is a small \muf, and its domain 
$d_\al=\abc{\jPi_\al}\sq\omi\ti\om$ is countable. 
If $\al<\ga$ then  put $\jpi_\al=\jPi_\al$. 
Suppose that $\ga\le\al<\la$. 
Then $\al\ge\gai$, and hence 
$\ang{0,m}\in \abc{\jPi_\al}$ by the choice of $\gai$.
Define a modified \muf\ $\jpi_\al$ such that
$\abc{\jpi_\al}=d_\al\cup\ans{\ang{\xi,m}}$,
if $\ang{\et,k}\in d_\al\bez\ans{\ang{\xi,m}}$ 
then  $\jpi_\al(\et,k)=\jPi_\al(\et,k)$, and finally
$\jpi_\al(\xi,m)=\jPi_\al(0,m)$. 

We claim that $\vjpi=\sis{\jpi_\al}{\al<\la}$ is a 
\muq, that is, if $\al<\ba<\la$ then 
$\jpi_\al\bssq\jpi_\ba$. 
This amounts to the folowing: if 
$\ang{\et,k}\in \abc{\jpi_\al}$ then 
$\jpi_\al(\et,k)\ssq\jpi_\ba(\et,k)$. 
Note that $\jpi_\al(\et,k)=\jPi_\al(\et,k)$ whenever 
$\ang{\et,k}\ne\ang{\xi,m}$. 
Thus it remains to check that 
$\jpi_\al(\xi,m)\ssq\jpi_\ba(\xi,m)$ given 
$\al<\ba<\la$ such that $\ang{\xi,m}\in \abc{\jpi_\al}$. 
If $\al<\ga$ then $\abc{\jpi_\al}=\abc{\jPi_\al}=d_\al$ 
by construction, and hence $\ang{\xi,m}\nin d_\al$ by 
the choice of $\xi$. 
It remains to consider the case $\ga\le\al<\la$.
Then $\ang{0,m}\in d_\al$ (see above), hence the 
pairs $\ang{0,m}$, $\ang{\xi,m}$ belong to 
$\abc{\jpi_\al}$ by construction, 
and then obviously belong to $\abc{\jpi_\ba}$ as $\al<\ba$. 
Now $\jpi_\al(\xi,m)=\jPi_\al(0,m)$ and 
$\jpi_\ba(\xi,m)=\jPi_\ba(0,m)$, and we have  
$\jpi_\al(\xi,m)\ssq\jpi_\ba(\xi,m)$ since $\vjPi$ 
is a \muq.

Now we claim that the \muq\ 
$\vjpi=\sis{\jpi_\al}{\al<\la}$ 
satisfies \ref{nwo1}, \ref{nwo2}, \ref{nwo3} with $\za=\ga$. 
%belongs to $U$.
%
If $\al<\la$ then the difference between 
$\jPi_\al$ and $\jpi_\al$
is located in the one-element domain 
$\ans{\ang{\xi,m}}$, therefore 
$\prl {\jpi_\al}{m}=\prl {\jPi_\al}{m}$. 
It follows that  
$\prl {\vjpi}{m}=\prl {(\vjPi\res\la)}{m}$, hence 
$\vjpi\in\wmf {m}$. 
We further have  
$\vjpi\res\ga=\vjPi\res\ga$ by construction. 
Thus $\vjPi\res\ga\su\vjpi$, hence 
$\vjPi\res\gai\su\vjpi$, and we have \ref{nwo1}. 

We also have \ref{nwo2} and \ref{nwo3} 
(with $\za=\ga$) by construction. 

Thus $\prg\vjpi{m}\in U$. 
In addition, $\prg\vjpi{m}$ extends 
$\prg{(\vjPi\res\ga)}m$, since even more  
$\vjPi\res\ga\su\vjpi$ by construction. 
But this contradicts to the Case 2 assumption.\vom

To conclude, either case leads to a contradiction, 
proving the theorem.
\epf

\parf{Proof of the main theorem}
\las{pmt}

\bpf[Theorem~\ref{mt}]
We consider a \dd\fP generic extension $\rL[G]$ of 
$\rL$ and present it in the form $\rL[G]=\rL[X]$ 
as in Section~\ref{kge}, where 
$X=\sis{x_{\xi k}}{\ang{\xi,k}\in\omi\ti\om}$, and each  
$x_{\xi k}=x_{\xi k}[G]$ is a real in $\dn\cap\rL[G]$.
We also consider the subextensions 
$\rL[\prl G m]=\rL[\sis{x_{\xi k}}{\xi<\omi\land k<m}]$ 
of $\rL[G]=\rL[X]$. 
Then \ref{mt1} and \ref{mt2} of Theorem~\ref{mt} hold by 
Corollary~\ref{mod3}, 
\ref{mt3} holds by Theorem~\ref{eet}, and finally 
\ref{mt4} holds by Theorem~\ref{nwot}. 
\epf

\np

%\clearpage

\bibliographystyle{plain}
{\small
%\bibliography{u}

%
}

\np

\normalsize
\printindex

\end{document}